\documentclass[11pt]{epiarticle}
\usepackage{epimath}
\usepackage{cancel}

\setpapertype{A4}

\usepackage[english]{babel}

\usepackage[dvips]{graphicx}     

\newtheorem{definition}{Definition}

\newtheorem{anomaly}{Anomaly}

\title{A density-based approach for non-heuristic approximations of prime counting functions}

\titlemark{A density-based approach for non-heuristic approximations of prime counting functions}
\author{Bhupinder Singh Anand}
\authormark{B.\ S.\ Anand}
\date{Update of \today} 
\journal{Preprint at arXiv: \href{http://arxiv.org/abs/1510.04225}{http://arxiv.org/abs/1510.04225}} 

\begin{document}

\maketitle


\begin{prelims}

\def\abstractname{Abstract}
\abstract{All the known approximations of $\pi(n)$ for finite values of $n$ are derived from real-valued functions that are asymptotic to $\pi(x)$, such as $\frac{x}{log_{e}x}$, $Li(x)$ and Riemann's function $R(x) = \sum_{n=1}^{\infty}\frac{\mu(n)}{(n)}li(x^{1/n})$. The degree of approximation for finite values of $n$ is determined only heuristically, by conjecturing upon an error term in the asymptotic relation that can be seen to yield a closer approximation than others to the actual values of $\pi(n)$. By considering the density of the set of all integers that are not divisible by the first $\pi(\sqrt {n})$ primes $p_{_{1}}, p_{_{2}}, \ldots, p_{_{\pi(\sqrt {n})}}$ we show that, for any $n$, the expected number of such integers in any interval of length $(p_{_{\pi(\sqrt{ n})+1}}^{2} - p_{_{\pi(\sqrt n)}}^{2})$ is $(p_{_{\pi(\sqrt{ n})+1}}^{2} - p_{_{\pi(\sqrt n)}}^{2})\prod_{i = 1}^{\pi(\sqrt{n})}(1 - \frac{1}{p_{_{i}}})$. We then show that a non-heuristic approximation---with a bounded cumulative error term $> 0$---for the number of primes less than or equal to $n$ is given for all $n$ by $\pi(n) \approx \sum_{j = 1}^{n}\prod_{i = 1}^{\pi(\sqrt{j})}(1 - \frac{1}{p_{_{i}}}) \sim a.\frac{n}{log_{e}n} \rightarrow \infty$ for some constant $a > 2.e^{-\gamma} \approx 1.12292 \ldots$. We further show that the expected number of Dirichlect and twin primes in the interval ($p_{_{\pi(\sqrt {n})}}^{2},\ p_{_{\pi(\sqrt{ n})+1}}^{2}$) can be estimated similarly; and conclude that the number of such primes $\leq n$ is, in each case, cumulatively approximated non-heuristically by a function that $\rightarrow \infty$.}

\keywords{Chebyshev's Theorem; complete system of incongruent residues; computational complexity; Dirichlect primes; Euler's constant $\gamma$; expected value; factorising is polynomial time; integer factorising algorithm; Mertens' theorem; mutually independent prime divisors; polynomial time algorithm; prime counting function; prime density; primes in an arithmetic progression; Prime Number Theorem; probability model; twin primes.}


\tableofcontents

\end{prelims}


\section{Eratosthenes sieve and the nature of divisibility}
\label{intro}

\begin{quote}
\footnotesize
``Prime numbers are the most basic objects in mathematics. They also are among the most mysterious, for after centuries of study, the structure of the set of prime numbers is still not well understood. Describing the distribution of primes is at the heart of much mathematics...".\footnote{Andrew Granville: from \href{http://www.ams.org/gnews\#!news_id=499}{this} AMS press release of 5 December 1997.}
\end{quote}

\noindent In this investigation we show how the usual, linearly displayed, Eratosthenes sieve argument reveals the structure of divisibility (and, ipso facto, of primality) more transparently when displayed as a 2-dimensional matrix representation of the residues $r_{i}(n)$, defined for all $n \geq 2$ and all $i \geq 2$ by:

\vspace{+1ex}
$n + r_{i}(n) \equiv 0\ (mod\ i)$, where $i > r_{i}(n) \geq 0$ \footnote{See \S \ref{sec:4.fig1}, Appendix II(A), Fig.7 and II(B), Fig.8.}.

\vspace{+1ex}
\noindent \normalsize \textbf{Density:} For instance, the residues $r_{i}(n)$ can be defined for all $n \geq 1$ as the values of the density-defining functions $R_{i}(n)$, defined for all $i \geq 1$ as illustrated below\footnote{For $R_{_{i}}$ and $r_{_{i}}$ read $R_{_{i}}(n)$ and $r_{_{i}}(n)$ respectively. See also Fig.7.}, where:  

\tiny
\vspace{+.5ex}
\begin{tabbing}
Divisorsxxxx: \= 00001 \= 00002 \= 00003 \= 00004 \= 00005 \= 00006 \= 00007 \= 00008 \= 00009 \= 000010 \= 000011 \= \ldots xxx \= 0n-9 \ldots \kill
Function: \> $R_{_{1}}$ \> $R_{_{2}}$ \> $R_{_{3}}$ \> $R_{_{4}}$ \> $R_{_{5}}$ \> $R_{_{6}}$ \> $R_{_{7}}$ \> $R_{_{8}}$ \> $R_{_{9}}$ \> $R_{_{10}}$ \> $R_{_{11}}$ \> \ldots \> $R_{_{n}}$ \\
\rule{113mm}{.5mm} \\
$\textcolor{black}{n = 1}$ \> \textcolor{black}{0} \> 1 \> 2 \> 3 \> 4 \> 5 \> 6 \> 7 \> 8 \> 9 \> 10 \> \ldots \> n-1 \\
$\textcolor{black}{n = 2}$ \> \textcolor{black}{0} \> \textcolor{black}{0} \> 1 \> 2 \> 3 \> 4 \> 5 \> 6 \> 7 \> 8 \> 9 \> \ldots \> n-2 \\
$\textcolor{black}{n = 3}$ \> \textcolor{black}{0} \> \textcolor{black}{1} \> \textcolor{black}{0} \> 1 \> 2 \> 3 \> 4 \> 5 \> 6 \> 7 \> 8 \> \ldots \> n-3 \\
$\textcolor{black}{n = 4}$ \> \textcolor{black}{0} \> 0 \> 2 \> \textcolor{black}{0} \> 1 \> 2 \> 3 \> 4 \> 5 \> 6 \> 7 \> \ldots \> n-4 \\
$\textcolor{black}{n = 5}$ \> \textcolor{black}{0} \> \textcolor{black}{1} \> \textcolor{black}{1} \> \textcolor{black}{3} \> \textcolor{black}{0} \> 1 \> 2 \> 3 \> 4 \> 5 \> 6 \> \ldots \> n-5 \\
$\textcolor{black}{n = 6}$ \> \textcolor{black}{0} \> 0 \> 0 \> 2 \> 4 \> \textcolor{black}{0} \> 1 \> 2 \> 3 \> 4 \> 5 \> \ldots \> n-6 \\
$\textcolor{black}{n = 7}$ \> \textcolor{black}{0} \> \textcolor{black}{1} \> \textcolor{black}{2} \> \textcolor{black}{1} \> \textcolor{black}{3} \> \textcolor{black}{5} \> \textcolor{black}{0} \> 1 \> 2 \> 3 \> 4 \> \ldots \> n-7 \\
$\textcolor{black}{n = 8}$ \> \textcolor{black}{0} \> 0 \> 1 \> 0 \> 2 \> 4 \> 6 \> \textcolor{black}{0} \> 1 \> 2 \> 3 \> \ldots \> n-8 \\
$\textcolor{black}{n = 9}$ \> \textcolor{black}{0} \> 1 \> 0 \> 3 \> 1 \> 3 \> 5 \> 7 \> \textcolor{black}{0} \> 1 \> 2 \> \ldots \> n-9 \\
$\textcolor{black}{n = 10}$ \> \textcolor{black}{0} \> 0 \> 2 \> 2 \> 0 \> 2 \> 4 \> 6 \> 8 \> \textcolor{black}{0} \> 1 \> \ldots \> n-10 \\
$\textcolor{black}{n = 11}$ \> \textcolor{black}{0} \> \textcolor{black}{1} \> \textcolor{black}{1} \> \textcolor{black}{1} \> \textcolor{black}{4} \> \textcolor{black}{1} \> \textcolor{black}{3} \> \textcolor{black}{5} \> \textcolor{black}{7} \> \textcolor{black}{9} \> \textcolor{black}{0} \> \ldots \> n-11 \\
\\
$\textcolor{black}{n}$ \> $\textcolor{black}{r_{_{1}}}$ \> $r_{_{2}}$ \> $r_{_{3}}$ \> $r_{_{4}}$ \> $r_{_{5}}$ \> $r_{_{6}}$ \> $r_{_{7}}$ \> $r_{_{8}}$ \> $r_{_{9}}$ \> $r_{_{10}}$ \> $r_{_{11}}$ \> \ldots \> \textcolor{black}{0} \\
\rule{113mm}{.5mm} \\ \end{tabbing}

\footnotesize
\noindent $\bullet$ Each function $R_{_{i}}(n)$ cycles through the values $(i-1,\ i-2,\ \ldots,\ 0)$ with period $i$;

\vspace{+1ex}
\noindent $\bullet$ For any $i \geq 2$ the density---over the set of natural numbers---of the set $\{n\}$ of integers that are divisible by $i$ is $\frac{1}{i}$; and the density of integers that are not divisible by $i$ is $\frac{i - 1}{i}$.

\vspace{+1ex}
\noindent \normalsize \textbf{Primality:} The residues $r_{i}(n)$ can also be viewed alternatively as values of the associated primality-defining sequences, $E(n) = \{r_{i}(n): i \geq 1\}$, defined for all $n \geq 1$, as illustrated below\footnote{See Fig.8.}, where:

\tiny
\vspace{+.5ex}
\begin{tabbing}
Divisorsxxxx: \= 00001 \= 00002 \= 00003 \= 00004 \= 00005 \= 00006 \= 00007 \= 00008 \= 00009 \= 000010 \= 000011 \= \ldots xxx \= 0n-9 \ldots \kill

Function: \> $R_{_{1}}$ \> $R_{_{2}}$ \> $R_{_{3}}$ \> $R_{_{4}}$ \> $R_{_{5}}$ \> $R_{_{6}}$ \> $R_{_{7}}$ \> $R_{_{8}}$ \> $R_{_{9}}$ \> $R_{_{10}}$ \> $R_{_{11}}$ \> \ldots \> $R_{_{n}}$ \\
\rule{113mm}{.5mm} \\
$\textcolor{cyan}{E(1)}$: \> \textcolor{cyan}{0} \> 1 \> 2 \> 3 \> 4 \> 5 \> 6 \> 7 \> 8 \> 9 \> 10 \> \ldots \> n-1 \\
$\textcolor{red}{E(2)}$: \> \textcolor{cyan}{0} \> \textcolor{cyan}{0} \> 1 \> 2 \> 3 \> 4 \> 5 \> 6 \> 7 \> 8 \> 9 \> \ldots \> n-2 \\
$\textcolor{red}{E(3)}$: \> \textcolor{cyan}{0} \> \textcolor{red}{1} \> \textcolor{cyan}{0} \> 1 \> 2 \> 3 \> 4 \> 5 \> 6 \> 7 \> 8 \> \ldots \> n-3 \\
$\textcolor{cyan}{E(4)}$: \> \textcolor{cyan}{0} \> \textcolor{cyan}{0} \> 2 \> \textcolor{cyan}{0} \> 1 \> 2 \> 3 \> 4 \> 5 \> 6 \> 7 \> \ldots \> n-4 \\
$\textcolor{red}{E(5)}$: \> \textcolor{cyan}{0} \> \textcolor{red}{1} \> \textcolor{red}{1} \> \textcolor{red}{3} \> \textcolor{cyan}{0} \> 1 \> 2 \> 3 \> 4 \> 5 \> 6 \> \ldots \> n-5 \\
$\textcolor{cyan}{E(6)}$: \> \textcolor{cyan}{0} \> \textcolor{cyan}{0} \> \textcolor{cyan}{0} \> 2 \> 4 \> \textcolor{cyan}{0} \> 1 \> 2 \> 3 \> 4 \> 5 \> \ldots \> n-6 \\
$\textcolor{red}{E(7)}$: \> \textcolor{cyan}{0} \> \textcolor{red}{1} \> \textcolor{red}{2} \> \textcolor{red}{1} \> \textcolor{red}{3} \> \textcolor{red}{5} \> \textcolor{cyan}{0} \> 1 \> 2 \> 3 \> 4 \> \ldots \> n-7 \\
$\textcolor{cyan}{E(8)}$: \> \textcolor{cyan}{0} \> \textcolor{cyan}{0} \> 1 \> \textcolor{cyan}{0} \> 2 \> 4 \> 6 \> \textcolor{cyan}{0} \> 1 \> 2 \> 3 \> \ldots \> n-8 \\
$\textcolor{cyan}{E(9)}$: \> \textcolor{cyan}{0} \> 1 \> \textcolor{cyan}{0} \> 3 \> 1 \> 3 \> 5 \> 7 \> \textcolor{cyan}{0} \> 1 \> 2 \> \ldots \> n-9 \\
$\textcolor{cyan}{E(10)}$: \> \textcolor{cyan}{0} \> \textcolor{cyan}{0} \> 2 \> 2 \> \textcolor{cyan}{0} \> 2 \> 4 \> 6 \> 8 \> \textcolor{cyan}{0} \> 1 \> \ldots \> n-10 \\
$\textcolor{red}{E(11)}$: \> \textcolor{cyan}{0} \> \textcolor{red}{1} \> \textcolor{red}{1} \> \textcolor{red}{1} \> \textcolor{red}{4} \> \textcolor{red}{1} \> \textcolor{red}{3} \> \textcolor{red}{5} \> \textcolor{red}{7} \> \textcolor{red}{9} \> \textcolor{cyan}{0} \> \ldots \> n-11 \\
\ldots \\
$E(n)$: \> $\textcolor{cyan}{r_{_{1}}}$ \> $r_{_{2}}$ \> $r_{_{3}}$ \> $r_{_{4}}$ \> $r_{_{5}}$ \> $r_{_{6}}$ \> $r_{_{7}}$ \> $r_{_{8}}$ \> $r_{_{9}}$ \> $r_{_{10}}$ \> $r_{_{11}}$ \> \ldots \> \textcolor{cyan}{0} \\ \rule{113mm}{.5mm} \\
\end{tabbing}

\footnotesize
\noindent $\bullet$ The sequences $E(n)$ highlighted in red correspond to a prime\footnote{Conventionally defined as integers that are not divisible by any smaller integer other than $1$.} $p$ (since $r_{i}(p) \neq 0$ for $1 < i < p$) in the usual, linearly displayed, Eratosthenes sieve:

\vspace{+1ex}
$\textcolor{cyan}{E(\cancel{1})},\ \textcolor{red}{E(2)},\ \textcolor{red}{E(3)},\ \textcolor{cyan}{E(\cancel{4})},\ \textcolor{red}{E(5)},\ \textcolor{cyan}{E(\cancel{6})},\ \textcolor{red}{E(7)},\ \textcolor{cyan}{E(\cancel{8})},\ \textcolor{cyan}{E(\cancel{9})},\ \textcolor{cyan}{E(\cancel{10})},\ \textcolor{red}{E(11)},\ \ldots$

\vspace{+1ex}

\noindent $\bullet$ The sequences highlighted in cyan identify a crossed out composite $n$ (since $r_{i}(n) = 0$ for some $i < i < n$) in the usual, linearly displayed, Eratosthenes sieve. 

\normalsize
\vspace{+2ex}
\noindent By considering the density of the set of all integers that are not divisible by the first $k$ primes $p_{_{1}}, p_{_{2}}, \ldots, p_{_{k}}$, we shall show that the expected number of such integers in any interval of length $(p_{_{\pi(\sqrt{ n})+1}}^{2} - p_{_{\pi(\sqrt n)}}^{2})$ is $\{(p_{_{\pi(\sqrt{ n})+1}}^{2} - p_{_{\pi(\sqrt n)}}^{2})\prod_{i = 1}^{k}(1 - \frac{1}{p_{_{i}}})\}$. 

\vspace{+1ex}
\noindent We shall conclude non-heuristically that:

\begin{quote}
$\bullet$ For each $n$, the expected number of primes in the interval $(p_{_{\pi(\sqrt {n})}}^{2},\ p_{_{\pi(\sqrt{ n})+1}}^{2})$ is $\{(p_{_{\pi(\sqrt{ n})+1}}^{2} - p_{_{\pi(\sqrt {n})}}^{2})\prod_{i = 1}^{\pi(\sqrt {n})}(1 - \frac{1}{p_{_{i}}})\}$. The number $\pi(n)$ of primes $\leq n$ is thus cumulatively approximated (Lemma \ref{sec:2.5.lem.1} and Corollary \ref{sec:2.2.lem.7.3.cor.1}) for $n \geq 4$ by $\pi(n) \approx \sum_{j = 1}^{n}\prod_{i = 1}^{\pi(\sqrt{j})}(1 - \frac{1}{p_{_{i}}}) \sim a.\frac{n}{log_{e}n} \rightarrow \infty$.

\vspace{+1ex}
$\bullet$ For each $n$, the expected number of Dirichlect primes---of the form $a+m.d$ for some natural number $m \geq 1$---in the interval $(p_{_{\pi(\sqrt {n})}}^{2},\ p_{_{\pi(\sqrt{ n})+1}}^{2})$  is $\{(p_{_{\pi(\sqrt{ n})+1}}^{2} - p_{_{\pi(\sqrt {n})}}^{2})\prod_{i = 1}^{k}{\frac{1}{q_{_{i}}^{\alpha_{_{i}}}}}.\prod_{i = 1}^{k}(1 - \frac{1}{q_{_{i}}})^{-1}.\prod_{j = 1}^{\pi(\sqrt{n})}(1 - \frac{1}{p_{_{j}}})\}$, where $1 \leq a < d = q_{_{1}}^{\alpha_{_{1}}}.q_{_{2}}^{\alpha_{_{2}}} \ldots q_{_{k}}^{\alpha_{_{k}}}$ and $(a, d) = 1$. The number $\pi_{_{(a, d)}}(n)$ of Dirichlect primes $\leq n$ is thus cumulatively approximated (Lemma \ref{dirichlect.lem}) for all $n \geq q_{_{k}}^{2}$ by $\pi_{_{(a, d)}}(n) \approx \prod_{i = 1}^{k}{\frac{1}{q_{_{i}}^{\alpha_{_{i}}}}}.\prod_{i = 1}^{k}(1 - \frac{1}{q_{_{i}}})^{-1}.\sum_{l = 1}^{n}\prod_{j = 1}^{\pi(\sqrt{l})}(1 - \frac{1}{p_{_{j}}}) \rightarrow \infty$.

\vspace{+1ex}
$\bullet$ For each $n$, the expected number of $\mathbb{TW}$ primes---such that $n$ is a prime and $n+2$ is either a prime or $p_{_{\pi(\sqrt{n})+1}}^{2}$---in the interval $(p_{_{\pi(\sqrt {n})}}^{2},\ p_{_{\pi(\sqrt{ n})+1}}^{2})$ is $\{(p_{_{\pi(\sqrt{ n})+1}}^{2} - p_{_{\pi(\sqrt {n})}}^{2})\prod_{i=2}^{\pi(\sqrt{n})}(1 - \frac{2}{p_{_{i}}})\}$. The number $\pi_{_{2}}(p_{_{k+1}}^{2})$ of twin primes $\leq p_{_{k+1}}^{2}$ is thus cumulatively approximated (Lemma \ref{sec:twin.lem.4.1.1}) for all $k \geq 1$ by $\pi_{_{2}}(p_{_{k+1}}^{2}) \approx \sum_{j=9}^{p_{_{k+1}}^{2}} \prod_{i=2}^{\pi(\sqrt{j})-1}(1 - \frac{2}{p_{_{i}}}) \rightarrow \infty$.
\end{quote}

\subsection{The functions $\pi(x)$ and $\frac{x}{log_{e}x}$: A historical perspective}
\label{intro.1}

To place this investigation in an appropriate historical perspective, we note that Adrien-Marie Legendre and Carl Friedrich Gauss are reported\footnote{\scriptsize{cf. Prime Number Theorem. (2014, June 10). In Wikipedia, The Free Encyclopedia. Retrieved 09:53, July 9, 2014, from \href{http://en.wikipedia.org/w/index.php?title=Prime_number_theorem&oldid=612391868}{http://en.wikipedia.org/w/index.php?title\=Prime\_number\_theorem\&oldid=612391868}}; \footnotesize{see also \cite{Gr95}.}} to have independently conjectured in 1796 that, if $\pi (x)$ denotes the number of primes less than $x$, then $\pi (x)$ is asymptotically equivalent to $\frac{x}{log_{e}x}$. 

\vspace{+1ex}
\noindent Around 1848/1850, Pafnuty Lvovich Chebyshev proved that $\pi(x) \asymp \frac{x}{log_{e}x}$, and confirmed that if $\pi (x)/\frac{x}{log_{e}x}$ has a limit, then it must be 1\footnote{\cite{Di52}, p.439; see also \cite{HW60}, p.9, Theorem 7 and p.345, \S22.4 for a proof of Chebychev's Theorem.}. 

\vspace{+1ex}
\begin{center}
\line(1,0){260} \\
\textbf{Fig.1: The asymptotic behaviour of the primes} \\
\line(1,0){260}

\includegraphics{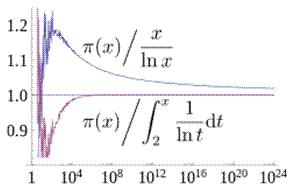}
\end{center}

\begin{quote}
\scriptsize Fig.1: Graph showing ratio of the prime-counting function $\pi(x)$ to two of its approximations, $\frac{x}{ln\ x}$ and $Li(x)$. As $x$ increases (note $x$ axis is logarithmic), both ratios tend towards $1$. The ratio for $\frac{x}{ln\ x}$ converges from above very slowly, while the ratio for $Li(x)$ converges more quickly from below.\footnote{\scriptsize cf. Prime Number Theorem. (2014, June 10). In Wikipedia, The Free Encyclopedia. Retrieved 09:53, July 9, 2014, from \href{http://en.wikipedia.org/w/index.php?title=Prime_number_theorem&oldid=612391868}{http://en.wikipedia.org/w/index.php?title\=Prime\_number\_theorem\&oldid=612391868}.}
\end{quote}

\vspace{+1ex}
\noindent The question of whether $\pi (x)/\frac{x}{log_{e}x}$ has a limit at all, or whether it oscillates, was answered---it has a limit---first by Jacques Hadamard and Charles Jean de la Vall\'{e}e Poussin independently in 1896, using advanced argumentation involving functions of a complex variable\footnote{\cite{Di52}, p.439; see also \cite{Ti51}, Chapter III, p.8 for details of Hadamard's and de la Vall\'{e}e Poussin's proofs of the Prime Number Theorem.}; and again independently by Paul Erd\"{o}s and Atle Selberg\footnote{See \cite{HW60}, p.360, Theorem 433 for a proof of Selberg's Theorem.} in 1949/1950, using only elementary---but still abstruse---methods without involving functions of a complex variable.

\subsection{A better heuristic approximation to $\pi(x)$: The integral $Li(x)$}
\label{intro.2}

We also note that, reportedly\footnote{cf. Prime Number Theorem. (2014, June 10). In Wikipedia, The Free Encyclopedia. Retrieved 09:53, July 9, 2014, from: \href{http://en.wikipedia.org/w/index.php?title=Prime_number_theorem&oldid=612391868}{http://en.wikipedia.org/w/index.php?title\=Prime\_number\_theorem\&oldid=612391868}.}:

\begin{quote}
\footnotesize
``In a handwritten note on a reprint of his 1838 paper `Sur l'usage des s\'{e}ries infinies dans la th\'{e}orie des nombres', which he mailed to Carl Friedrich Gauss, Peter Gustav Lejeune Dirichlect conjectured (under a slightly different form appealing to a series rather than an integral) that an even better approximation to $\pi(x)$ is given by the offset logarithmic integral $Li(x)$ defined by:

\begin{quote}
 $Li(x) = \int_{2}^{x}\frac{1}{log_{e}t}.dt = li(x) - li(2)$."\footnote{Where $li(x) = \int_{0}^{x}\frac{1}{log_{e}t}.dt$.}
\end{quote}
\end{quote}

\noindent We further note that in 1889 Jean de la Vall\'{e}e Poussin proved\footnote{\cite{Di52}, p.440.} (cf. Fig.1):

\begin{quote}
\footnotesize
``\ldots that $Li(x)$ represents $\pi(x)$ more exactly than $\frac{x}{log_{e}x}$ and its remaining approximations $\frac{x}{log_{e}x} + \frac{x}{log_{e}^{^{2}}x} + \ldots + \frac{(m-1)!x}{log_{e}^{^{m}}x}$."
\end{quote}

\subsection{All known approximations of $\pi(n)$ for finite values of $n$ are heuristic}

\begin{center}
\line(1,0){205} \\
\textbf{Fig.2: The distribution of the primes} \\
\line(1,0){205}

\vspace{+1ex}
\includegraphics{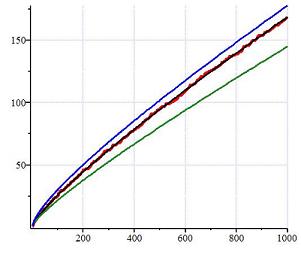}
\end{center}

\begin{quote}
\scriptsize Fig.2: The above graph compares the actual number $\pi(x)$ (red) of primes $\leq x$ with the distribution of primes as estimated variously by the functions $Li(x)$ (blue), $R(x)$ (black), and $\frac{x}{log_{e}x}$ (green), where $R(x)$ is Riemann's function $\sum_{n=1}^{\infty}\frac{\mu(n)}{(n)}li(x^{1/n})$.\footnote{cf. How Many Primes Are There? In \textit{The Prime Pages.} Retrieved 10:29, September 27, 2015, from: \\ \href{https://primes.utm.edu/howmany.html}{https://primes.utm.edu/howmany.html}.}
\end{quote}

\noindent We note that all the known approximations of $\pi(n)$ for finite values of $n$ are derived from real-valued functions that are only known to be asymptotic to $\pi(x)$, such as $\frac{x}{log_{e}x}$, $Li(x)$ and Riemann's function $R(x) = \sum_{n=1}^{\infty}\frac{\mu(n)}{(n)}li(x^{1/n})$. 

\vspace{+1ex}
\noindent Consequently, the degree of approximation for finite values of $n$ is determined only heuristically, by conjecturing upon an error term in the asymptotic relation that can be seen to yield the closest approximation upon comparison with the actual values of $\pi(n)$ within a finite range of values of $n$ (eg. Fig.2, where $n = 1000$.).

\subsection{A non-heuristic cumulative approximation of $\pi(n)$ for \textit{all} values of $n$}

\vspace{+1ex}
\noindent The question arises: Is there a function which approximates $\pi(n)$ non-heuristically for all values of $n$?

\begin{center}
\line(1,0){318} \\
\textbf{Fig.3: Density-based estimated distribution of the primes} \\
\line(1,0){318}

\vspace{+1ex}
\includegraphics{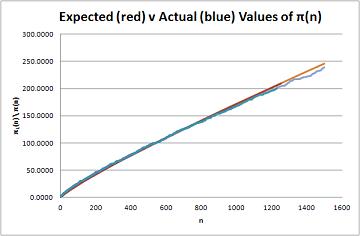}
\end{center}

\begin{quote}
\scriptsize Fig.3: The above graph compares the density-based estimated values (red) vs actual values (blue) of $\pi(n)$ for $4 \leq n \leq 1500$\footnote{See \S \ref{appendix.prim.dist.1500}, Appendix III for the values of the above plot.}, where the density-based estimated value $\pi_{_{L}}(n)$ of $\pi(n)$ is $\sum_{j=1}^{n}\prod_{i=1}^{\sqrt{j}}(1-\frac{1}{p_{_{i}}})$.
\end{quote}

\noindent In this investigation we shall address the above question by showing that the density\footnote{cf. \cite{St02}, Chapter 2, p.10.} of integers co-prime to the first $k$ primes, $p_{_{1}}, p_{_{2}}, \ldots, p_{_{k}}$, over the set of natural numbers, is:

\vspace{+1ex}
$\prod_{i = 1}^{k}(1 - \frac{1}{p_{_{i}}})$;

\vspace{+1ex}
\noindent and that the expected number of such integers in the interval ($a, b$) is thus:

\vspace{+1ex}
$(b-a)\prod_{i = 1}^{k}(1 - \frac{1}{p_{_{i}}})$,

\vspace{+1ex}
\noindent where the binomial standard deviation of the expected number of integers co-prime to $p_{_{1}}, p_{_{2}}, \ldots, p_{_{k}}$ in any interval of length ($b-a$) is:

\vspace{+2ex}
$\sqrt{(b-a)\prod_{i = 1}^{k}(1 - \frac{1}{p_{_{i}}})(1-\prod_{i = 1}^{k}(1 - \frac{1}{p_{_{i}}}))}$.

\vspace{+1ex}
\noindent Taking $(a, b)$ as the interval $(p_{_{\pi(\sqrt n)}}^{2},\ p_{_{\pi(\sqrt n)+1}}^{2})$, we conclude that a cumulative non-heuristic estimate of the number $\pi(p_{_{\pi(\sqrt{ n})+1}}^{2})$ of primes less than $p_{_{\pi(\sqrt{ n})+1}}^{2}$ is:

\vspace{+1ex}
$\pi_{_{L}}(p_{_{\pi(\sqrt{ n})+1}}^{2}) = \sum_{j=1}^{\pi(\sqrt {n})}\{(p_{_{j+1}}^{2} - p_{_{j}}^{2})\prod_{i = 1}^{j}(1 - \frac{1}{p_{_{i}}})\}$,

\vspace{+1ex}
\noindent with cumulative standard deviation:

\vspace{+2ex}
$\sum_{j=1}^{\pi(\sqrt {n})}\sqrt{(p_{_{j+1}}^{2} - p_{_{j}}^{2})\prod_{i = 1}^{j}(1 - \frac{1}{p_{_{i}}})(1-\prod_{i = 1}^{j}(1 - \frac{1}{p_{_{i}}}))}$.

\vspace{+1ex}
\noindent Moreover, a non-heuristic approximation for $\pi(n)$---with a bounded error term $> 0$ for all finite $n$---is given (Lemma \ref{sec:2.5.lem.1} and Corollary \ref{sec:2.2.lem.7.3.cor.1}) by the density-based prime counting function $\pi_{_{L}}(n)$ (cf. Fig.3):

\begin{quote}
$\pi(n) \approx \pi_{_{L}}(n) = \sum_{j = 1}^{n}\prod_{i = 1}^{\pi(\sqrt{j})}(1 - \frac{1}{p_{_{i}}}) \sim a.\frac{n}{log_{e}n} \rightarrow \infty$
\end{quote}

\vspace{+1ex}
\noindent for some constant $a > 2.e^{-\gamma} \approx 1.12292 \ldots$

\vspace{+1ex}
\noindent We then show how such a density-based approach to estimating prime counting functions non-heuristically yields elementary, density-based, proofs of Dirichlect's Theorem (Theorem \ref{dirichlect.thm}) and the Twin-Prime Conjecture (Theorem \ref{sec:twin.thm.1}).

\subsection{Expected number of primes in the interval ($p_{_{\pi(\sqrt n)}}^{2},\ p_{_{\pi(\sqrt n)+1}}^{2}$)}
\label{intro.1.3}

\begin{center}
\line(1,0){208}

\vspace{+.5ex}
\textbf{Fig.4: The graph of $y = \prod_{i = 1}^{\pi(\sqrt{x})}(1 - \frac{1}{p_{_{i}}})$}

\line(1,0){208}

\vspace{+4ex}
\begin{picture}(200,165)
\put(0,0){\line(1,0){245}}
\put(0,0){\line(0,1){180}}

\put(-17,2){\footnotesize{y $\uparrow$}}

\put(-2,80){\line(1,0){2}}
\put(-2,100){\line(1,0){2}}
\put(-2,120){\line(1,0){2}}
\put(-2,175){\line(1,0){2}}

\put(-15,76){\footnotesize{$\frac{8}{35}$}}
\put(-15,97){\footnotesize{$\frac{4}{15}$}}
\put(-15,117){\footnotesize{$\frac{1}{3}$}}
\put(-15,172){\footnotesize{$\frac{1}{2}$}}

\put(-22,-15){\footnotesize{x $\rightarrow$}}

\put(4,-4){\line(0,1){4}}
\put(8,-4){\line(0,1){4}}
\put(18,-4){\line(0,1){4}}
\put(50,-4){\line(0,1){4}}
\put(98,-4){\line(0,1){4}}
\put(242,-4){\line(0,1){4}}

\put(2,-15){\tiny{$2$}}
\put(6,-15){\tiny{$4$}}
\put(16,-15){\tiny{$9$}}
\put(48,-15){\tiny{$25$}}
\put(96,-15){\tiny{$49$}}
\put(240,-15){\tiny{$121$}}

\put(190,-10){\tiny{\textit{Not to scale}}}

\put(4,0){\framebox(14,175){}}
\put(8,0){\line(0,1){33}}
\put(8,50){\line(0,1){125}}
\put(5,40){\tiny $\frac{4}{3.5}$}
\put(18,0){\framebox(32,120)}{}
\put(30,40){\tiny $\frac{5}{5.3}$}
\put(50,0){\framebox(48,100){}}
\put(68,40){\tiny $\frac{6}{6.4}$}
\put(98,0){\framebox(144,80){}}
\put(140,40){\tiny $\frac{\pi(11^2 - 7^{2}) = 15}{\pi_{_{L}}(11^2 - 7^{2}) =16.4}$}
\put(242,0){\line(0,1){4}}
\end{picture}
\vspace{+3ex}

\begin{quote}
{\scriptsize Fig.4: Graph of $y = \prod_{i = 1}^{\pi(\sqrt{x})}(1 - \frac{1}{p_{_{i}}})$. The rectangles represent $(p_{_{j+1}}^{2} - p_{_{j}}^{2})\prod_{i = 1}^{j}(1 - \frac{1}{p_{_{i}}})$ for $j \geq 1$. Figures within each rectangle are the primes corresponding to the functions $\pi(n)$ and $\pi_{_{L}}(n)$ within the interval $(p_{_{j}}^{2},\ p_{_{j+1}}^{2})$ for $j \geq 2$. The area under the curve is $\pi_{_{L}}(x) = (x - p_{_{n}}^{2})\prod_{i = 1}^{n}(1 - \frac{1}{p_{_{i}}}) + \sum_{j = 1}^{n-1}(p_{_{j+1}}^{2} - p_{_{j}}^{2})\prod_{i = 1}^{j}(1 - \frac{1}{p_{_{i}}}) + 2$ (see Fig.5). }
\end{quote}
\end{center}

\vspace{+2ex}
\noindent More specifically, since $n$ is a prime if, and only if, it is not divisible by any prime $p \leq \sqrt{n}$:

\begin{quote}
\noindent (i) a non-heuristic estimate for the number of primes in the interval $(p_{_{\pi(\sqrt n)}}^{2},\ p_{_{\pi(\sqrt n)+1}}^{2})$ is (Theorem \ref{sec:2.5.cor.1}) the expected number of primes in the interval, given in terms of a prime counting function $\pi_{_{L}}(n)$ as (cf. Fig.4):

\vspace{+1ex}
\begin{quote}
$\pi_{_{L}}(p_{_{\pi(\sqrt n)+1}}^{2}) - \pi_{_{L}}(p_{_{\pi(\sqrt n)}}^{2}) = \{(p_{_{\pi(\sqrt n)+1}}^{2} - p_{_{\pi(\sqrt n)}}^{2})\prod_{i = 1}^{\pi(\sqrt{n})}(1 - \frac{1}{p_{_{i}}})\}$
\end{quote}
\end{quote}

\begin{quote}
\noindent (ii) and a cumulative non-heuristic approximation of the number $\pi(n)$ of primes less than or equal to $n$ (Lemma \ref{sec:2.5.lem.1} and Corollary \ref{sec:2.2.lem.7.3}) is the prime counting function $\pi_{_{L}}(n)$ (cf. Fig.5)\footnote{Compare \cite{HL23}, pp.36-37. See also \S \ref{appendix.prim.dist.1500}, Appendix III for the estimated values $\pi_{_{L}}(n)$, and the actual values $\pi(n)$, for $4 \leq n \leq 1500$.}:

\vspace{+1ex}
\begin{quote}
$\pi(n) \approx \pi_{_{L}}(n) = \sum_{j = 1}^{n}\prod_{i = 1}^{\pi(\sqrt{j})}(1 - \frac{1}{p_{_{i}}})$.
\end{quote}
\end{quote}

\begin{center}
\line(1,0){160}

\textbf{Fig.5: The graph of $y = \pi_{_{L}}(x)$}

\line(1,0){160}

\vspace{+4ex}
\begin{picture}(200,165)
\put(0,0){\line(1,0){245}}
\put(0,0){\line(0,1){165}}

\put(-17,2){\footnotesize{y $\uparrow$}}

\put(-2,17.5){\line(1,0){2}}
\put(-2,44){\line(1,0){2}}
\put(-2,76){\line(1,0){2}}
\put(-2,158){\line(1,0){2}}

\put(-20,15.5){\footnotesize{$3.5$}}
\put(-20,42){\footnotesize{$8.8$}}
\put(-25,74){\footnotesize{$15.2$}}
\put(-25,156){\footnotesize{$31.6$}}

\put(-22,-15){\footnotesize{x $\rightarrow$}}

\put(4,-4){\line(0,1){4}}
\put(8,-4){\line(0,1){4}}
\put(18,-4){\line(0,1){4}}
\put(50,-4){\line(0,1){4}}
\put(98,-4){\line(0,1){4}}
\put(242,-4){\line(0,1){4}}

\put(2,-15){\tiny{$2$}}
\put(6,-15){\tiny{$4$}}
\put(16,-15){\tiny{$9$}}
\put(48,-15){\tiny{$25$}}
\put(96,-15){\tiny{$49$}}
\put(240,-15){\tiny{$121$}}

\put(190,-10){\tiny{\textit{Not to scale}}}

\put(18,0){\line(0,1){17.5}}
\put(50,0){\line(0,1){44}}
\put(98,0){\line(0,1){76}}
\put(242,0){\line(0,1){158}}

\qbezier(8,0)(18,17.5)(18,17.5)
\qbezier(18,17.5)(50,44)(50,44)
\qbezier(50,44)(98,76)(98,76)
\qbezier(98,76)(242,158)(242,158)

\put(150,115){\tiny{$\frac{8}{35}$}}
\put(60,62){\tiny{$\frac{4}{15}$}}
\put(27,35){\tiny{$\frac{1}{3}$}}
\put(7.5,12){\tiny{$\frac{1}{2}$}}

\put(242,0){\line(0,1){4}}
\end{picture}
\end{center}

\vspace{+3ex}
\begin{quote}
{\scriptsize Fig.5: Graph of $y = \pi_{_{L}}(x) = (x - p_{_{n}}^{2})\prod_{i = 1}^{n}(1 - \frac{1}{p_{_{i}}}) + \sum_{j = 1}^{n-1}(p_{_{j+1}}^{2} - p_{_{j}}^{2})\prod_{i = 1}^{j}(1 - \frac{1}{p_{_{i}}}) + 2$ in the interval $(p_{_{n}}^{2},\ p_{_{n+1}}^{2})$. Note that the gradient in the interval $(p_{_{n}}^{2},\ p_{_{n+1}}^{2})$ is $\prod_{i = 1}^{n}(1 - \frac{1}{p_{_{i}}}) \rightarrow 0$.}
\end{quote}

\section{Density-based non-heuristic approximations of prime counting functions}
\label{intro.1.7}

\noindent In the rest of this investigation we formally consider elementary, density-based, arguments for:

\vspace{+1ex}
\begin{quote}
\noindent (i) \textit{Dirichlect's Theorem}: We show that the number $\pi_{_{(a, d)}}(n)$ of Dirichlect primes of the form $a+m.d$ which are less than or equal to $n$, where $a, d$ are co-prime and $1 \leq a < d = q_{_{1}}^{\alpha_{_{1}}}.q_{_{2}}^{\alpha_{_{2}}} \ldots q_{_{k}}^{\alpha_{_{k}}}$ ($q_{_{i}}$ prime), is non-heuristically approximated by the cumulative Dirichlect prime counting function $\pi_{_{D}}(n)$ (Definition \ref{sec:5.lem.6}), such that:

\begin{quote}
$\pi_{_{(a, d)}}(n) \approx \pi_{_{D}}(n) = \prod_{i = 1}^{k}{\frac{1}{q_{_{i}}^{\alpha_{_{i}}}}}.\prod_{i = 1}^{k}(1 - \frac{1}{q_{_{i}}})^{-1}.\pi_{_{L}}(n) \rightarrow \infty$.
\end{quote}

\noindent (ii) \textit{Twin Prime Theorem}: We show that there are infinitely many twin primes since a cumulative non-heuristic approximation of the number $\pi_{_{2}}(p_{_{k+1}}^{2})$ of twin primes $\leq p_{_{k+1}}^{2}$ for all $k \geq 1$ is:

\begin{quote}
$\sum_{j=9}^{p_{_{k+1}}^{2}} \prod_{i=2}^{\pi(\sqrt{j})-1}(1 - \frac{2}{p_{_{i}}}) \rightarrow \infty$.
\end{quote}
\end{quote}

\subsection{The residues $r_{i}(n)$.}
\label{res}

\vspace{+2ex}
\noindent We begin by formally defining the residues $r_{i}(n)$ for all $n \geq 2$ and all $i \geq 2$ as below\footnote{The residues $r_{i}(n)$ can also be graphically displayed variously as shown in the Appendix II in \S\ref{appendix}.}:

\begin{definition}
\label{sec:4.lem.1}
$n + r_{i}(n) \equiv 0\ (mod\ i)$ where $i > r_{i}(n) \geq 0$.  
\end{definition}

\noindent Since each residue $r_{i}(n)$ cycles over the $i$ values $(i-1, i-2, \ldots, 0)$, these values are all incongruent and form a complete system of residues\footnote{\cite{HW60}, p.49.} $mod\ i$.

\vspace{+1ex}
\noindent It immediately follows that:

\begin{lemma}
\label{sec:4.lem.1.1}
$r_{i}(n) = 0$ if, and only if, $i$ is a divisor of $n$. \hfill $\Box$
\end{lemma}

\subsection{The probability model $\mathbb{M}_{i} = \{(0, 1, 2, \ldots, i-1),\ r_{i}(n), \frac{1}{i}\}$}
\label{res.1}

\vspace{+2ex}
\noindent By the standard definition of the probability $\mathbb{P}(e)$ of an event $e$\footnote{See \cite{Ko56}, Chapter I, \S 1, Axiom III, pg.2.}, we have by Lemma \ref{sec:4.lem.1.1} that:

\begin{lemma}
\label{sec:4.lem.2}
For any $n \geq 2,\ i \geq 2$ and any given integer $i > u \geq 0$:

\vspace{+1ex}
$\bullet$ the probability $\mathbb{P}(r_{i}(n) = u)$ that $r_{i}(n) = u$ is $\frac{1}{i}$;

\vspace{+1ex}
$\bullet$ $\sum_{u=0}^{u=i-1}\mathbb{P}(r_{i}(n) = u) = 1$;

\vspace{+1ex}
$\bullet$ and the probability $\mathbb{P}(r_{i}(n) \neq u)$ that $r_{i}(n) \neq u$ is $1 - \frac{1}{i}$. \hfill $\Box$
\end{lemma}

\noindent By the standard definition of a probability model, we conclude that:

\begin{theorem}
\label{prob.model}
For any $i \geq 2$, $\mathbb{M}_{i} = \{(0, 1, 2, \ldots, i-1),\ r_{i}(n), \frac{1}{i}\}$ yields a probability model for each of the values of $r_{i}(n)$. \hfill $\Box$
\end{theorem}

\begin{corollary}
\label{sec:4.lem.2.cor.1}
For any $n \geq 2$ and any prime $p \geq 2$, the probability $\mathbb{P}(r_{p}(n) = 0)$ that $r_{p}(n) = 0$, and that $p$ divides $n$, is $\frac{1}{p}$; and the probability $\mathbb{P}(r_{p}(n) \neq 0)$ that $r_{p}(n) \neq 0$, and that $p$ does not divide $n$, is $1 - \frac{1}{p}$. \hfill $\Box$
\end{corollary}

\noindent We also note the standard definition\footnote{See \cite{Ko56}, Chapter VI, \S 1, Definition 1, pg.57 and \S 2, pg.58.}:

\begin{definition}
\label{sec:4.def.1}
Two events $e_{i}$ and $e_{j}$ are mutually independent for $i \neq j$ if, and only if, $\mathbb{P}(e_{i}\ \cap\ e_{j}) = \mathbb{P}(e_{i}).\mathbb{P}(e_{j})$.
\end{definition}

\vspace{+1ex}
\subsection{The prime divisors of any integer $n$ are mutually independent}
\label{res.2}

\vspace{+2ex}
\noindent We then have that:

\begin{lemma}
\label{sec:4.lem.3}
If $p_{_{i}}$ and $p_{_{j}}$ are two primes where $i \neq j$ then, for any $n \geq 2$, we have:

\vspace{+1ex}
\begin{quote}
$\mathbb{P}((r_{p_{_{i}}}(n) = u) \cap (r_{p_{_{j}}}(n) = v)) = \mathbb{P}(r_{p_{_{i}}}(n) = u).\mathbb{P}(r_{p_{_{j}}}(n) = v)$
\end{quote}

\vspace{+1ex}
\noindent where $p_{_{i}} > u \geq 0$ and $p_{_{j}} > v \geq 0$.
\end{lemma}

\noindent \textbf{Proof}: The $p_{_{i}}.p_{_{j}}$ numbers $v.p_{_{i}} + u.p_{_{j}}$, where $p_{_{i}} > u \geq 0$ and $p_{_{j}} > v \geq 0$, are all incongruent and form a complete system of residues\footnote{\cite{HW60}, p.52, Theorem 59.} $mod\ (p_{_{i}}.p_{_{j}})$. Hence:

\vspace{+1ex}
\begin{quote}
$\mathbb{P}((r_{p_{_{i}}}(n) = u) \cap (r_{p_{_{j}}}(n) = v)) = \frac{1}{p_{_{i}}.p_{_{j}}}$
\end{quote}

\vspace{+1ex}
\noindent By Lemma \ref{sec:4.lem.2}:

\vspace{+1ex}
\begin{quote}
$\mathbb{P}(r_{p_{_{i}}}(n) = u).\mathbb{P}(r_{p_{_{j}}}(n) = v) = (\frac{1}{p_{_{i}}})(\frac{1}{p_{_{j}}})$. 
\end{quote}

\vspace{+1ex}
\noindent The lemma follows. \hfill $\Box$

\vspace{+1ex}
\noindent If $u = 0$ and $v = 0$ in Lemma \ref{sec:4.lem.3}, so that both $p_{_{i}}$ and $p_{_{j}}$ are prime divisors of $n$, we immediately conclude by Definition \ref{sec:4.def.1} that:

\begin{corollary}
\label{sec:4.lem.2.cor.2}
$\mathbb{P}((r_{p_{_{i}}}(n) = 0) \cap (r_{p_{_{j}}}(n) = 0)) = \mathbb{P}(r_{p_{_{i}}}(n) = 0).\mathbb{P}(r_{p_{_{j}}}(n) = 0)$. \hfill $\Box$
\end{corollary}

\noindent We can also express this as:

\begin{corollary}
\label{sec:4.lem.2.cor.3}
$\mathbb{P}(p_{_{i}} | n\ \cap\ p_{_{j}} | n) = \mathbb{P}(p_{_{i}} | n).\mathbb{P}(p_{_{j}} | n)$. \hfill $\Box$ 
\end{corollary}

\noindent We thus conclude that:

\begin{theorem}
\label{sec:4.thm.1}
The prime divisors of any integer $n$ are mutually independent. \hfill $\Box$
\end{theorem}

\vspace{+1ex}
\subsubsection{Integer Factorising cannot be polynomial-time}
\label{int-fact}

\vspace{+2ex}
\noindent We digress briefly from our investigation of prime counting functions to note that Theorem \ref{sec:4.thm.1} immediately yields the computational complexity consequence\footnote{cf. \cite{Cook}.} that no deterministic algorithm\footnote{A deterministic algorithm computes a mathematical function which has a unique value for any input in its domain, and the algorithm is a process that produces this particular value as output.} can compute a factor of any randomly given integer $n$ in polynomial time\footnote{cf. \cite{Cook}, p.1; also \cite{Br00}, p.1, fn.1.}.

\vspace{+1ex}
\noindent We note the standard definition\footnote{cf. \cite{Cook}, p.1; also \cite{Br00}, p.1, fn.1: ``For a polynomial-time algorithm the expected running time should be a polynomial in the length of the input, i.e. $O((log N)^{c})$ for some constant $c$".}:

\begin{definition}
\label{sec:1.def.1}
A deterministic algorithm computes a number-theoretical function $f(n)$ in polynomial-time if there exists $k$ such that, for all inputs $n$, the algorithm computes $f(n)$ in $\leq (log_{e}\ n)^{k} + k$ steps. 
\end{definition}

\noindent It then follows from Theorem \ref{sec:4.thm.1} that:

\begin{corollary}
\label{sec:1.cor.1}
Any deterministic algorithm that always computes a prime factor of $n$ cannot be polynomial-time.
\end{corollary}

\noindent \textbf{Proof}: Any computational process that successfully identifies a prime divisor of $n$ must necessarily appeal to at least one logical operation for identifying such a factor.

\vspace{+1ex}
\noindent Since $n$ is a prime if, and only if, it is not divisible by any prime $p \leq \sqrt{n}$, and $n$ may be the square of a prime, it follows from Theorem \ref{sec:4.thm.1} that we necessarily require at least one logical operation for each prime $p \leq \sqrt{n}$ in order to logically determine whether $p$ is a prime divisor of $n$. 

\vspace{+1ex}
\noindent Since the number of such primes is of the order $O(\sqrt{n}/log_{e}\ n)$, the number of computations required by any deterministic algorithm that always computes a prime factor of $n$ cannot be polynomial-time---i.e. of order $O((log_{e}\ n)^{c})$ for any $c$---in the length of the input $n$. The corollary follows. \hfill $\Box$

\vspace{+1ex}
\subsection{Density of integers not divisible by primes $Q = \{q_{_{1}}, q_{_{2}}, ..., q_{_{k}}\}$}
\label{sec:2.2}

Reverting back to our consideration of prime distribution, we conclude from Lemma \ref{sec:4.lem.2} and Lemma \ref{sec:4.lem.3} that:

\begin{lemma}
\label{sec:2.2.lem.0}
The density of the set of all integers that are not divisible by any of a given set of primes $Q = \{q_{_{1}}, q_{_{2}}, ..., q_{_{k}}\}$ is:

\vspace{+2ex}
$\prod_{q \in Q}(1-1/q)$. \hfill $\Box$
\end{lemma}

\vspace{+1ex}
\noindent It follows that:

\begin{lemma}
\label{sec:2.2.lem.2}
The expected number of integers in any interval (a,b) that are not divisible by any of a given set of primes $Q = \{q_{_{1}}, q_{_{2}}, ..., q_{_{k}}\}$ is:

\vspace{+2ex}
$(b-a)\prod_{q \in Q}(1-1/q)$. \hfill $\Box$
\end{lemma}

\vspace{+1ex}
\subsection{The function $\pi_{_{H}}(n)$}
\label{sec:Hpi}

\vspace{+2ex}
\noindent In particular, the expected number $\pi_{_{H}}(n)$ of integers $\leq n$ that are not divisible by any of the first $k$ primes $p_{_{1}}, p_{_{2}}, ..., p_{_{k}}$ is:

\begin{corollary}
\label{sec:2.2.lem.2.def.1}
$\pi_{_{H}}(n) = n.\prod_{i = 1}^{k}(1 - \frac{1}{p_{_{i}}})$.
\end{corollary}

\noindent It follows that:

\begin{corollary}
\label{sec:2.5.cor.2.H}
\noindent The expected number of primes $\leq p_{_{\pi(\sqrt n)+1}}^{2}$ is:

\vspace{+1ex}
$\pi_{_{H}}(p_{_{\pi(\sqrt n)+1}}^{2}) = p_{_{\pi(\sqrt n) +1}}^{2}\prod_{i = 1}^{\pi(\sqrt n)}(1 - \frac{1}{p_{_{i}}})$

\vspace{+1ex}
\noindent with cumulative standard deviation:

\vspace{+1.5ex}
$p_{_{\pi(\sqrt n)+1}}\sqrt{\prod_{i = 1}^{\pi(\sqrt n)}(1 - \frac{1}{p_{_{i}}})(1-\prod_{i = 1}^{\pi(\sqrt n)}(1 - \frac{1}{p_{_{i}}}))}$. \hfill $\Box$
\end{corollary}

\noindent We conclude that $\pi_{_{H}}(n)$ is the non-heuristic approximation of the number of primes $\leq n$\footnote{Fig.12 in \S \ref{appendix.prim.dist.1500} compares the values of $\pi(n)$ and $\pi_{_{H}}(n)$ for $4 \leq n \leq 1500$.}:

\begin{lemma}
\label{sec:2.5.lem.H.1}
$\pi(n) \approx \pi_{_{H}}(n) = n.\prod_{i = 1}^{\pi(\sqrt{n})}(1 - \frac{1}{p_{_{i}}})$.
\end{lemma}

\vspace{+1ex}
\subsection{The function $\pi_{_{L}}(n)$}
\label{sec:Lpi}

\vspace{+2ex}
\noindent It also follows immediately from Theorem \ref{sec:2.2.lem.2} that:

\begin{corollary}
\label{sec:2.5.cor.1}
\noindent The expected number of primes in the interval ($p_{_{\pi(\sqrt n)}}^{2},\ p_{_{\pi(\sqrt n)+1}}^{2}$) is:

\vspace{+1ex}
$(p_{_{\pi(\sqrt n)+1}}^{2} - p_{_{\pi(\sqrt n)}}^{2})\prod_{i = 1}^{\pi(\sqrt n)}(1 - \frac{1}{p_{_{i}}})$

\vspace{+1ex}
\noindent with standard binomial deviation:

\vspace{+1.5ex}
$\sqrt{(p_{_{\pi(\sqrt n)+1}}^{2} - p_{_{\pi(\sqrt n)}}^{2})\prod_{i = 1}^{\pi(\sqrt n)}(1 - \frac{1}{p_{_{i}}})(1-\prod_{i = 1}^{\pi(\sqrt n)}(1 - \frac{1}{p_{_{i}}}))}$. \hfill $\Box$
\end{corollary}

\vspace{+1ex}
\noindent It further follows from Lemma \ref{sec:2.2.lem.2} and Corollary \ref{sec:2.5.cor.1} that:

\begin{corollary}
\label{sec:2.5.cor.2}
\noindent The number $\pi(p_{_{\pi(\sqrt n)+1}}^{2})$ of primes less than $p_{_{\pi(\sqrt n)+1}}^{2}$ is cumulatively approximated by:

\vspace{+1ex}
$\pi_{_{L}}(p_{_{\pi(\sqrt n)+1}}^{2}) = \sum_{j=1}^{\pi(\sqrt n)}\{(p_{_{j+1}}^{2} - p_{_{j}}^{2})\prod_{i = 1}^{j}(1 - \frac{1}{p_{_{i}}})\}$

\vspace{+1ex}
\noindent with cumulative standard deviation:

\vspace{+1.5ex}
$\sum_{j=1}^{\pi(\sqrt n)}\sqrt{(p_{_{j+1}}^{2} - p_{_{j}}^{2})\prod_{i = 1}^{j}(1 - \frac{1}{p_{_{i}}})(1-\prod_{i = 1}^{j}(1 - \frac{1}{p_{_{i}}}))}$. \hfill $\Box$
\end{corollary}

\noindent We conclude that $\pi_{_{L}}(n)$ is the cumulative non-heuristic approximation of the number of primes $\leq n$\footnote{Fig.12 in \S \ref{appendix.prim.dist.1500}, and Fig.13 in \S \ref{hyp.obs}, comparatively analyse the values of $\pi(n)$ and $\pi_{_{L}}(n)$ for $4 \leq n \leq 1500$.}:

\begin{lemma}
\label{sec:2.5.lem.1}
$\pi(n) \approx \pi_{_{L}}(n) = \sum_{j = 1}^{n}\prod_{i = 1}^{\pi(\sqrt{j})}(1 - \frac{1}{p_{_{i}}})$.
\end{lemma}

\noindent It immediately follows from Lemma \ref{sec:2.5.lem.H.1} and Lemma \ref{sec:2.5.lem.1} that\footnote{We show in Appendix I, Lemma \ref{anomaly.lem}, that $\pi_{_{L}}(n)$ is a better approximation of $\pi(n)$ than $\pi_{_{H}}(n)$ for $n \geq 9$.}:

\begin{corollary}
\label{sec:2.5.lem.1.cor}
$\pi_{_{L}}(n) > \pi_{_{H}}(n)$ for all $n \geq 9$.
\end{corollary}

\vspace{+1ex}
\subsection{The interval $(p_{_{n}}^{2},\ p_{_{n+1}}^{2})$}
\label{sec:2.2.0.1}

\vspace{+2ex}
\noindent It follows immediately from the definition of $\pi(x)$ as the number of primes less than or equal to $x$ that: 

\begin{lemma}
\label{sec:2.2.lem.7.1}
$\prod_{i = 1}^{\pi(\sqrt{x})}(1 - \frac{1}{p_{_{i}}}) = \prod_{i = 1}^{\pi(\sqrt{x+1})}(1 - \frac{1}{p_{_{i}}})$ for $p_{n}^{2} \leq x < p_{n+1}^{2}$. \hfill $\Box$
\end{lemma}

\noindent We can thus generalise the number-theoretic function of Lemma \ref{sec:2.5.lem.1} as the real-valued function:

\begin{definition}
\label{sec:2.2.lem.7.2}
$\pi_{_{L}}(x) = \pi_{_{L}}(p_{_{n}}^{2}) + (x - p_{_{n}}^{2})\prod_{i = 1}^{n}(1 - \frac{1}{p_{_{i}}})$ for $p_{n}^{2} \leq x < p_{n+1}^{2}$. \hfill $\Box$
\end{definition}

\noindent We note that the graph of $\pi_{_{L}}(x)$ in the interval $(p_{_{n}}^{2},\ p_{_{n+1}}^{2})$ for $n \geq 1$ is now a straight line with gradient $\prod_{i = 1}^{n}(1 - \frac{1}{p_{_{i}}})$, as illustrated in \S \ref{intro.1.3}, Fig.5 where we defined $\pi_{_{L}}(x)$ equivalently by:

\vspace{+1ex}
$\pi_{_{L}}(x) = (x - p_{_{n}}^{2})\prod_{i = 1}^{n}(1 - \frac{1}{p_{_{i}}}) + \sum_{j = 1}^{n-1}(p_{_{j+1}}^{2} - p_{_{j}}^{2})\prod_{i = 1}^{j}(1 - \frac{1}{p_{_{i}}}) + 2$

\vspace{+1ex}
\subsection{The functions $\pi_{_{L}}(x)/\frac{x}{log_{e}x}$ and $\pi_{_{H}}(x)/\frac{x}{log_{e}x}$}
\label{sec:2.2.0.7}

\vspace{+2ex}
\noindent We consider next the function $\pi_{_{L}}(x)/\frac{x}{log_{e}x}$ in the interval $(p_{_{n}}^{2},\ p_{_{n+1}}^{2})$:

\vspace{+1ex}
$\pi_{_{L}}(x)/\frac{x}{log_{e}x} = (\pi_{_{L}}(p_{_{n}}^{2}) + (x - p_{_{n}}^{2})\prod_{i = 1}^{n}(1 - \frac{1}{p_{_{i}}}))/\frac{x}{log_{e}x}$

\vspace{+1ex}
\noindent This now yields the derivative $(\pi_{_{L}}(x).\frac{log_{e}x}{x})'$ in the interval $(p_{_{n}}^{2},\ p_{_{n+1}}^{2})$ as:

\vspace{+1ex}
$\pi_{_{L}}(x).(\frac{log_{e}x}{x})' + (\pi_{_{L}}(x))'.\frac{log_{e}x}{x}$

\vspace{+1ex}
$(\pi_{_{L}}(p_{_{n}}^{2}) + (x - p_{_{n}}^{2})\prod_{i = 1}^{n}(1 - \frac{1}{p_{_{i}}})).(\frac{log_{e}x}{x})' + (\pi_{_{L}}(p_{_{n}}^{2}) + (x - p_{_{n}}^{2})\prod_{i = 1}^{n}(1 - \frac{1}{p_{_{i}}}))'.\frac{log_{e}x}{x}$

\vspace{+1ex}
$(\pi_{_{L}}(p_{_{n}}^{2}) + (x - p_{_{n}}^{2})\prod_{i = 1}^{n}(1 - \frac{1}{p_{_{i}}})).(\frac{1}{x^{2}} - \frac{log_{e}x}{x^{2}}) + (\prod_{i = 1}^{n}(1 - \frac{1}{p_{_{i}}})).\frac{log_{e}x}{x}$

\vspace{+1ex}
\noindent Since $p_{n}^{2} \leq x < p_{n+1}^{2}$, by Mertens'\footnote{\cite{HW60}, Theorem 429, p.351.} and Chebyshev's Theorems we can express the above as:

\vspace{+1ex}
$\sim (\pi_{_{L}}(p_{_{n}}^{2}) + \frac{e^{-\gamma}(x - p_{_{n}}^{2})}{log_{e}n}).(\frac{1}{x^{2}} - \frac{log_{e}x}{x^{2}}) + \frac{e^{-\gamma}.log_{e}x}{x.log_{e}n}$

\vspace{+1ex}
$\sim (\frac{\pi_{_{L}}(p_{_{n}}^{2})}{x} + \frac{e^{-\gamma}}{log_{e}n}(1 - \frac{p_{_{n}}^{2}}{x})).\frac{(1 - log_{e}x)}{x} + \frac{e^{-\gamma}.log_{e}x}{x.log_{e}n}$

\vspace{+1ex}
$\sim (\frac{\pi_{_{L}}(p_{_{n}}^{2})}{p_{_{n}}^{2}}.\frac{p_{_{n}}^{2}}{x} + \frac{e^{-\gamma}}{log_{e}n}(1 - \frac{p_{_{n}}^{2}}{x})).\frac{(1 - 2.log_{e}p_{_{n}})}{p_{_{n}}^{2}} + \frac{2.e^{-\gamma}.log_{e}p_{_{n}}}{p_{_{n}}^{2}.log_{e}n}$

\vspace{+1ex}
\noindent Since each term $\rightarrow 0$ as $n \rightarrow \infty$, we conclude that the function $\pi_{_{L}}(x)/\frac{x}{log_{e}x}$ does not oscillate but tends to a limit as $x \rightarrow \infty$ since:

\begin{lemma}
\label{sec:2.2.lem.7.3}
$(\pi_{_{L}}(x)/\frac{x}{log_{e}x})' \in o(1)$. \hfill $\Box$
\end{lemma}

\noindent We further conclude that:

\begin{corollary}
\label{sec:2.2.lem.7.3.cor.1}
$\pi_{_{L}}(n) = \sum_{j = 1}^{n}\prod_{i = 1}^{\pi(\sqrt{j})}(1 - \frac{1}{p_{_{i}}}) \sim a.\frac{n}{log_{e}n}$ for some constant $a$. \hfill $\Box$
\end{corollary}

\noindent We note that $a > 2.e^{-\gamma}$\footnote{ Where $2.e^{-\lambda} \approx 1.12292 \ldots$; \cite{Gr95}, p.13.}, since $\prod_{i = 1}^{\pi(\sqrt{j})}(1 - \frac{1}{p_{_{i}}}) \geq \prod_{i = 1}^{\pi(\sqrt n)}(1 - \frac{1}{p_{_{i}}})$ for all $1 \leq j \leq n$, and it follows from Definition \ref{sec:2.2.lem.2.def.1} that:

\begin{corollary}
\label{sec:2.2.lem.7.3.cor.2}
$\pi_{_{H}}(n) = n.\prod_{i = 1}^{\pi(\sqrt n)}(1 - \frac{1}{p_{_{i}}}) \sim 2.e^{-\gamma}.\frac{n}{log_{e}n}$\footnote{By Mertens' Theorem; since $log_{_{e}}\pi(\sqrt n) \sim (log_{_{e}}\sqrt n - log_{_{e}}log_{_{e}}\sqrt n)$ by the Prime Number Theorem.}. \hfill $\Box$
\end{corollary}

\vspace{+1ex}
\subsection{Primes in an arithmetic progression}
\label{dirichlect}

\vspace{+2ex}
\noindent We consider now Dirichlect's Theorem, which is the assertion that if $a$ and $d$ are co-prime and $1 \leq a < d$, then the arithmetic progression $a + m.d$, where $m \geq 1$, contains an infinitude of (Dirichlect) primes. 

\vspace{+1ex}
\noindent We first note that Lemma \ref{sec:4.lem.3} can be extended to prime powers in general\footnote{\textit{Hint}: The following arguments may be easier to follow if we visualise the residues $r_{p_{_{i}}^{\alpha}}(n)$ and $r_{p_{_{i}}^{\beta}}(n)$ as they would occur in \S\ref{sec:4.fig1}, Fig.7 and Fig.8.}:

\begin{lemma}
\label{sec:5.lem.1}
If $p_{_{i}}$ and $p_{_{j}}$ are two primes where $i \neq j$ then, for any $n \geq 2,\ \alpha, \beta \geq 1$, we have:

\vspace{+1ex}
\begin{quote}
$\mathbb{P}((r_{p_{_{i}}^{\alpha}}(n) = u) \cap (r_{p_{_{j}}^{\beta}}(n) = v)) = \mathbb{P}(r_{p_{_{i}}^{\alpha}}(n) = u).\mathbb{P}(r_{p_{_{j}}^{\beta}}(n) = v)$
\end{quote}

\vspace{+1ex}
\noindent where $p_{_{i}}^{\alpha} > u \geq 0$ and $p_{_{j}}^{\beta} > v \geq 0$.
\end{lemma}

\noindent \textbf{Proof}: The $p_{_{i}}^{\alpha}.p_{_{j}}^{\beta}$ numbers $v.p_{_{i}}^{\alpha} + u.p_{_{j}}^{\beta}$, where $p_{_{i}}^{\alpha} > u \geq 0$ and $p_{_{j}}^{\beta} > v \geq 0$, are all incongruent and form a complete system of residues\footnote{\cite{HW60}, p.52, Theorem 59.} $mod\ (p_{_{i}}^{\alpha}.p_{_{j}}^{\beta})$. Hence:

\vspace{+1ex}
\begin{quote}
$\mathbb{P}((r_{p_{_{i}}^{\alpha}}(n) = u) \cap (r_{p_{_{j}}^{\beta}}(n) = v)) = \frac{1}{p_{_{i}}^{\alpha}.p_{_{j}}^{\beta}}$
\end{quote}

\vspace{+1ex}
\noindent By Lemma \ref{sec:4.lem.2}:

\vspace{+1ex}
\begin{quote}
$\mathbb{P}(r_{p_{_{i}}^{\alpha}}(n) = u).\mathbb{P}(r_{p_{_{j}}^{\beta}}(n) = v) = (\frac{1}{p_{_{i}}^{\alpha}})(\frac{1}{p_{_{j}}^{\beta}})$. 
\end{quote}

\vspace{+1ex}
\noindent The lemma follows. \hfill $\Box$

\vspace{+1ex}
\noindent If $u = 0$ and $v = 0$ in Lemma \ref{sec:5.lem.1}, so that both $p_{_{i}}$ and $p_{_{j}}$ are prime divisors of $n$, we immediately conclude by Definition \ref{sec:4.def.1} that:

\begin{corollary}
\label{sec:5.lem.1.cor.1}
$\mathbb{P}((r_{p_{_{i}}^{\alpha}}(n) = 0) \cap (r_{p_{_{j}}^{\beta}}(n) = 0)) = \mathbb{P}(r_{p_{_{i}}^{\alpha}}(n) = 0).\mathbb{P}(r_{p_{_{j}}{\beta}}(n) = 0)$. \hfill $\Box$
\end{corollary}

\noindent We can also express this as:

\begin{corollary}
\label{sec:5.lem.1.cor.2}
$\mathbb{P}(p_{_{i}}^{\alpha} | n\ \cap\ p_{_{j}}^{\beta} | n) = \mathbb{P}(p_{_{i}}^{\alpha} | n).\mathbb{P}(p_{_{j}}^{\beta} | n)$. \hfill $\Box$ 
\end{corollary}

\noindent We thus conclude that:

\begin{theorem}
\label{sec:5.thm.1}
For any two primes $p \neq q$ and natural numbers $n, \alpha, \beta \geq 1$, whether or not $p^{\alpha}$ divides $n$ is independent of whether or not $q^{\beta}$ divides $n$. \hfill $\Box$
\end{theorem}

\vspace{+1ex}
\subsubsection{The probability that $n$ is a prime of the form $a + m.d$}
\label{dirichlect:probability}

\vspace{+2ex}
\noindent We note next that:

\begin{lemma}
\label{sec:5.lem.2}
For any co-prime natural numbers $1 \leq a < d = q_{_{1}}^{\alpha_{_{1}}}.q_{_{2}}^{\alpha_{_{2}}} \ldots q_{_{k}}^{\alpha_{_{k}}}$ where:

\vspace{+1ex}
\begin{quote}
 $q_{_{1}} < q_{_{2}} < \ldots < q_{_{k}}$ are primes and $\alpha_{_{1}}, \alpha_{_{2}} \ldots \alpha_{_{k}} \geq 1$ are natural numbers;
\end{quote}

\vspace{+1ex}
\noindent the natural number $n$ is of the form $a+m.d$ for some natural number $m \geq 1$ if, and only if:

\vspace{+1ex}
\begin{quote}
$a + r_{_{q_{_{i}}^{\alpha_{_{i}}}}}(n) \equiv 0\ (mod\ q_{_{i}}^{\alpha_{_{i}}})$ for all $1 \leq i \leq k$
\end{quote}

\vspace{+1ex}
\noindent where $0 \leq r_{_{i}}(n) < i$ is defined for all $i > 1$ by:

\vspace{+1ex}
\begin{quote}
$n + r_{_{i}}(n) \equiv 0\ (mod\ i)$ .
\end{quote}
\end{lemma}

\noindent \textbf{Proof}: First, if $n$ is of the form $a+m.d$ for some natural number $m \geq 1$, where $1 \leq a < d = q_{_{1}}^{\alpha_{_{1}}}.q_{_{2}}^{\alpha_{_{2}}} \ldots q_{_{k}}^{\alpha_{_{k}}}$, then:

\vspace{+1ex}
\begin{quote}
$\begin{array}{rlcll}
& n & \equiv & a\ (mod\ d) & \\
and: & n + r_{_{q_{_{i}}^{\alpha_{_{i}}}}}(n) & \equiv & 0\ (mod\ q_{_{i}}^{\alpha_{_{i}}}) & for\ all\ 1 \leq i \leq k \\
whence: & a + r_{_{q_{_{i}}^{\alpha_{_{i}}}}}(n) & \equiv & 0\ (mod\ q_{_{i}}^{\alpha_{_{i}}}) & for\ all\ 1 \leq i \leq k
\end{array}$
\end{quote}

\vspace{+1ex}
\noindent Second:

\vspace{+1ex}
\begin{quote}
$\begin{array}{rlcll}
If: & a + r_{_{q_{_{i}}^{\alpha_{_{i}}}}}(n) & \equiv & 0\ (mod\ q_{_{i}}^{\alpha_{_{i}}}) & for\ all\ 1 \leq i \leq k \\
and: & n + r_{_{q_{_{i}}^{\alpha_{_{i}}}}}(n) & \equiv & 0\ (mod\ q_{_{i}}^{\alpha_{_{i}}}) & for\ all\ 1 \leq i \leq k \\
then: & n - a & \equiv & 0\ (mod\ q_{_{i}}^{\alpha_{_{i}}}) & for\ all\ 1 \leq i \leq k \\
whence: & n & \equiv & a\ (mod\ d) & 
\end{array}$
\end{quote}

\vspace{+1ex}
\noindent The Lemma follows. \hfill $\Box$

\vspace{+1ex}
\noindent By Lemma \ref{sec:4.lem.2}, it follows that:

\begin{corollary}
\label{sec:5.lem.3}
The probability that $a + r_{_{q_{_{i}}^{\alpha_{_{i}}}}}(n) \equiv 0\ (mod\ q_{_{i}}^{\alpha_{_{i}}})$ for any $1 \leq i \leq k$ is $\frac{1}{q_{_{i}}^{\alpha_{_{i}}}}$. \hfill $\Box$
\end{corollary}

\vspace{+1ex}
\noindent By Theorem \ref{sec:5.thm.1}, it further follows that:

\begin{corollary}
\label{sec:5.lem.4}
The joint probability that $a + r_{_{q_{_{i}}^{\alpha_{_{i}}}}}(n) \equiv 0\ (mod\ q_{_{i}}^{\alpha_{_{i}}})$ for all $1 \leq i \leq k$ is $\prod_{i = 1}^{k}{\frac{1}{q_{_{i}}^{\alpha_{_{i}}}}}$. \hfill $\Box$
\end{corollary}

\vspace{+1ex}
\noindent We conclude by Lemma \ref{sec:5.lem.2} that:

\begin{corollary}
\label{sec:5.lem.4}
The probability that $n$ is of the form $a+m.d$ for some natural number $m \geq 1$, where $1 \leq a < d = q_{_{1}}^{\alpha_{_{1}}}.q_{_{2}}^{\alpha_{_{2}}} \ldots q_{_{k}}^{\alpha_{_{k}}}$ is $\prod_{i = 1}^{k}{\frac{1}{q_{_{i}}^{\alpha_{_{i}}}}}$. \hfill $\Box$
\end{corollary}

\vspace{+1ex}
\noindent It follows that:

\begin{corollary}
\label{sec:5.lem.5}
The density of Dirichlect integers, defined as numbers of the form $a+m.d$ for some natural number $m \geq 1$ which are not divisible by any given set of primes $\mathbb{R} = \{r_{_{1}}, r_{_{2}}, \ldots, r_{_{l}}\}$, where $1 \leq a < d = q_{_{1}}^{\alpha_{_{1}}}.q_{_{2}}^{\alpha_{_{2}}} \ldots q_{_{k}}^{\alpha_{_{l}}}$ is:

\vspace{+1ex}
\begin{quote}
$\prod_{i = 1}^{k}{\frac{1}{q_{_{i}}^{\alpha_{_{i}}}}}.\prod_{r \in \mathbb{R}\ \&\ r \neq q_{_{i}}}(1 - \frac{1}{r})$.
\end{quote}
\end{corollary}

\noindent \textbf{Proof}: Since $a, d$ are co-prime, we have by Lemma \ref{sec:5.lem.2} that if $n$ is of the form $a+m.d$ for some natural number $m \geq 1$, where $1 \leq a < d = q_{_{1}}^{\alpha_{_{1}}}.q_{_{2}}^{\alpha_{_{2}}} \ldots q_{_{k}}^{\alpha_{_{k}}}$, we have that:

\vspace{+1ex}
\begin{quote}
$\begin{array}{rlcll}
& n & \equiv & a\ (mod\ q_{_{i}}) & for\ all\ 1 \leq i \leq k \\
whilst: & n + r_{_{i}}(n) & \equiv & 0\ (mod\ i) & for\ all\ 1 \leq i \\
whence: & a + r_{_{q_{_{i}}}}(n) & \equiv & 0\ (mod\ q_{_{i}}) & for\ all\ 1 \leq i \leq k \\
& r_{_{q_{_{i}}}}(n) & \neq & 0 & for\ all\ 1 \leq i \leq k \\
and: & q_{_{i}} & \not | & n & for\ all\ 1 \leq i \leq k
\end{array}$
\end{quote}

\vspace{+1ex}
\noindent Hence, if $n$ is of the form $a+m.d$ for some natural number $m \geq 1$, where $1 \leq a < d = q_{_{1}}^{\alpha_{_{1}}}.q_{_{2}}^{\alpha_{_{2}}} \ldots q_{_{k}}^{\alpha_{_{k}}}$ and $(a, d) = 1$, the probability that $q_{_{i}} \not | n$ for all $1 \leq i \leq k$ is $1$.

\vspace{+1ex}
\noindent By Lemma \ref{sec:2.2.lem.0}, Theorem \ref{sec:2.2.lem.2} and Theorem \ref{sec:5.thm.1}, the density of Dirichlect numbers of the form $a+m.d$ which are not divisible by any given set of primes $\mathbb{R} = \{r_{_{1}}, r_{_{2}}, \ldots, r_{_{l}}\}$ is thus:

\vspace{+1ex}
\begin{quote}
$\prod_{i = 1}^{k}{\frac{1}{q_{_{i}}^{\alpha_{_{i}}}}}.\prod_{r \in \mathbb{R}\ \&\ r \neq q_{_{i}}}{(1 - \frac{1}{r})}$ 
\end{quote}

\vspace{+1ex}
\noindent The Corollary follows. \hfill $\Box$

\begin{corollary}
\label{sec:5.lem.6}
The expected number of Dirichlect integers in any interval $(a,b)$ is:

\vspace{+1ex}
\begin{quote}
$(b-a)\prod_{i = 1}^{k}{\frac{1}{q_{_{i}}^{\alpha_{_{i}}}}}.\prod_{i = 1}^{k}(1 - \frac{1}{q_{_{i}}})^{-1}.\prod_{r \in \mathbb{R}}(1 - \frac{1}{r})$. \hfill $\Box$
\end{quote}
\end{corollary}

\vspace{+1ex}
\subsubsection{An elementary density-based proof of Dirichlect's Theorem}
\label{dirichlect:proof}

\vspace{+2ex}
\noindent Since $n$ is a prime if, and only if, it is not divisible by any prime $p \leq \sqrt{n}$, it follows that the number $\pi_{_{(a, d)}}(n)$ of Dirichlect primes, of the form $a+m.d$ for some natural number $m \geq 1$ and $1 \leq a < d = q_{_{1}}^{\alpha_{_{1}}}.q_{_{2}}^{\alpha_{_{2}}} \ldots q_{_{k}}^{\alpha_{_{k}}}$, that are less than or equal to any $n \geq q_{_{k}}^{2}$ is cumulatively approximated by the non-heuristic Dirichlect prime counting function:

\begin{definition}
\label{sec:5.lem.6}
$\pi_{_{D}}(n) = \sum_{l = 1}^{n}(\prod_{i = 1}^{k}{\frac{1}{q_{_{i}}^{\alpha_{_{i}}}}}.\prod_{i = 1}^{k}(1 - \frac{1}{q_{_{i}}})^{-1}.\prod_{j = 1}^{\pi(\sqrt{l})}(1 - \frac{1}{p_{_{j}}})$.
\end{definition}

\noindent We conclude that:

\begin{lemma}
\label{dirichlect.lem}
 $\pi_{_{(a, d)}}(n) \approx \pi_{_{D}}(n) \rightarrow \infty$ as $n \rightarrow \infty$.
\end{lemma}

\noindent \textbf{Proof}: If $a, d$ are co-prime and $1 \leq a < d = q_{_{1}}^{\alpha_{_{1}}}.q_{_{2}}^{\alpha_{_{2}}} \ldots q_{_{k}}^{\alpha_{_{k}}}$, we have for any $n \geq q_{_{k}}^{2}$:

\vspace{+1ex}
\begin{quote}
$\begin{array}{lcl}
\pi_{_{D}}(n) & = & \sum_{l = 1}^{n}(\prod_{i = 1}^{k}{\frac{1}{q_{_{i}}^{\alpha_{_{i}}}}}.\prod_{i = 1}^{k}(1 - \frac{1}{q_{_{i}}})^{-1}.\prod_{j = 1}^{\pi(\sqrt{l})}(1 - \frac{1}{p_{_{j}}}))\\ \\

& = & \prod_{i = 1}^{k}{\frac{1}{q_{_{i}}^{\alpha_{_{i}}}}}.\prod_{i = 1}^{k}(1 - \frac{1}{q_{_{i}}})^{-1}.\sum_{l = 1}^{n}\prod_{j = 1}^{\pi(\sqrt{l})}(1 - \frac{1}{p_{_{j}}}) \\ \\

& \geq & \prod_{i = 1}^{k}{\frac{1}{q_{_{i}}^{\alpha_{_{i}}}}}.\prod_{i = 1}^{k}(1 - \frac{1}{q_{_{i}}})^{-1}.n.\prod_{j = 1}^{\pi(\sqrt{n})}(1 - \frac{1}{p_{_{j}}}) 
\end{array}$
\end{quote}

\vspace{+1ex}
\noindent The lemma follows since, by Mertens' Theorem, $\prod_{p \leq x} (1 - \frac{1}{p}) \sim \frac{e^{-\lambda}}{log_{e}x}$, we have that:

\vspace{+1ex}
\begin{quote}
$n.\prod_{j = 1}^{\pi(\sqrt{n})}(1 - \frac{1}{p_{_{j}}}) \sim \frac{2e^{-\gamma}n}{log_{e}(n)} \rightarrow \infty$ as $n \rightarrow \infty$. \hfill $\Box$
\end{quote}

\vspace{+1ex}
\noindent Since $p_{_{n+1}}^{2} - p_{_{n}}^{2} \rightarrow \infty$ as $n \rightarrow \infty$, we conclude that:

\begin{theorem}
\label{dirichlect.thm}
There are an infinity of primes in any arithmetic progression $a+m.d$ where $(a,d) = 1$\footnote{Compare \cite{HW60}, p.13, Theorem 15*.}. \hfill $\Box$
\end{theorem}

\vspace{+1ex}
\subsection{An elementary density-based proof that there are infinitely many twin-primes}
\label{sec:twin.1}

\vspace{+2ex}
\noindent We define $\pi_{_{2}}(n)$ as the number of integers $p \leq n$ such that both $p$ and $p+2$ are prime.

\vspace{+1ex}
\noindent In order to estimate $\pi_{_{2}}(n)$, we first define:

\begin{definition}
\label{sec:twin.def.2}
An integer $n$ is a $\mathbb{TW}(k)$ integer if, and only if, $r_{p_{_{i}}}(n) \neq 0$ and $r_{p_{_{i}}}(n) \neq 2$ for all $1 \leq i \leq k$, where $0 \leq r_{_{i}}(n) < i$ is defined for all $i > 1$ by:

\vspace{+1ex}
\begin{quote}
$n + r_{_{i}}(n) \equiv 0\ (mod\ i)$ .
\end{quote}
\end{definition}

\noindent We note that:

\begin{lemma}
\label{sec:twin.lem.1}
If $n$ is a $\mathbb{TW}(k)$ integer, then both $n$ and $n+2$ are not divisible by any of the first $k$ primes $\{p_{_{1}}, p_{_{2}}, \ldots, p_{_{k}}\}$.
\end{lemma}

\noindent \textbf{Proof}: The lemma follows immediately from Definition \ref{sec:twin.def.2}, Definition \ref{sec:4.lem.1} and Lemma \ref{sec:4.lem.1.1}. \hfill $\Box$

\vspace{+1ex}
\noindent Since each residue $r_{i}(n)$ cycles over the $i$ values $(i-1, i-2, \ldots, 0)$, these values are all incongruent and form a complete system of residues $mod\ i$.

\vspace{+1ex}
\noindent It thus follows from Definition \ref{sec:twin.def.2} that the density of $\mathbb{TW}(k)$ integers over the set of natural numbers is:

\begin{lemma}
\label{sec:twin.lem.3}
$\mathbb{D}(\mathbb{TW}(k)) = \prod_{i=2}^{k}(1 - \frac{2}{p_{_{i}}})$. \hfill $\Box$
\end{lemma}

\vspace{+1ex}
\noindent We also have that:

\begin{lemma}
\label{sec:twin.lem.2}
If $p_{_{k}}^{2} \leq n \leq p_{_{k+1}}^{2}$ is a $\mathbb{TW}(k)$ integer, then $n$ is a prime and either $n+2$ is also a prime, or $n+2 = p_{_{k +1}}^{2}$.
\end{lemma}

\noindent \textbf{Proof}: By Definition \ref{sec:twin.def.2} and Definition \ref{sec:4.lem.1}:

\vspace{+1ex}
$\begin{array}{lcl}
r_{p_{_{i}}}(n) & \neq & 2\ for\ all\ 1 \leq i \leq k \\

n + 2 & \neq & \lambda.p_{_{i}}\ for\ all\ 1 \leq i \leq k ,\ \lambda \geq 1
\end{array}$

\vspace{+1ex}
\noindent Hence $n$ is prime; and either $n+2$ is divisible by $p_{_{k+1}}$, in which case $n+2 = p_{_{k+1}}^{2}$, or it is a prime. \hfill $\Box$

\vspace{+1ex}
\noindent If we define $\pi_{_{\mathbb{TW}(k)}}(n)$ as the number of $\mathbb{TW}(k)$ integers $\leq n$, by Lemma \ref{sec:twin.lem.3} the expected number of $\mathbb{TW}(k)$ integers in any interval $(a, b)$ is given---with a binomial standard deviation---by:

\begin{lemma}
\label{sec:twin.lem.4}
$\pi_{_{\mathbb{TW}(k)}}(b) - \pi_{_{\mathbb{TW}(k)}}(a) \approx (b-a)\prod_{i=2}^{k}(1 - \frac{2}{p_{_{i}}})$. \hfill $\Box$
\end{lemma}

\vspace{+1ex}
\noindent Since $n$ is a prime if, and only if, it is not divisible by any prime $p \leq \sqrt{n}$, it follows from Lemma \ref{sec:twin.lem.2} that $\pi_{_{\mathbb{TW}(k)}}(p_{_{k+1}}^{2}) - \pi_{_{\mathbb{TW}(k)}}(p_{_{k}}^{2})$ is at most one less than the number of twin-primes in the interval $(p_{_{k+1}}^{2} - p_{_{k}}^{2})$.

\begin{lemma}
\label{sec.twin.lem.0}
$\pi_{_{\mathbb{TW}(k)}}(p_{_{k+1}}^{2}) - \pi_{_{\mathbb{TW}(k)}}(p_{_{k}}^{2}) + 1 \geq \pi_{_{2}}(p_{_{k+1}}^{2}) - \pi_{_{2}}(p_{_{k)}}^{2}) \geq \pi_{_{\mathbb{TW}(k)}}(p_{_{k+1}}^{2}) - \pi_{_{\mathbb{TW}(k)}}(p_{_{k}}^{2})$
\end{lemma} 

\vspace{+1ex}
\noindent Now, by Lemma \ref{sec:twin.lem.4} the expected number of $\mathbb{TW}(k)$ integers in the interval $(p_{_{k+1}}^{2} - p_{_{k}}^{2})$ is given by:

\begin{lemma}
\label{sec:twin.lem.4.1}
$\pi_{_{\mathbb{TW}(k)}}(p_{_{k+1}}^{2}) - \pi_{_{\mathbb{TW}(k)}}(p_{_{k}}^{2}) \approx (p_{_{k+1}}^{2} - p_{_{k}}^{2})\prod_{i=2}^{k}(1 - \frac{2}{p_{_{i}}})$. \hfill $\Box$
\end{lemma}

\vspace{+1ex}
\noindent We conclude that the number $\pi_{_{2}}$ of twin primes $\leq p_{_{k+1}}^{2}$ is given by the cumulative non-heuristic approximation: 

\begin{lemma}
\label{sec:twin.lem.4.1.1}
$\sum_{j=1}^{k} (\pi_{_{2}}(p_{_{j+1}}^{2}) - \pi_{_{2}}(p_{_{j}}^{2})) = \pi_{_{2}}(p_{_{k+1}}^{2}) \approx \sum_{j=1}^{k} (p_{_{j+1}}^{2} - p_{_{j}}^{2})\prod_{i=2}^{j}(1 - \frac{2}{p_{_{i}}})$. \hfill $\Box$
\end{lemma}

\noindent We further conclude that:

\begin{theorem}
\label{sec:twin.thm.1}
$\pi_{_{2}}(n) \rightarrow \infty$ as $n \rightarrow \infty$.
\end{theorem}

\noindent \textbf{Proof}: We have that, for $k \geq 2$:

\vspace{+1ex}
$\begin{array}{rcl}

\sum_{j=1}^{k} (p_{_{j+1}}^{2} - p_{_{j}}^{2})\prod_{i=2}^{j}(1 - \frac{2}{p_{_{i}}}) & = & \sum_{j=9}^{p_{_{k+1}}^{2}} \prod_{i=2}^{\pi(\sqrt{j})-1}(1 - \frac{2}{p_{_{i}}}) \\ \\

& \geq & (p_{_{k+1}}^{2}-9).\prod_{i=2}^{k}(1 - \frac{2}{p_{_{i}}}) \\ \\

& \geq & (p_{_{k+1}}^{2}-9).\prod_{i=2}^{k}(1 - \frac{1}{p_{_{i}}})(1 - \frac{1}{(p_{_{i}}-1)}) \\ \\

& \geq & (p_{_{k+1}}^{2}-9).\prod_{i=2}^{k}(1 - \frac{1}{p_{_{i}}})(1 - \frac{1}{p_{_{i-1}}}) \\ \\

& \geq & (p_{_{k+1}}^{2}-9).\prod_{i=2}^{k}(1 - \frac{1}{p_{_{i-1}}})^{2} \\ \\

& \geq & (p_{_{k+1}}^{2}-9).\prod_{i=1}^{k}(1 - \frac{1}{p_{_{i}}})^{2} \\ \\
\end{array}$

\vspace{+1ex}
\noindent Now, by Mertens' Theorem, we have that:

\vspace{+1ex}
$\begin{array}{lcl}
(p_{_{k+1}}^{2}-9).\prod_{i=1}^{k}(1 - \frac{1}{p_{_{i}}})^{2} & \sim & (p_{_{k+1}}^{2}-9).(\frac{e^{-\gamma}}{log_{e}k})^{2} \\ \\

& \rightarrow & \infty\ as\ n \rightarrow \infty
\end{array}$

\vspace{+1ex}
\noindent The theorem follows by Lemma \ref{sec:twin.lem.4.1.1}. \hfill $\Box$

\subsection{The Generalised Prime Counting Function: $\sum_{j = 1}^{n} \prod_{i = a}^{\pi(\sqrt{j})}(1 - \frac{b}{p_{_{i}}})$}
\label{limit.infinite.product}

\vspace{+2ex}
\noindent We note that the argument of Theorem \ref{sec:twin.thm.1} in \S \ref{sec:twin.1} is a special case of the limiting behaviour of the Generalised Prime Counting Function $\sum_{j = 1}^{n} \prod_{i = a}^{\pi(\sqrt{j})}(1 - \frac{b}{p_{_{i}}})$, which estimates the number of integers $\leq n$ such that there are $b$ values that cannot occur amongst the residues $r_{p_{_{i}}}(n)$ for $a \leq i \leq \pi(\sqrt{j})$\footnote{Thus $b = 1$ yields an estimate for the number of primes $\leq n$, and $b = 2$ an estimate for the number of TW primes (Definition \ref{sec:twin.def.2}) $\leq n$.}:

\begin{theorem}
\label{thm.infinite.product}
$\sum_{j = 1}^{n} \prod_{i = a}^{\pi(\sqrt{j})}(1 - \frac{b}{p_{_{i}}}) \rightarrow \infty$ as $n \rightarrow \infty$ if $p_{_{a}} > b \geq 1$. 
\end{theorem}

\noindent \textbf{Proof}: For $p_{_{a}} > b \geq 1$, we have that:

\vspace{+2ex}
$\begin{array}{rcl}
\sum_{j = 1}^{n} \prod_{i = a}^{\pi(\sqrt{j})}(1 - \frac{b}{p_{_{i}}}) & \geq & \sum_{j = p_{_{a}}^{2}}^{n} \prod_{i = a}^{\pi(\sqrt{j})}(1 - \frac{b}{p_{_{i}}}) \\ \\
& \geq & \sum_{j = p_{_{a}}^{2}}^{n} \prod_{i = a}^{\pi(\sqrt{n})}(1 - \frac{b}{p_{_{i}}}) \\ \\
& \geq & (n - p_{_{a}}^{2}).\prod_{i = a}^{\pi(\sqrt{n})}(1 - \frac{b}{p_{_{i}}}) \\ \\
& \geq & (n - p_{_{a}}^{2}).\prod_{i = a}^{n}(1 - \frac{b}{p_{_{i}}})
\end{array}$

\vspace{+2ex}
\noindent The theorem follows if:

\vspace{+2ex}
$log_{_{e}}(n - p_{_{a}}^{2}) + \sum_{i = a}^{n}log_{_{e}}(1 - \frac{b}{p_{_{i}}}) \rightarrow \infty$

\vspace{+1ex}
\begin{quote}
\noindent (i) We note first the standard result for $|x| < 1$ that:

\vspace{+1ex}
\begin{quote}
$log_{_{e}}(1 - x) = - \sum_{m = 1}^{\infty}\frac{x^{m}}{m}$
\end{quote}

\vspace{+1ex}
\noindent For any $p_{_{i}} > b \geq 1$, we thus have:

\vspace{+1ex}
\begin{quote}
$log_{_{e}}(1 - \frac{b}{p_{_{i}}}) = - \sum_{m = 1}^{\infty}\frac{(b/p_{_{i}})^{m}}{m} = -\frac{b}{p_{_{i}}} - \sum_{m = 2}^{\infty}\frac{(b/p_{_{i}})^{m}}{m}$
\end{quote}

\vspace{+1ex}
\noindent Hence:

\vspace{+1ex}
\begin{quote}
$\sum_{i = a}^{n}log_{_{e}}(1 - \frac{b}{p_{_{i}}}) = -\sum_{i = a}^{n}(\frac{b}{p_{_{i}}}) - \sum_{i = a}^{n}(\sum_{m = 2}^{\infty}\frac{(b/p_{_{i}})^{m}}{m})$
\end{quote}

\vspace{+1ex}
\noindent (ii) We note next that, for all $i \geq a$:

\vspace{+1ex}
\begin{quote}
$c < (1 - \frac{b}{p_{_{a}}}) \rightarrow c < (1 - \frac{b}{p_{_{i}}})$
\end{quote}

\vspace{+1ex}
\noindent It follows for any such $c$ that:

\vspace{+1ex}
\begin{quote}
$\sum_{m = 2}^{\infty}\frac{(b/p_{_{i}})^{m}}{m} \leq \sum_{m = 2}^{\infty}(\frac{b}{p_{_{i}}})^{m} = \frac{(b/p_{_{i}})^{2}}{1 - b/p_{_{i}}} \leq \frac{b^{2}}{c.p_{_{i}}^{2}}$
\end{quote}

\vspace{+1ex}
\noindent Since:

\vspace{+1ex}
\begin{quote}
$\sum_{i = 1}^{\infty}\frac{1}{p_{_{i}}^{2}} = O(1)$
\end{quote}

\vspace{+1ex}
\noindent it further follows that:

\vspace{+1ex}
\begin{quote}
$\sum_{i = a}^{n}(\sum_{m = 2}^{\infty}\frac{(b/p_{_{i}})^{m}}{m}) \leq \sum_{i = a}^{n}(\frac{b^{2}}{c.p_{_{i}}^{2}}) = O(1)$
\end{quote}

\vspace{+1ex}
\noindent (iii) From the standard result\footnote{\cite{HW60}, p.351, Theorem 427.}:

\vspace{+1ex}
\begin{quote}
$\sum_{p \leq x}\frac{1}{p} = log_{e}log_{e}x + O(1) + o(1)$
\end{quote}

\vspace{+1ex}
\noindent it then follows that:

\vspace{+1ex}
\begin{quote}
$\begin{array}{rcl}
\sum_{i = a}^{n}log_{_{e}}(1 - \frac{b}{p_{_{i}}}) & \geq & -\sum_{i = a}^{n}(\frac{b}{p_{_{i}}}) - O(1) \\ \\

& \geq & - b.(log_{e}log_{e}n + O(1) + o(1)) - O(1)
\end{array}$
\end{quote}
\end{quote}

\vspace{+1ex}
\noindent The theorem follows since:

\vspace{+1ex}
$log_{_{e}}(n - p_{_{a}}^{2}) - b.(log_{e}log_{e}n + O(1) + o(1)) - O(1) \rightarrow \infty$

\vspace{+1ex}
\noindent and so:

\vspace{+1ex}
$log_{_{e}}(n - p_{_{a}}^{2}) + \sum_{i = a}^{n}log_{_{e}}(1 - \frac{b}{p_{_{i}}}) \rightarrow \infty$ \hfill $\Box$

\section{Appendix I: An anomaly}
\label{intro.1.4}

\begin{center}
\line(1,0){205}

\vspace{+.5ex}
\textbf{Fig.6: The graph of $y = \prod_{i = 1}^{\pi(\sqrt{x})}(1 - \frac{1}{p_{_{i}}})$}

\line(1,0){205}

\vspace{+4ex}
\begin{picture}(200,165)
\put(0,0){\line(1,0){245}}
\put(0,0){\line(0,1){180}}

\put(-17,2){\footnotesize{y $\uparrow$}}

\put(-2,80){\line(1,0){2}}
\put(-2,100){\line(1,0){2}}
\put(-2,120){\line(1,0){2}}
\put(-2,175){\line(1,0){2}}

\put(-15,76){\footnotesize{$\frac{8}{35}$}}
\put(-15,97){\footnotesize{$\frac{4}{15}$}}
\put(-15,117){\footnotesize{$\frac{1}{3}$}}
\put(-15,172){\footnotesize{$\frac{1}{2}$}}

\put(-22,-15){\footnotesize{x $\rightarrow$}}

\put(4,-4){\line(0,1){4}}
\put(8,-4){\line(0,1){4}}
\put(18,-4){\line(0,1){4}}
\put(50,-4){\line(0,1){4}}
\put(98,-4){\line(0,1){4}}
\put(242,-4){\line(0,1){4}}

\put(2,-15){\tiny{$1$}}
\put(6,-15){\tiny{$4$}}
\put(16,-15){\tiny{$9$}}
\put(48,-15){\tiny{$25$}}
\put(96,-15){\tiny{$49$}}
\put(240,-15){\tiny{$121$}}

\put(190,-10){\tiny{\textit{Not to scale}}}

\put(4,120){\framebox(14,55){}}
\put(20,175){\tiny $A$}
\put(7,40){\tiny $\frac{4}{3}$}
\qbezier[60](18,0)(18,80)(18,120)

\put(4,100){\framebox(46,20)}{}
\put(52,120){\tiny $B$}
\put(30,40){\tiny $\frac{9}{6.7}$}
\qbezier[50](50,0)(50,90)(50,100)

\put(4,80){\framebox(94,20){}}
\put(100,100){\tiny $C$}
\put(68,40){\tiny $\frac{15}{11.2}$}
\qbezier[40](98,0)(98,40)(98,80)

\put(4,0){\framebox(238,80){}}
\put(244,80){\tiny $D$}
\put(140,40){\tiny $\frac{\pi(11^{2}) = 30}{\pi_{_{H}}(11^2) =25.1}$}
\put(242,0){\line(0,1){4}}
\end{picture}
\vspace{+3ex}

\begin{quote}
{\scriptsize Fig.6: Graph of $y = \prod_{i = 1}^{\pi(\sqrt{x})}(1 - \frac{1}{p_{_{i}}})$\footnote{Compare with the graph of the same function, $y = \prod_{i = 1}^{\pi(\sqrt{x})}(1 - \frac{1}{p_{_{i}}})$, in \S \ref{intro.1.3}, Fig.4.}. The overlapping rectangles $A, B, C, D, \ldots$ represent $\pi_{_{H}}(p_{_{j+1}}^{2}) = p_{_{j+1}}^{2}.\prod_{i = 1}^{j}(1 - \frac{1}{p_{_{i}}})$ for $j \geq 1$. Figures within each rectangle are the primes and estimated primes corresponding to the functions $\pi(n)$ and $\pi_{_{H}}(n)$, respectively, within the interval $(1,\ p_{_{j+1}}^{2})$ for $j \geq 2$}.
\end{quote}
\end{center}

\footnotesize
\noindent Although we have shown in this investigation that both $\pi_{_{H}}(n)$ (Lemma \ref{sec:2.5.lem.H.1}) and $\pi_{_{L}}(n)$ (Lemma \ref{sec:2.5.lem.1}) yield non-heuristic approximations for $\pi(n)$\footnote{Fig.12 in \S \ref{appendix.prim.dist.1500} details the values of $\pi_{_{L}}(n), \pi(n)$ and $\pi_{_{H}}(n)$ for $4 \leq n \leq 1500$.}, Fig.6 suggests anomalously that, prima facie, $\pi(p_{_{n}}^{2}) > \pi_{_{H}}(p_{_{n}}^{2})$ for all $n > 1$, since:

\begin{quote}
$\bullet$ The number of primes in any interval $(p_{_{j}}^{2},\ p_{_{j+1}}^{2})$ for any given $j < n$ is constant as $n \rightarrow \infty$; but

\vspace{+1ex}
$\bullet$ The contribution of the expected number of primes in the interval $(p_{_{j}}^{2},\ p_{_{j+1}}^{2})$ to the total expected number (see Corollary \ref{sec:2.5.cor.2.H}), $\pi_{_{H}}(p_{_{n+1}}^{2}) = p_{_{n+1}}^{2}.\prod_{i = 1}^{j}(1 - \frac{1}{p_{_{i}}})$, of primes in the interval $(1,\ p_{_{n+1}}^{2})$ decreases monotonically.
\end{quote} 

\noindent In other words, an apparent anomaly surfaces when we express $\pi(n)$ and the function $\pi_{_{H}}(n)$ in terms of the number of primes determined by each function respectively in each interval $(p_{_{n}}^{2},\ p_{_{n+1}}^{2})$ as follows:

\vspace{+1ex}
$\begin{array}{lcl}
\pi(p_{_{n+1}}^{2}) & = & \sum_{j=1}^{n}(\pi(p_{_{j+1}}^{2}) - \pi(p_{_{j}}^{2})) + \pi(p_{_{1}}^{2}) \\ \\

\pi_{_{H}}(p_{_{n+1}}^{2}) & = & p_{_{n+1}}^{2}.\prod_{i = 1}^{\pi(\sqrt{p_{_{n+1}}^{2}})}(1 - \frac{1}{p_{_{i}}}) \\ \\

& = & (\sum_{j=1}^{n}(p_{_{j+1}}^{2} - p_{_{j}}^{2}) + p_{_{1}}^{2}).\prod_{i = 1}^{n}(1 - \frac{1}{p_{_{i}}}) \\ \\

& = & \sum_{j=1}^{n}(p_{_{j+1}}^{2}.\prod_{i = 1}^{n}(1 - \frac{1}{p_{_{i}}}) - p_{_{j}}^{2}.\prod_{i = 1}^{n}(1 - \frac{1}{p_{_{i}}})) + p_{_{1}}^{2}.\prod_{i = 1}^{n}(1 - \frac{1}{p_{_{i}}})
\end{array}$

\vspace{+1.5ex}
\noindent Now, by Lemma \ref{sec:2.5.lem.1}, $\pi(n) \approx \pi_{_{L}}(n)$. Hence, for any given $k > 1$:

\vspace{+1.5ex}
$\begin{array}{cll}
\pi(p_{_{k+1}}^{2}) - \pi(p_{_{k}}^{2}) & \approx & \pi_{_{L}}(p_{_{k+1}}^{2}) - \pi_{_{L}}(p_{_{k}}^{2})\ as\ k \rightarrow \infty \\ \\

\pi_{_{L}}(p_{_{k+1}}^{2}) - \pi_{_{L}}(p_{_{k}}^{2}) & = & (p_{_{k+1}}^{2} - p_{_{k}}^{2}).\prod_{i = 1}^{k}(1 - \frac{1}{p_{_{i}}}) \\

& \geq & ((p_{_{k}}+2)^{2} - p_{_{k}}^{2}).\prod_{i = 1}^{k}(1 - \frac{1}{p_{_{i}}}) \\

& \geq & 4(p_{_{k}}+1).\prod_{i = 1}^{k}(1 - \frac{1}{p_{_{i}}}) \\

& \in & O(\frac{p_{_{k}}}{log_{e}p_{_{k}}})\ as\ k \rightarrow \infty\\

& \rightarrow & \infty\ as\ k \rightarrow \infty
\end{array}$

\vspace{+1ex}
\noindent whilst, for any given $k > 1$:

\vspace{+1ex}
\begin{quote}
$p_{_{k+1}}^{2}.\prod_{i = 1}^{n}(1 - \frac{1}{p_{_{i}}}) - p_{_{k}}^{2}.\prod_{i = 1}^{n}(1 - \frac{1}{p_{_{i}}}) \rightarrow 0$ as $n \rightarrow \infty$.
\end{quote}

\noindent We conclude that\footnote{Fig.12 in \S \ref{appendix.prim.dist.1500} compares the values of $\pi_{_{L}}(n), \pi(n)$ and $\pi_{_{H}}(n)$ for $4 \leq n \leq 1500$.}:

\begin{lemma}
\label{anomaly.lem}
$\pi_{_{L}}(n)$ is a better approximation of $\pi(n)$ than $\pi_{_{H}}(n)$ for all $n \geq 9$.
\end{lemma}

\vspace{+1ex}
\noindent However, prima facie, we then have the anomaly:

\begin{anomaly}
\label{anomaly}
$\pi(n) \approx \pi_{_{L}}(n) > \pi_{_{H}}(n) \sim 2e^{-\lambda}\frac{n}{log_{_{e}}n}$.
\end{anomaly}

\section{Appendix II: The residue function $r_{i}(n)$}
\label{appendix}

\noindent We graphically illustrate how the residues $r_{i}(n)$ occur naturally as values of:

\footnotesize
\vspace{+1ex}
A: The natural-number based residue functions $R_{i}(n)$;

\vspace{+1ex}
B: The natural-number based residue sequences $E(n)$;

\vspace{+1ex}
\noindent and as the output of:

\vspace{+1ex}
C: The natural-number based algorithm $E_{\mathbb{N}}$;

\vspace{+1ex}
D: The prime-number based algorithm  $E_{\mathbb{P}}$;

\vspace{+1ex}
E: The prime-number based algorithm  $E_{\mathbb{Q}}$.

\subsection*{A: The natural-number based density-defining functions $R_{i}(n)$}
\label{sec:4.fig1}

\footnotesize
\vspace{+1ex}
\noindent \textbf{Density:} For instance, the residues $r_{i}(n)$ can be defined for all $n \geq 1$ as the values of the density-defining functions $R_{i}(n)$, defined for all $i \geq 1$, as illustrated below in Fig.7, where:

\begin{quote}
\noindent $\bullet$ Each function $R_{_{i}}(n)$ cycles through the values $(i-1,\ i-2,\ \ldots,\ 0)$ with period $i$;

\vspace{+1ex}
\noindent $\bullet$ For any $i \geq 2$ the density---over the set of natural numbers---of the set $\{n\}$ of integers that are divisible by $i$ is $\frac{1}{i}$; and the density of integers that are not divisible by $i$ is $\frac{i - 1}{i}$.
\end{quote}

\vspace{+1ex}
\noindent \line(1,0){285} \\
\noindent \textbf{\small{Fig.7: The natural-number based residue functions $R_{i}(n)$}} \\
\noindent \line(1,0){285}

\tiny
\begin{tabbing}
Divisors: \= 001 \= 002 \= 003 \= 004 \= 005 \= 006 \= 007 \= 008 \= 009 \= 0010 \= 0011 \= \ldots \= 0n-9 \ldots \kill
Function: \> $R_{_{1}}$ \> $R_{_{2}}$ \> $R_{_{3}}$ \> $R_{_{4}}$ \> $R_{_{5}}$ \> $R_{_{6}}$ \> $R_{_{7}}$ \> $R_{_{8}}$ \> $R_{_{9}}$ \> $R_{_{10}}$ \> $R_{_{11}}$ \> \ldots \> $R_{_{n}}$ \\
\rule{75mm}{.5mm} \\
$\textcolor{black}{n = 1}$ \> \textcolor{black}{0} \> 1 \> 2 \> 3 \> 4 \> 5 \> 6 \> 7 \> 8 \> 9 \> 10 \> \ldots \> n-1 \\
$\textcolor{black}{n = 2}$ \> \textcolor{black}{0} \> \textcolor{black}{0} \> 1 \> 2 \> 3 \> 4 \> 5 \> 6 \> 7 \> 8 \> 9 \> \ldots \> n-2 \\
$\textcolor{black}{n = 3}$ \> \textcolor{black}{0} \> \textcolor{black}{1} \> \textcolor{black}{0} \> 1 \> 2 \> 3 \> 4 \> 5 \> 6 \> 7 \> 8 \> \ldots \> n-3 \\
$\textcolor{black}{n = 4}$ \> \textcolor{black}{0} \> 0 \> 2 \> \textcolor{black}{0} \> 1 \> 2 \> 3 \> 4 \> 5 \> 6 \> 7 \> \ldots \> n-4 \\
$\textcolor{black}{n = 5}$ \> \textcolor{black}{0} \> \textcolor{black}{1} \> \textcolor{black}{1} \> \textcolor{black}{3} \> \textcolor{black}{0} \> 1 \> 2 \> 3 \> 4 \> 5 \> 6 \> \ldots \> n-5 \\
$\textcolor{black}{n = 6}$ \> \textcolor{black}{0} \> 0 \> 0 \> 2 \> 4 \> \textcolor{black}{0} \> 1 \> 2 \> 3 \> 4 \> 5 \> \ldots \> n-6 \\
$\textcolor{black}{n = 7}$ \> \textcolor{black}{0} \> \textcolor{black}{1} \> \textcolor{black}{2} \> \textcolor{black}{1} \> \textcolor{black}{3} \> \textcolor{black}{5} \> \textcolor{black}{0} \> 1 \> 2 \> 3 \> 4 \> \ldots \> n-7 \\
$\textcolor{black}{n = 8}$ \> \textcolor{black}{0} \> 0 \> 1 \> 0 \> 2 \> 4 \> 6 \> \textcolor{black}{0} \> 1 \> 2 \> 3 \> \ldots \> n-8 \\
$\textcolor{black}{n = 9}$ \> \textcolor{black}{0} \> 1 \> 0 \> 3 \> 1 \> 3 \> 5 \> 7 \> \textcolor{black}{0} \> 1 \> 2 \> \ldots \> n-9 \\
$\textcolor{black}{n = 10}$ \> \textcolor{black}{0} \> 0 \> 2 \> 2 \> 0 \> 2 \> 4 \> 6 \> 8 \> \textcolor{black}{0} \> 1 \> \ldots \> n-10 \\
$\textcolor{black}{n = 11}$ \> \textcolor{black}{0} \> \textcolor{black}{1} \> \textcolor{black}{1} \> \textcolor{black}{1} \> \textcolor{black}{4} \> \textcolor{black}{1} \> \textcolor{black}{3} \> \textcolor{black}{5} \> \textcolor{black}{7} \> \textcolor{black}{9} \> \textcolor{black}{0} \> \ldots \> n-11 \\
\\
$\textcolor{black}{n}$ \> $\textcolor{black}{r_{_{1}}}$ \> $r_{_{2}}$ \> $r_{_{3}}$ \> $r_{_{4}}$ \> $r_{_{5}}$ \> $r_{_{6}}$ \> $r_{_{7}}$ \> $r_{_{8}}$ \> $r_{_{9}}$ \> $r_{_{10}}$ \> $r_{_{11}}$ \> \ldots \> \textcolor{black}{0} \\
\end{tabbing}

\noindent \tiny{Fig.7: The natural-number based residue functions $R_{i}(n)$}

\subsection*{B: The natural-number based primality-defining sequences $E(n)$}

\footnotesize
\vspace{+1ex}
\noindent \textbf{Primality:} The residues $r_{i}(n)$ can also be viewed alternatively as values of the associated primality-defining sequences, $E(n) = \{r_{i}(n): i \geq 1\}$, defined for all $n \geq 1$, as illustrated below in Fig.8, where:

\begin{quote}
\noindent $\bullet$ The sequences $E(n)$ highlighted in red correspond to a prime\footnote{Conventionally defined as integers that are not divisible by any smaller integer other than $1$.} $p$ (since $r_{i}(p) \neq 0$ for $1 < i < p$) in the usual, linearly displayed, Eratosthenes sieve:

\vspace{+1ex}
$\textcolor{cyan}{E(\cancel{1})},\ \textcolor{red}{E(2)},\ \textcolor{red}{E(3)},\ \textcolor{cyan}{E(\cancel{4})},\ \textcolor{red}{E(5)},\ \textcolor{cyan}{E(\cancel{6})},\ \textcolor{red}{E(7)},\ \textcolor{cyan}{E(\cancel{8})},\ \textcolor{cyan}{E(\cancel{9})},\ \textcolor{cyan}{E(\cancel{10})},\ \textcolor{red}{E(11)},\ \ldots$

\vspace{+1ex}

\noindent $\bullet$ The sequences highlighted in cyan identify a crossed out composite $n$ (since $r_{i}(n) = 0$ for some $i < i < n$) in the usual, linearly displayed, Eratosthenes sieve. 

\vspace{+1ex}
\noindent $\bullet$ The `boundary' residues $r_{1}(n) = 0$ and $r_{n}(n) = 0$ are identified in cyan. 
\end{quote}

\noindent \line(1,0){285} \\
\noindent \textbf{\small{Fig.8: The natural-number based residue sequences $E(n)$}} \\
\noindent \line(1,0){285}

\tiny
\begin{tabbing}
Divisors: \= 001 \= 002 \= 003 \= 004 \= 005 \= 006 \= 007 \= 008 \= 009 \= 0010 \= 0011 \= \ldots \= 0n-9 \ldots \kill

Function: \> $R_{_{1}}$ \> $R_{_{2}}$ \> $R_{_{3}}$ \> $R_{_{4}}$ \> $R_{_{5}}$ \> $R_{_{6}}$ \> $R_{_{7}}$ \> $R_{_{8}}$ \> $R_{_{9}}$ \> $R_{_{10}}$ \> $R_{_{11}}$ \> \ldots \> $R_{_{n}}$ \\
\rule{75mm}{.5mm} \\
$\textcolor{cyan}{E(1)}$: \> \textcolor{cyan}{0} \> 1 \> 2 \> 3 \> 4 \> 5 \> 6 \> 7 \> 8 \> 9 \> 10 \> \ldots \> n-1 \\
$\textcolor{red}{E(2)}$: \> \textcolor{cyan}{0} \> \textcolor{cyan}{0} \> 1 \> 2 \> 3 \> 4 \> 5 \> 6 \> 7 \> 8 \> 9 \> \ldots \> n-2 \\
$\textcolor{red}{E(3)}$: \> \textcolor{cyan}{0} \> \textcolor{red}{1} \> \textcolor{cyan}{0} \> 1 \> 2 \> 3 \> 4 \> 5 \> 6 \> 7 \> 8 \> \ldots \> n-3 \\
$\textcolor{cyan}{E(4)}$: \> \textcolor{cyan}{0} \> \textcolor{cyan}{0} \> 2 \> \textcolor{cyan}{0} \> 1 \> 2 \> 3 \> 4 \> 5 \> 6 \> 7 \> \ldots \> n-4 \\
$\textcolor{red}{E(5)}$: \> \textcolor{cyan}{0} \> \textcolor{red}{1} \> \textcolor{red}{1} \> \textcolor{red}{3} \> \textcolor{cyan}{0} \> 1 \> 2 \> 3 \> 4 \> 5 \> 6 \> \ldots \> n-5 \\
$\textcolor{cyan}{E(6)}$: \> \textcolor{cyan}{0} \> \textcolor{cyan}{0} \> \textcolor{cyan}{0} \> 2 \> 4 \> \textcolor{cyan}{0} \> 1 \> 2 \> 3 \> 4 \> 5 \> \ldots \> n-6 \\
$\textcolor{red}{E(7)}$: \> \textcolor{cyan}{0} \> \textcolor{red}{1} \> \textcolor{red}{2} \> \textcolor{red}{1} \> \textcolor{red}{3} \> \textcolor{red}{5} \> \textcolor{cyan}{0} \> 1 \> 2 \> 3 \> 4 \> \ldots \> n-7 \\
$\textcolor{cyan}{E(8)}$: \> \textcolor{cyan}{0} \> \textcolor{cyan}{0} \> 1 \> \textcolor{cyan}{0} \> 2 \> 4 \> 6 \> \textcolor{cyan}{0} \> 1 \> 2 \> 3 \> \ldots \> n-8 \\
$\textcolor{cyan}{E(9)}$: \> \textcolor{cyan}{0} \> 1 \> \textcolor{cyan}{0} \> 3 \> 1 \> 3 \> 5 \> 7 \> \textcolor{cyan}{0} \> 1 \> 2 \> \ldots \> n-9 \\
$\textcolor{cyan}{E(10)}$: \> \textcolor{cyan}{0} \> \textcolor{cyan}{0} \> 2 \> 2 \> \textcolor{cyan}{0} \> 2 \> 4 \> 6 \> 8 \> \textcolor{cyan}{0} \> 1 \> \ldots \> n-10 \\
$\textcolor{red}{E(11)}$: \> \textcolor{cyan}{0} \> \textcolor{red}{1} \> \textcolor{red}{1} \> \textcolor{red}{1} \> \textcolor{red}{4} \> \textcolor{red}{1} \> \textcolor{red}{3} \> \textcolor{red}{5} \> \textcolor{red}{7} \> \textcolor{red}{9} \> \textcolor{cyan}{0} \> \ldots \> n-11 \\
\ldots \\
$E(n)$: \> $\textcolor{cyan}{r_{_{1}}}$ \> $r_{_{2}}$ \> $r_{_{3}}$ \> $r_{_{4}}$ \> $r_{_{5}}$ \> $r_{_{6}}$ \> $r_{_{7}}$ \> $r_{_{8}}$ \> $r_{_{9}}$ \> $r_{_{10}}$ \> $r_{_{11}}$ \> \ldots \> \textcolor{cyan}{0} \\
\ldots
\end{tabbing}

\noindent \tiny{Fig.8: The natural-number based residue sequences $E(n)$}

\subsection*{C: The output of a natural-number based algorithm $E_{\mathbb{N}}$}
\label{sec:2.fig}

\footnotesize
\vspace{+1ex}
\noindent We give below in Fig.9 the output for $1 \leq n \leq 11$ of a natural-number based algorithm $E_{\mathbb{N}}$ that computes the values $r_{i}(n)$ of the sequence $E_{\mathbb{N}}(n)$ for only $1 \leq i \leq n$ for any given $n$.

\vspace{+1ex}
\noindent \line(1,0){308} \\
\noindent \textbf{\small{Fig.9: The output of the natural-number based algorithm $E_{\mathbb{N}}$}} \\
\noindent \line(1,0){308}

\tiny
\begin{tabbing}
Divisors: \= 001 \= 002 \= 003 \= 004 \= 005 \= 006 \= 007 \= 008 \= 009 \= 0010 \= 0011 \= \ldots \= 0n-9 \ldots \kill
Divisors: \> 1 \> 2 \> 3 \> 4 \> 5 \> 6 \> 7 \> 8 \> 9 \> 10 \> 11 \> \ldots \> n \ldots \\
\rule{75mm}{.5mm} \\
$E_{\mathbb{N}}(1)$: \> 0 \\
$E_{\mathbb{N}}(2)$: \> 0 \> 0 \\
$E_{\mathbb{N}}(3)$: \> 0 \> 1 \> 0 \\
$E_{\mathbb{N}}(4)$: \> 0 \> 0 \> 2 \> 0 \\
$E_{\mathbb{N}}(5)$: \> 0 \> 1 \> 1 \> 3 \> 0 \\
$E_{\mathbb{N}}(6)$: \> 0 \> 0 \> 0 \> 2 \> 4 \> 0 \\
$E_{\mathbb{N}}(7)$: \> 0 \> 1 \> 2 \> 1 \> 3 \> 5 \> 0 \\
$E_{\mathbb{N}}(8)$: \> 0 \> 0 \> 1 \> 0 \> 2 \> 4 \> 6 \> 0 \\
$E_{\mathbb{N}}(9)$: \> 0 \> 1 \> 0 \> 3 \> 1 \> 3 \> 5 \> 7 \> 0 \\
$E_{\mathbb{N}}(10)$: \> 0 \> 0 \> 2 \> 2 \> 0 \> 2 \> 4 \> 6 \> 8 \> 0 \\
$E_{\mathbb{N}}(11)$: \> 0 \> 1 \> 1 \> 1 \> 4 \> 1 \> 3 \> 5 \> 7 \> 9 \> 0 \\
\ldots \\
$E_{\mathbb{N}}(n)$: \> $r_{_{1}}$ \> $r_{_{2}}$ \> $r_{_{3}}$ \> $r_{_{4}}$ \> $r_{_{5}}$ \> $r_{_{6}}$ \> $r_{_{7}}$ \> $r_{_{8}}$ \> $r_{_{9}}$ \> $r_{_{10}}$ \> $r_{_{11}}$ \> \ldots \> 0 \\
\ldots
\end{tabbing}

\noindent \tiny{Fig.9: The output of the natural-number based algorithm $\textcolor{black}{E_{\mathbb{N}}}$}

\subsection*{D: The output of the prime-number based algorithm  $E_{\mathbb{P}}$}
\label{sec:3.fig}

\footnotesize
\vspace{+1ex}
\noindent Fig.10 gives the output for $2 \leq n \leq 31$ of a prime-number based algorithm $E_{\mathbb{Q}}$ that computes the values $q_{_{i}}(n) = r_{p_{_{i}}}(n)$ of the sequence $E_{\mathbb{P}}(n)$ for only each prime $2 \leq p_{_{i}} \leq n$ for any given $n$.

\vspace{+1ex}
\noindent \line(1,0){305} \\
\noindent \textbf{\small{Fig.10: The output of the prime-number based algorithm $E_{\mathbb{P}}$}} \\
\noindent \line(1,0){305}

\tiny
\begin{tabbing}
Divisors: \= 001 \= 002 \= 003 \= 004 \= 005 \= 006 \= 007 \= 008 \= 009 \= 0010 \= 0011 \= \ldots \= 0n-9 \ldots \kill
Prime: \> $p_{_{1}}$ \> $p_{_{2}}$ \> $p_{_{3}}$ \> $p_{_{4}}$ \> $p_{_{5}}$ \> $p_{_{6}}$ \> $p_{_{7}}$ \> $p_{_{8}}$ \> $p_{_{9}}$ \> $p_{_{10}}$ \> $p_{_{11}}$ \> \ldots \> $p_{_{n}}$ \ldots \\
Divisor: \> 2 \> 3 \> 5 \> 7 \> 11 \> 13 \> 17 \> 19 \> 23 \> 29 \> 31 \> \ldots \> $p_{_{n}}$ \ldots \\
\rule{77mm}{.5mm} \\
$E_{\mathbb{P}}(2)$: \> 0 \\
$E_{\mathbb{P}}(3)$: \> 1 \> 0 \\
$E_{\mathbb{P}}(4)$: \> 0 \> 2 \\
$E_{\mathbb{P}}(5)$: \> 1 \> 1 \> 0 \\
$E_{\mathbb{P}}(6)$: \> 0 \> 0 \> 4 \\
$E_{\mathbb{P}}(7)$: \> 1 \> 2 \> 3 \> 0 \\
$E_{\mathbb{P}}(8)$: \> 0 \> 1 \> 2 \> 6 \\
$E_{\mathbb{P}}(9)$: \> 1 \> 0 \> 1 \> 5 \\
$E_{\mathbb{P}}(10)$: \> 0 \> 2 \> 0 \> 4 \\
$E_{\mathbb{P}}(11)$: \> 1 \> 1 \> 4 \> 3 \> 0 \\
$E_{\mathbb{P}}(12)$: \> 0 \> 0 \> 3 \> 2 \> 10 \\
$E_{\mathbb{P}}(13)$: \> 1 \> 2 \> 2 \> 1 \> 9 \> 0 \\
$E_{\mathbb{P}}(14)$: \> 0 \> 1 \> 1 \> 0 \> 8 \> 12 \\
$E_{\mathbb{P}}(15)$: \> 1 \> 0 \> 0 \> 6 \> 7 \> 11 \\
$E_{\mathbb{P}}(16)$: \> 0 \> 2 \> 4 \> 5 \> 6 \> 10 \\
$E_{\mathbb{P}}(17)$: \> 1 \> 1 \> 3 \> 4 \> 5 \> 9 \> 0 \\
$E_{\mathbb{P}}(18)$: \> 0 \> 0 \> 2 \> 3 \> 4 \> 8 \> 16 \\
$E_{\mathbb{P}}(19)$: \> 1 \> 2 \> 1 \> 2 \> 3 \> 7 \> 15 \> 0 \\
$E_{\mathbb{P}}(20)$: \> 0 \> 1 \> 0 \> 1 \> 2 \> 6 \> 14 \> 18 \\
$E_{\mathbb{P}}(21)$: \> 1 \> 0 \> 4 \> 0 \> 1 \> 5 \> 13 \> 17 \\
$E_{\mathbb{P}}(22)$: \> 0 \> 2 \> 3 \> 6 \> 0 \> 4 \> 12 \> 16 \\
$E_{\mathbb{P}}(23)$: \> 1 \> 1 \> 2 \> 5 \> 10 \> 3 \> 11 \> 15 \> 0 \\
$E_{\mathbb{P}}(24)$: \> 0 \> 0 \> 1 \> 4 \> 9 \> 2 \> 10 \> 14 \> 22 \\
$E_{\mathbb{P}}(25)$: \> 1 \> 2 \> 0 \> 3 \> 8 \> 1 \> 9 \> 13 \> 21 \\
$E_{\mathbb{P}}(26)$: \> 0 \> 1 \> 4 \> 2 \> 7 \> 0 \> 8 \> 12 \> 20 \\
$E_{\mathbb{P}}(27)$: \> 1 \> 0 \> 3 \> 1 \> 6 \> 12 \> 7 \> 11 \> 19 \\
$E_{\mathbb{P}}(28)$: \> 0 \> 2 \> 2 \> 0 \> 5 \> 11 \> 6 \> 10 \> 18 \\
$E_{\mathbb{P}}(29)$: \> 1 \> 1 \> 1 \> 6 \> 4 \> 10 \> 5 \> 9 \> 17 \> 0 \\
$E_{\mathbb{P}}(30)$: \> 0 \> 0 \> 0 \> 5 \> 3 \> 9 \> 4 \> 8 \> 16 \> 28 \\
$E_{\mathbb{P}}(31)$: \> 1 \> 2\> 4 \> 4 \> 2 \> 8 \> 3 \> 7 \> 15 \> 27 \> 0\\
\ldots \\
$E_{\mathbb{P}}(n)$: \> $q_{_{1}}$ \> $q_{_{2}}$ \> $q_{_{3}}$ \> $q_{_{4}}$ \> $q_{_{5}}$ \> $q_{_{6}}$ \> $q_{_{7}}$ \> $q_{_{8}}$ \> $q_{_{9}}$ \> $q_{_{10}}$ \> $q_{_{11}}$ \> \ldots \> 0 \\
\ldots
\end{tabbing}

\noindent \tiny{Fig.10: The output of the prime-number based algorithm $\textcolor{black}{E_{\mathbb{P}}}$}

\subsection*{E: The output of the prime-number based algorithms $\textcolor{cyan}{E_{\mathbb{P}}}$ and $\textcolor{black}{E_{\mathbb{Q}}}$}
\label{sec:5.fig}

\footnotesize
\vspace{+1ex}
\noindent We give below in Fig.11 the output for $2 \leq n \leq 121$ of the two prime-number based algorithms $E_{\mathbb{P}}$ (whose output $\{q_{_{i}}(n) = r_{p_{_{i}}}(n): 1 \leq i \leq \pi(n)\}$ is shown only partially, partly in cyan) and $E_{\mathbb{Q}}$ (whose output $q_{_{i}}(n) = \{r_{p_{_{i}}}(n): 1 \leq i \leq \pi(\sqrt{n})\}$ is highlighted in black and red, the latter indicating the generation of a prime sequence and, ipso facto, definition of the corresponding prime\footnote{For informal reference and perspective, formal definitions of both the  prime-number based algorithms $E_{\mathbb{P}}$ and $E_{\mathbb{Q}}$ are given in this work in progress \href{http://alixcomsi.com/40_Factorising_Update.pdf}{\textit{Factorising all $m \leq n$ is of order $\Theta(\sum_{i=2}^{n} \pi (\sqrt{i}))$}.}}.

\vspace{+1ex}
\noindent \line(1,0){349} \\
\noindent \textbf{\small{Fig.11: The output of the prime-number based algorithms $\textcolor{cyan}{E_{\mathbb{P}}}$ and $\textcolor{black}{E_{\mathbb{Q}}}$}} \\
\noindent \line(1,0){349}

\tiny
\begin{tabbing}
Divisors: \= 001 \= 002 \= 003 \= 004 \= 005 \= 006 \= 007 \= 008 \= 009 \= 0010 \= 0011 \= \ldots \= 0n-9 \ldots \kill
Prime: \> $\textcolor{black}{p_{_{1}}}$ \> $\textcolor{black}{p_{_{2}}}$ \> $\textcolor{black}{p_{_{3}}}$ \> $\textcolor{black}{p_{_{4}}}$ \> $\textcolor{black}{p_{_{5}}}$ \> $\textcolor{cyan}{p_{_{6}}}$ \> $\textcolor{cyan}{p_{_{7}}}$ \> $\textcolor{cyan}{p_{_{8}}}$ \> $\textcolor{cyan}{p_{_{9}}}$ \> $\textcolor{cyan}{p_{_{10}}}$ \> $\textcolor{cyan}{p_{_{11}}}$ \> \textcolor{cyan}{\ldots} \> $\textcolor{cyan}{p_{_{n}}}$ \textcolor{cyan}{\ldots} \\
Divisor: \> \textcolor{black}{2} \> \textcolor{black}{3} \> \textcolor{black}{5} \> \textcolor{black}{7} \> \textcolor{black}{11} \> \textcolor{cyan}{13} \> \textcolor{cyan}{17} \> \textcolor{cyan}{19} \> \textcolor{cyan}{23} \> \textcolor{cyan}{29} \> \textcolor{cyan}{31} \> \textcolor{cyan}{\ldots} \> $\textcolor{cyan}{p_{n}}$ \textcolor{cyan}{\ldots} \\
Function: \> $\textcolor{black}{Q_{_{1}}}$ \> $\textcolor{black}{Q_{_{2}}}$ \> $\textcolor{black}{Q_{_{3}}}$ \> $\textcolor{black}{Q_{_{4}}}$ \> $\textcolor{black}{Q_{_{5}}}$ \> $\textcolor{cyan}{Q_{_{6}}}$ \> $\textcolor{cyan}{Q_{_{7}}}$ \> $\textcolor{cyan}{Q_{_{8}}}$ \> $\textcolor{cyan}{Q_{_{9}}}$ \> $\textcolor{cyan}{Q_{_{10}}}$ \> $\textcolor{cyan}{Q_{_{11}}}$ \> \textcolor{cyan}{\ldots} \\
\rule{77mm}{.5mm} \\
$\textcolor{red}{E_{\mathbb{Q}}(2)}$: \> \textcolor{black}{0} \> (Prime by definition)\\
$\textcolor{red}{E_{\mathbb{Q}}(3)}$: \> \textcolor{red}{1} \> \textcolor{cyan}{0} \\
$\textcolor{black}{E_{\mathbb{Q}}(4)}$: \> \textcolor{black}{0} \> \textcolor{cyan}{2} \\
$\textcolor{red}{E_{\mathbb{Q}}(5)}$: \> \textcolor{red}{1} \> \textcolor{cyan}{1} \> \textcolor{cyan}{0} \\
$\textcolor{black}{E_{\mathbb{Q}}(6)}$: \> \textcolor{black}{0} \> \textcolor{cyan}{0} \> \textcolor{cyan}{4} \\
$\textcolor{red}{E_{\mathbb{Q}}(7)}$: \> \textcolor{red}{1} \> \textcolor{cyan}{2} \> \textcolor{cyan}{3} \> \textcolor{cyan}{0} \\
$\textcolor{black}{E_{\mathbb{Q}}(8)}$: \> \textcolor{black}{0} \> \textcolor{cyan}{1} \> \textcolor{cyan}{2} \> \textcolor{cyan}{6} \\
$\textcolor{black}{E_{\mathbb{Q}}(9)}$: \> \textcolor{black}{1} \> \textcolor{black}{0} \> \textcolor{cyan}{1} \> \textcolor{cyan}{5} \\
$\textcolor{black}{E_{\mathbb{Q}}(10)}$: \> \textcolor{black}{0} \> \textcolor{black}{2} \> \textcolor{cyan}{0} \> \textcolor{cyan}{4} \\
$\textcolor{red}{E_{\mathbb{Q}}(11)}$: \> \textcolor{red}{1} \> \textcolor{red}{1} \> \textcolor{cyan}{4} \> \textcolor{cyan}{3} \> \textcolor{cyan}{0} \\
$\textcolor{black}{E_{\mathbb{Q}}(12)}$: \> \textcolor{black}{0} \> \textcolor{black}{0} \> \textcolor{cyan}{3} \> \textcolor{cyan}{2} \> \textcolor{cyan}{10} \\
$\textcolor{red}{E_{\mathbb{Q}}(13)}$: \> \textcolor{red}{1} \> \textcolor{red}{2} \> \textcolor{cyan}{2} \> \textcolor{cyan}{1} \> \textcolor{cyan}{9} \> \textcolor{cyan}{0} \\
$\textcolor{black}{E_{\mathbb{Q}}(14)}$: \> \textcolor{black}{0} \> \textcolor{black}{1} \> \textcolor{cyan}{1} \> \textcolor{cyan}{0} \> \textcolor{cyan}{8} \> \textcolor{cyan}{12} \\
$\textcolor{black}{E_{\mathbb{Q}}(15)}$: \> \textcolor{black}{1} \> \textcolor{black}{0} \> \textcolor{cyan}{0} \> \textcolor{cyan}{6} \> \textcolor{cyan}{7} \> \textcolor{cyan}{11} \\
$\textcolor{black}{E_{\mathbb{Q}}(16)}$: \> \textcolor{black}{0} \> \textcolor{black}{2} \> \textcolor{cyan}{4} \> \textcolor{cyan}{5} \> \textcolor{cyan}{6} \> \textcolor{cyan}{10} \\
$\textcolor{red}{E_{\mathbb{Q}}(17)}$: \> \textcolor{red}{1} \> \textcolor{red}{1} \> \textcolor{cyan}{3} \> \textcolor{cyan}{4} \> \textcolor{cyan}{5} \> \textcolor{cyan}{9} \> \textcolor{cyan}{0} \\
$\textcolor{black}{E_{\mathbb{Q}}(18)}$: \> \textcolor{black}{0} \> \textcolor{black}{0} \> \textcolor{cyan}{2} \> \textcolor{cyan}{3} \> \textcolor{cyan}{4} \> \textcolor{cyan}{8} \> \textcolor{cyan}{16} \\
$\textcolor{red}{E_{\mathbb{Q}}(19)}$: \> \textcolor{red}{1} \> \textcolor{red}{2} \> \textcolor{cyan}{1} \> \textcolor{cyan}{2} \> \textcolor{cyan}{3} \> \textcolor{cyan}{7} \> \textcolor{cyan}{15} \> \textcolor{cyan}{0} \\
$\textcolor{black}{E_{\mathbb{Q}}(20)}$: \> \textcolor{black}{0} \> \textcolor{black}{1} \> \textcolor{cyan}{0} \> \textcolor{cyan}{1} \> \textcolor{cyan}{2} \> \textcolor{cyan}{6} \> \textcolor{cyan}{14} \> \textcolor{cyan}{18} \\
$\textcolor{black}{E_{\mathbb{Q}}(21)}$: \> \textcolor{black}{1} \> \textcolor{black}{0} \> \textcolor{cyan}{4} \> \textcolor{cyan}{0} \> \textcolor{cyan}{1} \> \textcolor{cyan}{5} \> \textcolor{cyan}{13} \> \textcolor{cyan}{17} \\
$\textcolor{black}{E_{\mathbb{Q}}(22)}$: \> \textcolor{black}{0} \> \textcolor{black}{2} \> \textcolor{cyan}{3} \> \textcolor{cyan}{6} \> \textcolor{cyan}{0} \> \textcolor{cyan}{4} \> \textcolor{cyan}{12} \> \textcolor{cyan}{16} \\
$\textcolor{red}{E_{\mathbb{Q}}(23)}$: \> \textcolor{red}{1} \> \textcolor{red}{1} \> \textcolor{cyan}{2} \> \textcolor{cyan}{5} \> \textcolor{cyan}{10} \> \textcolor{cyan}{3} \> \textcolor{cyan}{11} \> \textcolor{cyan}{15} \> \textcolor{cyan}{0} \\
$\textcolor{black}{E_{\mathbb{Q}}(24)}$: \> \textcolor{black}{0} \> \textcolor{black}{0} \> \textcolor{cyan}{1} \> \textcolor{cyan}{4} \> \textcolor{cyan}{9} \> \textcolor{cyan}{2} \> \textcolor{cyan}{10} \> \textcolor{cyan}{14} \> \textcolor{cyan}{22} \\
$\textcolor{black}{E_{\mathbb{Q}}(25)}$: \> \textcolor{black}{1} \> \textcolor{black}{2} \> \textcolor{black}{0} \> \textcolor{cyan}{3} \> \textcolor{cyan}{8} \> \textcolor{cyan}{1} \> \textcolor{cyan}{9} \> \textcolor{cyan}{13} \> \textcolor{cyan}{21} \\
$\textcolor{black}{E_{\mathbb{Q}}(26)}$: \> \textcolor{black}{0} \> \textcolor{black}{1} \> \textcolor{black}{4} \> \textcolor{cyan}{2} \> \textcolor{cyan}{7} \> \textcolor{cyan}{0} \> \textcolor{cyan}{8} \> \textcolor{cyan}{12} \> \textcolor{cyan}{20} \\
$\textcolor{black}{E_{\mathbb{Q}}(27)}$: \> \textcolor{black}{1} \> \textcolor{black}{0} \> \textcolor{black}{3} \> \textcolor{cyan}{1} \> \textcolor{cyan}{6} \> \textcolor{cyan}{12} \> \textcolor{cyan}{7} \> \textcolor{cyan}{11} \> \textcolor{cyan}{19} \\
$\textcolor{black}{E_{\mathbb{Q}}(28)}$: \> \textcolor{black}{0} \> \textcolor{black}{2} \> \textcolor{black}{2} \> \textcolor{cyan}{0} \> \textcolor{cyan}{5} \> \textcolor{cyan}{11} \> \textcolor{cyan}{6} \> \textcolor{cyan}{10} \> \textcolor{cyan}{18} \\
$\textcolor{red}{E_{\mathbb{Q}}(29)}$: \> \textcolor{red}{1} \> \textcolor{red}{1} \> \textcolor{red}{1} \> \textcolor{cyan}{6} \> \textcolor{cyan}{4} \> \textcolor{cyan}{10} \> \textcolor{cyan}{5} \> \textcolor{cyan}{9} \> \textcolor{cyan}{17} \> \textcolor{cyan}{0} \\
$\textcolor{black}{E_{\mathbb{Q}}(30)}$: \> \textcolor{black}{0} \> \textcolor{black}{0} \> \textcolor{black}{0} \> \textcolor{cyan}{5} \> \textcolor{cyan}{3} \> \textcolor{cyan}{9} \> \textcolor{cyan}{4} \> \textcolor{cyan}{8} \> \textcolor{cyan}{16} \> \textcolor{cyan}{28} \\
$\textcolor{red}{E_{\mathbb{Q}}(31)}$: \> \textcolor{red}{1} \> \textcolor{red}{2} \> \textcolor{red}{4} \> \textcolor{cyan}{4} \> \textcolor{cyan}{2} \> \textcolor{cyan}{8} \> \textcolor{cyan}{3} \> \textcolor{cyan}{7} \> \textcolor{cyan}{15} \> \textcolor{cyan}{27} \> \textcolor{cyan}{0} \\
$\textcolor{black}{E_{\mathbb{Q}}(32)}$: \> \textcolor{black}{0} \> \textcolor{black}{1} \> \textcolor{black}{3} \> \textcolor{cyan}{3} \> \textcolor{cyan}{1} \> \textcolor{cyan}{7} \> \textcolor{cyan}{2} \> \textcolor{cyan}{6} \> \textcolor{cyan}{14} \> \textcolor{cyan}{26} \> \textcolor{cyan}{30} \\
$\textcolor{black}{E_{\mathbb{Q}}(33)}$: \> \textcolor{black}{1} \> \textcolor{black}{0} \> \textcolor{black}{2} \> \textcolor{cyan}{2} \> \textcolor{cyan}{0} \> \textcolor{cyan}{6} \> \textcolor{cyan}{1} \> \textcolor{cyan}{5} \> \textcolor{cyan}{13} \> \textcolor{cyan}{25} \> \textcolor{cyan}{29} \\
$\textcolor{black}{E_{\mathbb{Q}}(34)}$: \> \textcolor{black}{0} \> \textcolor{black}{2} \> \textcolor{black}{1} \> \textcolor{cyan}{1} \> \textcolor{cyan}{10} \> \textcolor{cyan}{5} \> \textcolor{cyan}{0} \> \textcolor{cyan}{4} \> \textcolor{cyan}{12} \> \textcolor{cyan}{24} \> \textcolor{cyan}{28} \\
$\textcolor{black}{E_{\mathbb{Q}}(35)}$: \> \textcolor{black}{1} \> \textcolor{black}{1} \> \textcolor{black}{0} \> \textcolor{cyan}{0} \> \textcolor{cyan}{9} \> \textcolor{cyan}{4} \> \textcolor{cyan}{16} \> \textcolor{cyan}{3} \> \textcolor{cyan}{11} \> \textcolor{cyan}{23} \> \textcolor{cyan}{27} \\
$\textcolor{black}{E_{\mathbb{Q}}(36)}$: \> \textcolor{black}{0} \> \textcolor{black}{0} \> \textcolor{black}{4} \> \textcolor{cyan}{6} \> \textcolor{cyan}{8} \> \textcolor{cyan}{3} \> \textcolor{cyan}{15} \> \textcolor{cyan}{2} \> \textcolor{cyan}{10} \> \textcolor{cyan}{22} \> \textcolor{cyan}{26} \\
$\textcolor{red}{E_{\mathbb{Q}}(37)}$: \> \textcolor{red}{1} \> \textcolor{red}{2} \> \textcolor{red}{3} \> \textcolor{cyan}{5} \> \textcolor{cyan}{7} \> \textcolor{cyan}{2} \> \textcolor{cyan}{14} \> \textcolor{cyan}{1} \> \textcolor{cyan}{9} \> \textcolor{cyan}{21} \> \textcolor{cyan}{25} \\
$\textcolor{black}{E_{\mathbb{Q}}(38)}$: \> \textcolor{black}{0} \> \textcolor{black}{1} \> \textcolor{black}{2} \> \textcolor{cyan}{4} \> \textcolor{cyan}{6} \> \textcolor{cyan}{1} \> \textcolor{cyan}{13} \> \textcolor{cyan}{0} \> \textcolor{cyan}{8} \> \textcolor{cyan}{20} \> \textcolor{cyan}{24} \\
$\textcolor{black}{E_{\mathbb{Q}}(39)}$: \> \textcolor{black}{1} \> \textcolor{black}{0} \> \textcolor{black}{1} \> \textcolor{cyan}{3} \> \textcolor{cyan}{5} \>\textcolor{cyan}{0} \> \textcolor{cyan}{12} \> \textcolor{cyan}{18} \> \textcolor{cyan}{7} \> \textcolor{cyan}{19} \> \textcolor{cyan}{23} \\
$\textcolor{black}{E_{\mathbb{Q}}(40)}$: \> \textcolor{black}{0} \> \textcolor{black}{2} \> \textcolor{black}{0} \> \textcolor{cyan}{2} \> \textcolor{cyan}{4} \> \textcolor{cyan}{12} \> \textcolor{cyan}{11} \> \textcolor{cyan}{17} \> \textcolor{cyan}{6} \> \textcolor{cyan}{18} \> \textcolor{cyan}{22} \\
$\textcolor{red}{E_{\mathbb{Q}}(41)}$: \> \textcolor{red}{1} \> \textcolor{red}{1} \> \textcolor{red}{4} \> \textcolor{cyan}{1} \> \textcolor{cyan}{3} \> \textcolor{cyan}{11} \> \textcolor{cyan}{10} \> \textcolor{cyan}{16} \> \textcolor{cyan}{5} \> \textcolor{cyan}{17} \> \textcolor{cyan}{21} \\
$\textcolor{black}{E_{\mathbb{Q}}(42)}$: \> \textcolor{black}{0} \> \textcolor{black}{0} \> \textcolor{black}{3} \> \textcolor{cyan}{0} \> \textcolor{cyan}{2} \> \textcolor{cyan}{10} \> \textcolor{cyan}{9} \> \textcolor{cyan}{15} \> \textcolor{cyan}{4} \> \textcolor{cyan}{16} \> \textcolor{cyan}{20} \\
$\textcolor{red}{E_{\mathbb{Q}}(43)}$: \> \textcolor{red}{1} \> \textcolor{red}{2} \> \textcolor{red}{2} \> \textcolor{cyan}{6} \> \textcolor{cyan}{1} \> \textcolor{cyan}{9} \> \textcolor{cyan}{8} \> \textcolor{cyan}{14} \> \textcolor{cyan}{3} \> \textcolor{cyan}{15} \> \textcolor{cyan}{19} \\
$\textcolor{black}{E_{\mathbb{Q}}(44)}$: \> \textcolor{black}{0} \> \textcolor{black}{1} \> \textcolor{black}{1} \> \textcolor{cyan}{5} \> \textcolor{cyan}{0} \> \textcolor{cyan}{8} \> \textcolor{cyan}{7} \> \textcolor{cyan}{13} \> \textcolor{cyan}{2} \> \textcolor{cyan}{14} \> \textcolor{cyan}{18} \\
$\textcolor{black}{E_{\mathbb{Q}}(45)}$: \> \textcolor{black}{1} \> \textcolor{black}{0} \> \textcolor{black}{0} \> \textcolor{cyan}{4} \> \textcolor{cyan}{10} \> \textcolor{cyan}{7} \> \textcolor{cyan}{6} \> \textcolor{cyan}{12} \> \textcolor{cyan}{1} \> \textcolor{cyan}{13} \> \textcolor{cyan}{17} \\
$\textcolor{black}{E_{\mathbb{Q}}(46)}$: \> \textcolor{black}{0} \> \textcolor{black}{2} \> \textcolor{black}{4} \> \textcolor{cyan}{3} \> \textcolor{cyan}{9} \> \textcolor{cyan}{6} \> \textcolor{cyan}{5} \> \textcolor{cyan}{11} \> \textcolor{cyan}{0} \> \textcolor{cyan}{12} \> \textcolor{cyan}{16} \\
$\textcolor{red}{E_{\mathbb{Q}}(47)}$: \> \textcolor{red}{1} \> \textcolor{red}{1} \> \textcolor{red}{3} \> \textcolor{cyan}{2} \> \textcolor{cyan}{8} \> \textcolor{cyan}{5} \> \textcolor{cyan}{4} \> \textcolor{cyan}{10} \> \textcolor{cyan}{22} \> \textcolor{cyan}{11} \> \textcolor{cyan}{15} \\
$\textcolor{black}{E_{\mathbb{Q}}(48)}$: \> \textcolor{black}{0} \> \textcolor{black}{0} \> \textcolor{black}{2} \> \textcolor{cyan}{1} \> \textcolor{cyan}{7} \> \textcolor{cyan}{4} \> \textcolor{cyan}{3} \> \textcolor{cyan}{9} \> \textcolor{cyan}{21} \> \textcolor{cyan}{10} \> \textcolor{cyan}{14} \\
$\textcolor{black}{E_{\mathbb{Q}}(49)}$: \> \textcolor{black}{1} \> \textcolor{black}{2} \> \textcolor{black}{1} \> \textcolor{black}{0} \> \textcolor{cyan}{6} \> \textcolor{cyan}{3} \> \textcolor{cyan}{2} \> \textcolor{cyan}{8} \> \textcolor{cyan}{20} \> \textcolor{cyan}{9} \> \textcolor{cyan}{13} \\
$\textcolor{black}{E_{\mathbb{Q}}(50)}$: \> \textcolor{black}{0} \> \textcolor{black}{1} \> \textcolor{black}{0} \> \textcolor{black}{6} \> \textcolor{cyan}{5} \> \textcolor{cyan}{2} \> \textcolor{cyan}{1} \> \textcolor{cyan}{7} \> \textcolor{cyan}{19} \> \textcolor{cyan}{8} \> \textcolor{cyan}{12} \\
$\textcolor{black}{E_{\mathbb{Q}}(51)}$: \> \textcolor{black}{1} \> \textcolor{black}{0} \> \textcolor{black}{4} \> \textcolor{black}{5} \> \textcolor{cyan}{4} \> \textcolor{cyan}{1} \> \textcolor{cyan}{0} \> \textcolor{cyan}{6} \> \textcolor{cyan}{18} \> \textcolor{cyan}{7} \> \textcolor{cyan}{11} \\
$\textcolor{black}{E_{\mathbb{Q}}(52)}$: \> \textcolor{black}{0} \> \textcolor{black}{2} \> \textcolor{black}{3} \> \textcolor{black}{4} \> \textcolor{cyan}{3} \> \textcolor{cyan}{0} \> \textcolor{cyan}{16} \> \textcolor{cyan}{5} \> \textcolor{cyan}{17} \> \textcolor{cyan}{6} \> \textcolor{cyan}{10} \\
$\textcolor{red}{E_{\mathbb{Q}}(53)}$: \> \textcolor{red}{1} \> \textcolor{red}{1} \> \textcolor{red}{2} \> \textcolor{red}{3} \> \textcolor{cyan}{2} \> \textcolor{cyan}{12} \> \textcolor{cyan}{15} \> \textcolor{cyan}{4} \> \textcolor{cyan}{16} \> \textcolor{cyan}{5} \> \textcolor{cyan}{9} \\
$\textcolor{black}{E_{\mathbb{Q}}(54)}$: \> \textcolor{black}{0} \> \textcolor{black}{0} \> \textcolor{black}{1} \> \textcolor{black}{2} \> \textcolor{cyan}{1} \> \textcolor{cyan}{11} \> \textcolor{cyan}{14} \> \textcolor{cyan}{3} \> \textcolor{cyan}{15} \> \textcolor{cyan}{4} \> \textcolor{cyan}{8} \\
$\textcolor{black}{E_{\mathbb{Q}}(55)}$: \> \textcolor{black}{1} \> \textcolor{black}{2} \> \textcolor{black}{0} \> \textcolor{black}{1} \> \textcolor{cyan}{0} \> \textcolor{cyan}{10} \> \textcolor{cyan}{13} \> \textcolor{cyan}{2} \> \textcolor{cyan}{14} \> \textcolor{cyan}{3} \> \textcolor{cyan}{7} \\
$\textcolor{black}{E_{\mathbb{Q}}(56)}$: \> \textcolor{black}{0} \> \textcolor{black}{1} \> \textcolor{black}{4} \> \textcolor{black}{0} \> \textcolor{cyan}{10} \> \textcolor{cyan}{9} \> \textcolor{cyan}{12} \> \textcolor{cyan}{1} \> \textcolor{cyan}{13} \> \textcolor{cyan}{2} \> \textcolor{cyan}{6} \\
$\textcolor{black}{E_{\mathbb{Q}}(57)}$: \> \textcolor{black}{1} \> \textcolor{black}{0} \> \textcolor{black}{3} \> \textcolor{black}{6} \> \textcolor{cyan}{9} \> \textcolor{cyan}{8} \> \textcolor{cyan}{11} \> \textcolor{cyan}{0} \> \textcolor{cyan}{12} \> \textcolor{cyan}{1} \> \textcolor{cyan}{5} \\
$\textcolor{black}{E_{\mathbb{Q}}(58)}$: \> \textcolor{black}{0} \> \textcolor{black}{2} \> \textcolor{black}{2} \> \textcolor{black}{5} \> \textcolor{cyan}{8} \> \textcolor{cyan}{7} \> \textcolor{cyan}{10} \> \textcolor{cyan}{18} \> \textcolor{cyan}{11} \> \textcolor{cyan}{0} \> \textcolor{cyan}{4} \\
$\textcolor{red}{E_{\mathbb{Q}}(59)}$: \> \textcolor{red}{1} \> \textcolor{red}{1} \> \textcolor{red}{1} \> \textcolor{red}{4} \> \textcolor{cyan}{7} \> \textcolor{cyan}{6} \> \textcolor{cyan}{9} \> \textcolor{cyan}{17} \> \textcolor{cyan}{10} \> \textcolor{cyan}{28} \> \textcolor{cyan}{3} \\
$\textcolor{black}{E_{\mathbb{Q}}(60)}$: \> \textcolor{black}{0} \> \textcolor{black}{0} \> \textcolor{black}{0} \> \textcolor{black}{3} \> \textcolor{cyan}{6} \> \textcolor{cyan}{5} \> \textcolor{cyan}{8} \> \textcolor{cyan}{16} \> \textcolor{cyan}{9} \> \textcolor{cyan}{27} \> \textcolor{cyan}{2} \\
$\textcolor{red}{E_{\mathbb{Q}}(61)}$: \> \textcolor{red}{1} \> \textcolor{red}{2} \> \textcolor{red}{4} \> \textcolor{red}{2} \> \textcolor{cyan}{5} \> \textcolor{cyan}{4} \> \textcolor{cyan}{7} \> \textcolor{cyan}{15} \> \textcolor{cyan}{8} \> \textcolor{cyan}{26} \> \textcolor{cyan}{1} \\
$\textcolor{black}{E_{\mathbb{Q}}(62)}$: \> \textcolor{black}{0} \> \textcolor{black}{1} \> \textcolor{black}{3} \> \textcolor{black}{1} \> \textcolor{cyan}{4} \> \textcolor{cyan}{3} \> \textcolor{cyan}{6} \> \textcolor{cyan}{14} \> \textcolor{cyan}{7} \> \textcolor{cyan}{25} \> \textcolor{cyan}{0} \\
$\textcolor{black}{E_{\mathbb{Q}}(63)}$: \> \textcolor{black}{1} \> \textcolor{black}{0} \> \textcolor{black}{2} \> \textcolor{black}{0} \> \textcolor{cyan}{3} \> \textcolor{cyan}{2} \> \textcolor{cyan}{5} \> \textcolor{cyan}{13} \> \textcolor{cyan}{6} \> \textcolor{cyan}{24} \> \textcolor{cyan}{30} \\
$\textcolor{black}{E_{\mathbb{Q}}(64)}$: \> \textcolor{black}{0} \> \textcolor{black}{2} \> \textcolor{black}{1} \> \textcolor{black}{6} \> \textcolor{cyan}{2} \> \textcolor{cyan}{1} \> \textcolor{cyan}{4} \> \textcolor{cyan}{12} \> \textcolor{cyan}{5} \> \textcolor{cyan}{23} \> \textcolor{cyan}{29} \\
$\textcolor{black}{E_{\mathbb{Q}}(65)}$: \> \textcolor{black}{1} \> \textcolor{black}{1} \> \textcolor{black}{0} \> \textcolor{black}{5} \> \textcolor{cyan}{1} \> \textcolor{cyan}{0} \> \textcolor{cyan}{3} \> \textcolor{cyan}{11} \> \textcolor{cyan}{4} \> \textcolor{cyan}{22} \> \textcolor{cyan}{28} \\
$\textcolor{black}{E_{\mathbb{Q}}(66)}$: \> \textcolor{black}{0} \> \textcolor{black}{0} \> \textcolor{black}{4} \> \textcolor{black}{4} \> \textcolor{cyan}{0} \> \textcolor{cyan}{12} \> \textcolor{cyan}{2} \> \textcolor{cyan}{10} \> \textcolor{cyan}{3} \> \textcolor{cyan}{21} \> \textcolor{cyan}{27} \\
$\textcolor{red}{E_{\mathbb{Q}}(67)}$: \> \textcolor{red}{1} \> \textcolor{red}{2} \> \textcolor{red}{3} \> \textcolor{red}{3} \> \textcolor{cyan}{10} \> \textcolor{cyan}{11} \> \textcolor{cyan}{1} \> \textcolor{cyan}{9} \> \textcolor{cyan}{2} \> \textcolor{cyan}{20} \> \textcolor{cyan}{26} \\
$\textcolor{black}{E_{\mathbb{Q}}(68)}$: \> \textcolor{black}{0} \> \textcolor{black}{1} \> \textcolor{black}{2} \> \textcolor{black}{2} \> \textcolor{cyan}{9} \> \textcolor{cyan}{10} \> \textcolor{cyan}{0} \> \textcolor{cyan}{8} \> \textcolor{cyan}{1} \> \textcolor{cyan}{19} \> \textcolor{cyan}{25} \\
$\textcolor{black}{E_{\mathbb{Q}}(69)}$: \> \textcolor{black}{1} \> \textcolor{black}{0} \> \textcolor{black}{1} \> \textcolor{black}{1} \> \textcolor{cyan}{8} \> \textcolor{cyan}{9} \> \textcolor{cyan}{16} \> \textcolor{cyan}{7} \> \textcolor{cyan}{0} \> \textcolor{cyan}{18} \> \textcolor{cyan}{24} \\
$\textcolor{black}{E_{\mathbb{Q}}(70)}$: \> \textcolor{black}{0} \> \textcolor{black}{2} \> \textcolor{black}{0} \> \textcolor{black}{0} \> \textcolor{cyan}{7} \> \textcolor{cyan}{8} \> \textcolor{cyan}{15} \> \textcolor{cyan}{6} \> \textcolor{cyan}{22} \> \textcolor{cyan}{17} \> \textcolor{cyan}{23} \\
$\textcolor{red}{E_{\mathbb{Q}}(71)}$: \> \textcolor{red}{1} \> \textcolor{red}{1} \> \textcolor{red}{4} \> \textcolor{red}{6} \> \textcolor{cyan}{6} \> \textcolor{cyan}{7} \> \textcolor{cyan}{14} \> \textcolor{cyan}{5} \> \textcolor{cyan}{21} \> \textcolor{cyan}{16} \> \textcolor{cyan}{22} \\
$\textcolor{black}{E_{\mathbb{Q}}(72)}$: \> \textcolor{black}{0} \> \textcolor{black}{0} \> \textcolor{black}{3} \> \textcolor{black}{5} \> \textcolor{cyan}{5} \> \textcolor{cyan}{6} \> \textcolor{cyan}{13} \> \textcolor{cyan}{4} \> \textcolor{cyan}{20} \> \textcolor{cyan}{15} \> \textcolor{cyan}{21} \\
$\textcolor{red}{E_{\mathbb{Q}}(73)}$: \> \textcolor{red}{1} \> \textcolor{red}{2} \> \textcolor{red}{2} \> \textcolor{red}{4} \> \textcolor{cyan}{4} \> \textcolor{cyan}{5} \> \textcolor{cyan}{12} \> \textcolor{cyan}{3} \> \textcolor{cyan}{19} \> \textcolor{cyan}{14} \> \textcolor{cyan}{20} \\
$\textcolor{black}{E_{\mathbb{Q}}(74)}$: \> \textcolor{black}{0} \> \textcolor{black}{1} \> \textcolor{black}{1} \> \textcolor{black}{3} \> \textcolor{cyan}{3} \> \textcolor{cyan}{4} \> \textcolor{cyan}{11} \> \textcolor{cyan}{2} \> \textcolor{cyan}{18} \> \textcolor{cyan}{13} \> \textcolor{cyan}{19} \\
$\textcolor{black}{E_{\mathbb{Q}}(75)}$: \> \textcolor{black}{1} \> \textcolor{black}{0} \> \textcolor{black}{0} \> \textcolor{black}{2} \> \textcolor{cyan}{2} \> \textcolor{cyan}{3} \> \textcolor{cyan}{10} \> \textcolor{cyan}{1} \> \textcolor{cyan}{17} \> \textcolor{cyan}{12} \> \textcolor{cyan}{18} \\
$\textcolor{black}{E_{\mathbb{Q}}(76)}$: \> \textcolor{black}{0} \> \textcolor{black}{2} \> \textcolor{black}{4} \> \textcolor{black}{1} \> \textcolor{cyan}{1} \> \textcolor{cyan}{2} \> \textcolor{cyan}{9} \> \textcolor{cyan}{0} \> \textcolor{cyan}{16} \> \textcolor{cyan}{11} \> \textcolor{cyan}{17} \\
$\textcolor{black}{E_{\mathbb{Q}}(77)}$: \> \textcolor{black}{1} \> \textcolor{black}{1} \> \textcolor{black}{3} \> \textcolor{black}{0} \> \textcolor{cyan}{0} \> \textcolor{cyan}{1} \> \textcolor{cyan}{8} \> \textcolor{cyan}{18} \> \textcolor{cyan}{15} \> \textcolor{cyan}{10} \> \textcolor{cyan}{16} \\
$\textcolor{black}{E_{\mathbb{Q}}(78)}$: \> \textcolor{black}{0} \> \textcolor{black}{0} \> \textcolor{black}{2} \> \textcolor{black}{6} \> \textcolor{cyan}{10} \> \textcolor{cyan}{0} \> \textcolor{cyan}{7} \> \textcolor{cyan}{17} \> \textcolor{cyan}{14} \> \textcolor{cyan}{9} \> \textcolor{cyan}{15} \\
$\textcolor{red}{E_{\mathbb{Q}}(79)}$: \> \textcolor{red}{1} \> \textcolor{red}{2} \> \textcolor{red}{1} \> \textcolor{red}{5} \> \textcolor{cyan}{9} \> \textcolor{cyan}{12} \> \textcolor{cyan}{6} \> \textcolor{cyan}{16} \> \textcolor{cyan}{13} \> \textcolor{cyan}{8} \> \textcolor{cyan}{14} \\
$\textcolor{black}{E_{\mathbb{Q}}(80)}$: \> \textcolor{black}{0} \> \textcolor{black}{1} \> \textcolor{black}{0} \> \textcolor{black}{4} \> \textcolor{cyan}{8} \> \textcolor{cyan}{11} \> \textcolor{cyan}{5} \> \textcolor{cyan}{15} \> \textcolor{cyan}{12} \> \textcolor{cyan}{7} \> \textcolor{cyan}{13} \\
$\textcolor{black}{E_{\mathbb{Q}}(81)}$: \> \textcolor{black}{1} \> \textcolor{black}{0} \> \textcolor{black}{4} \> \textcolor{black}{3} \> \textcolor{cyan}{7} \> \textcolor{cyan}{10} \> \textcolor{cyan}{4} \> \textcolor{cyan}{14} \> \textcolor{cyan}{11} \> \textcolor{cyan}{6} \> \textcolor{cyan}{12} \\
$\textcolor{black}{E_{\mathbb{Q}}(82)}$: \> \textcolor{black}{0} \> \textcolor{black}{2} \> \textcolor{black}{3} \> \textcolor{black}{2} \> \textcolor{cyan}{6} \> \textcolor{cyan}{9} \> \textcolor{cyan}{3} \> \textcolor{cyan}{13} \> \textcolor{cyan}{10} \> \textcolor{cyan}{5} \> \textcolor{cyan}{11} \\
$\textcolor{red}{E_{\mathbb{Q}}(83)}$: \> \textcolor{red}{1} \> \textcolor{red}{1} \> \textcolor{red}{2} \> \textcolor{red}{1} \> \textcolor{cyan}{5} \> \textcolor{cyan}{8} \> \textcolor{cyan}{2} \> \textcolor{cyan}{12} \> \textcolor{cyan}{9} \> \textcolor{cyan}{4} \> \textcolor{cyan}{10} \\
$\textcolor{black}{E_{\mathbb{Q}}(84)}$: \> \textcolor{black}{0} \> \textcolor{black}{0} \> \textcolor{black}{1} \> \textcolor{black}{0} \> \textcolor{cyan}{4} \> \textcolor{cyan}{7} \> \textcolor{cyan}{1} \> \textcolor{cyan}{11} \> \textcolor{cyan}{8} \> \textcolor{cyan}{3} \> \textcolor{cyan}{9} \\
$\textcolor{black}{E_{\mathbb{Q}}(85)}$: \> \textcolor{black}{1} \> \textcolor{black}{2} \> \textcolor{black}{0} \> \textcolor{black}{6} \> \textcolor{cyan}{3} \> \textcolor{cyan}{6} \> \textcolor{cyan}{0} \> \textcolor{cyan}{10} \> \textcolor{cyan}{7} \> \textcolor{cyan}{2} \> \textcolor{cyan}{8} \\
$\textcolor{black}{E_{\mathbb{Q}}(86)}$: \> \textcolor{black}{0} \> \textcolor{black}{1} \> \textcolor{black}{4} \> \textcolor{black}{5} \> \textcolor{cyan}{2} \> \textcolor{cyan}{5} \> \textcolor{cyan}{16} \> \textcolor{cyan}{9} \> \textcolor{cyan}{6} \> \textcolor{cyan}{1} \> \textcolor{cyan}{7} \\
$\textcolor{black}{E_{\mathbb{Q}}(87)}$: \> \textcolor{black}{1} \> \textcolor{black}{0} \> \textcolor{black}{3} \> \textcolor{black}{4} \> \textcolor{cyan}{1} \> \textcolor{cyan}{4} \> \textcolor{cyan}{15} \> \textcolor{cyan}{8} \> \textcolor{cyan}{5} \> \textcolor{cyan}{0} \> \textcolor{cyan}{6} \\
$\textcolor{black}{E_{\mathbb{Q}}(88)}$: \> \textcolor{black}{0} \> \textcolor{black}{2} \> \textcolor{black}{2} \> \textcolor{black}{3} \> \textcolor{cyan}{0} \> \textcolor{cyan}{3} \> \textcolor{cyan}{14} \> \textcolor{cyan}{7} \> \textcolor{cyan}{4} \> \textcolor{cyan}{28} \> \textcolor{cyan}{5} \\
$\textcolor{red}{E_{\mathbb{Q}}(89)}$: \> \textcolor{red}{1} \> \textcolor{red}{1} \> \textcolor{red}{1} \> \textcolor{red}{2} \> \textcolor{cyan}{10} \> \textcolor{cyan}{2} \> \textcolor{cyan}{13} \> \textcolor{cyan}{6} \> \textcolor{cyan}{3} \> \textcolor{cyan}{27} \> \textcolor{cyan}{4} \\
$\textcolor{black}{E_{\mathbb{Q}}(90)}$: \> \textcolor{black}{0} \> \textcolor{black}{0} \> \textcolor{black}{0} \> \textcolor{black}{1} \> \textcolor{cyan}{9} \> \textcolor{cyan}{1} \> \textcolor{cyan}{12} \> \textcolor{cyan}{5} \> \textcolor{cyan}{2} \> \textcolor{cyan}{26} \> \textcolor{cyan}{3} \\
$\textcolor{black}{E_{\mathbb{Q}}(91)}$: \> \textcolor{black}{1} \> \textcolor{black}{2} \> \textcolor{black}{4} \> \textcolor{black}{0} \> \textcolor{cyan}{8} \> \textcolor{cyan}{0} \> \textcolor{cyan}{11} \> \textcolor{cyan}{4} \> \textcolor{cyan}{1} \> \textcolor{cyan}{25} \> \textcolor{cyan}{2} \\
$\textcolor{black}{E_{\mathbb{Q}}(92)}$: \> \textcolor{black}{0} \> \textcolor{black}{1} \> \textcolor{black}{3} \> \textcolor{black}{6} \> \textcolor{cyan}{7} \> \textcolor{cyan}{12} \> \textcolor{cyan}{10} \> \textcolor{cyan}{3} \> \textcolor{cyan}{0} \> \textcolor{cyan}{24} \> \textcolor{cyan}{1} \\
$\textcolor{black}{E_{\mathbb{Q}}(93)}$: \> \textcolor{black}{1} \> \textcolor{black}{0} \> \textcolor{black}{2} \> \textcolor{black}{5} \> \textcolor{cyan}{6} \> \textcolor{cyan}{11} \> \textcolor{cyan}{9} \> \textcolor{cyan}{2} \> \textcolor{cyan}{22} \> \textcolor{cyan}{23} \> \textcolor{cyan}{0} \\
$\textcolor{black}{E_{\mathbb{Q}}(94)}$: \> \textcolor{black}{0} \> \textcolor{black}{2} \> \textcolor{black}{1} \> \textcolor{black}{4} \> \textcolor{cyan}{5} \> \textcolor{cyan}{10} \> \textcolor{cyan}{8} \> \textcolor{cyan}{1} \> \textcolor{cyan}{21} \> \textcolor{cyan}{22} \> \textcolor{cyan}{30} \\
$\textcolor{black}{E_{\mathbb{Q}}(95)}$: \> \textcolor{black}{1} \> \textcolor{black}{1} \> \textcolor{black}{0} \> \textcolor{black}{3} \> \textcolor{cyan}{4} \> \textcolor{cyan}{9} \> \textcolor{cyan}{7} \> \textcolor{cyan}{0} \> \textcolor{cyan}{20} \> \textcolor{cyan}{21} \> \textcolor{cyan}{29} \\
$\textcolor{black}{E_{\mathbb{Q}}(96)}$: \> \textcolor{black}{0} \> \textcolor{black}{0} \> \textcolor{black}{4} \> \textcolor{black}{2} \> \textcolor{cyan}{3} \> \textcolor{cyan}{8} \> \textcolor{cyan}{6} \> \textcolor{cyan}{18} \> \textcolor{cyan}{19} \> \textcolor{cyan}{20} \> \textcolor{cyan}{28} \\
$\textcolor{red}{E_{\mathbb{Q}}(97)}$: \> \textcolor{red}{1} \> \textcolor{red}{2} \> \textcolor{red}{3} \> \textcolor{red}{1} \> \textcolor{cyan}{2} \> \textcolor{cyan}{7} \> \textcolor{cyan}{5} \> \textcolor{cyan}{17} \> \textcolor{cyan}{18} \> \textcolor{cyan}{19} \> \textcolor{cyan}{27} \\
$\textcolor{black}{E_{\mathbb{Q}}(98)}$: \> \textcolor{black}{0} \> \textcolor{black}{1} \> \textcolor{black}{2} \> \textcolor{black}{0} \> \textcolor{cyan}{1} \> \textcolor{cyan}{6} \> \textcolor{cyan}{4} \> \textcolor{cyan}{16} \> \textcolor{cyan}{17} \> \textcolor{cyan}{18} \> \textcolor{cyan}{26} \\
$\textcolor{black}{E_{\mathbb{Q}}(99)}$: \> \textcolor{black}{1} \> \textcolor{black}{0} \> \textcolor{black}{1} \> \textcolor{black}{6} \> \textcolor{cyan}{0} \> \textcolor{cyan}{5} \> \textcolor{cyan}{3} \> \textcolor{cyan}{15} \> \textcolor{cyan}{16} \> \textcolor{cyan}{17} \> \textcolor{cyan}{25} \\
$\textcolor{black}{E_{\mathbb{Q}}(100)}$: \> \textcolor{black}{0} \> \textcolor{black}{2} \> \textcolor{black}{0} \> \textcolor{black}{5} \> \textcolor{cyan}{10} \> \textcolor{cyan}{4} \> \textcolor{cyan}{2} \> \textcolor{cyan}{14} \> \textcolor{cyan}{15} \> \textcolor{cyan}{16} \> \textcolor{cyan}{24} \\
$\textcolor{red}{E_{\mathbb{Q}}(101)}$: \> \textcolor{red}{1} \> \textcolor{red}{1} \> \textcolor{red}{4} \> \textcolor{red}{4} \> \textcolor{cyan}{9} \> \textcolor{cyan}{3} \> \textcolor{cyan}{1} \> \textcolor{cyan}{13} \> \textcolor{cyan}{14} \> \textcolor{cyan}{15} \> \textcolor{cyan}{23} \\
$\textcolor{black}{E_{\mathbb{Q}}(102)}$: \> \textcolor{black}{0} \> \textcolor{black}{0} \> \textcolor{black}{3} \> \textcolor{black}{3} \> \textcolor{cyan}{8} \> \textcolor{cyan}{2} \> \textcolor{cyan}{0} \> \textcolor{cyan}{12} \> \textcolor{cyan}{13} \> \textcolor{cyan}{14} \> \textcolor{cyan}{22} \\
$\textcolor{red}{E_{\mathbb{Q}}(103)}$: \> \textcolor{red}{1} \> \textcolor{red}{2} \> \textcolor{red}{2} \> \textcolor{red}{2} \> \textcolor{cyan}{7} \> \textcolor{cyan}{1} \> \textcolor{cyan}{16} \> \textcolor{cyan}{11} \> \textcolor{cyan}{12} \> \textcolor{cyan}{13} \> \textcolor{cyan}{21} \\
$\textcolor{black}{E_{\mathbb{Q}}(104)}$: \> \textcolor{black}{0} \> \textcolor{black}{1} \> \textcolor{black}{1} \> \textcolor{black}{1} \> \textcolor{cyan}{6} \> \textcolor{cyan}{0} \> \textcolor{cyan}{15} \> \textcolor{cyan}{10} \> \textcolor{cyan}{11} \> \textcolor{cyan}{12} \> \textcolor{cyan}{20} \\
$\textcolor{black}{E_{\mathbb{Q}}(105)}$: \> \textcolor{black}{1} \> \textcolor{black}{0} \> \textcolor{black}{0} \> \textcolor{black}{0} \> \textcolor{cyan}{5} \> \textcolor{cyan}{12} \> \textcolor{cyan}{14} \> \textcolor{cyan}{9} \> \textcolor{cyan}{10} \> \textcolor{cyan}{11} \> \textcolor{cyan}{19} \\
$\textcolor{black}{E_{\mathbb{Q}}(106)}$: \> \textcolor{black}{0} \> \textcolor{black}{2} \> \textcolor{black}{4} \> \textcolor{black}{6} \> \textcolor{cyan}{4} \> \textcolor{cyan}{11} \> \textcolor{cyan}{13} \> \textcolor{cyan}{8} \> \textcolor{cyan}{9} \> \textcolor{cyan}{10} \> \textcolor{cyan}{18} \\
$\textcolor{red}{E_{\mathbb{Q}}(107)}$: \> \textcolor{red}{1} \> \textcolor{red}{1} \> \textcolor{red}{3} \> \textcolor{red}{5} \> \textcolor{cyan}{3} \> \textcolor{cyan}{10} \> \textcolor{cyan}{12} \> \textcolor{cyan}{7} \> \textcolor{cyan}{8} \> \textcolor{cyan}{9} \> \textcolor{cyan}{17} \\
$\textcolor{black}{E_{\mathbb{Q}}(108)}$: \> \textcolor{black}{0} \> \textcolor{black}{0} \> \textcolor{black}{2} \> \textcolor{black}{4} \> \textcolor{cyan}{2} \> \textcolor{cyan}{9} \> \textcolor{cyan}{11} \> \textcolor{cyan}{6} \> \textcolor{cyan}{7} \> \textcolor{cyan}{8} \> \textcolor{cyan}{16} \\
$\textcolor{red}{E_{\mathbb{Q}}(109)}$: \> \textcolor{red}{1} \> \textcolor{red}{2} \> \textcolor{red}{1} \> \textcolor{red}{3} \> \textcolor{cyan}{1} \> \textcolor{cyan}{8} \> \textcolor{cyan}{10} \> \textcolor{cyan}{5} \> \textcolor{cyan}{6} \> \textcolor{cyan}{7} \> \textcolor{cyan}{15} \\
$\textcolor{black}{E_{\mathbb{Q}}(110)}$: \> \textcolor{black}{0} \> \textcolor{black}{1} \> \textcolor{black}{0} \> \textcolor{black}{2} \> \textcolor{cyan}{0} \> \textcolor{cyan}{7} \> \textcolor{cyan}{9} \> \textcolor{cyan}{4} \> \textcolor{cyan}{5} \> \textcolor{cyan}{6} \> \textcolor{cyan}{14} \\
$\textcolor{black}{E_{\mathbb{Q}}(111)}$: \> \textcolor{black}{1} \> \textcolor{black}{0} \> \textcolor{black}{4} \> \textcolor{black}{1} \> \textcolor{cyan}{10} \> \textcolor{cyan}{6} \> \textcolor{cyan}{8} \> \textcolor{cyan}{3} \> \textcolor{cyan}{4} \> \textcolor{cyan}{5} \> \textcolor{cyan}{13} \\
$\textcolor{black}{E_{\mathbb{Q}}(112)}$: \> \textcolor{black}{0} \> \textcolor{black}{2} \> \textcolor{black}{3} \> \textcolor{black}{0} \> \textcolor{cyan}{9} \> \textcolor{cyan}{5} \> \textcolor{cyan}{7} \> \textcolor{cyan}{2} \> \textcolor{cyan}{3} \> \textcolor{cyan}{4} \> \textcolor{cyan}{12} \\
$\textcolor{red}{E_{\mathbb{Q}}(113)}$: \> \textcolor{red}{1} \> \textcolor{red}{1} \> \textcolor{red}{2} \> \textcolor{red}{6} \> \textcolor{cyan}{8} \> \textcolor{cyan}{4} \> \textcolor{cyan}{6} \> \textcolor{cyan}{1} \> \textcolor{cyan}{2} \> \textcolor{cyan}{3} \> \textcolor{cyan}{11} \\
$\textcolor{black}{E_{\mathbb{Q}}(114)}$: \> \textcolor{black}{0} \> \textcolor{black}{0} \> \textcolor{black}{1} \> \textcolor{black}{5} \> \textcolor{cyan}{7} \> \textcolor{cyan}{3} \> \textcolor{cyan}{5} \> \textcolor{cyan}{0} \> \textcolor{cyan}{1} \> \textcolor{cyan}{2} \> \textcolor{cyan}{10} \\
$\textcolor{black}{E_{\mathbb{Q}}(115)}$: \> \textcolor{black}{1} \> \textcolor{black}{2} \> \textcolor{black}{0} \> \textcolor{black}{4} \> \textcolor{cyan}{6} \> \textcolor{cyan}{2} \> \textcolor{cyan}{4} \> \textcolor{cyan}{18} \> \textcolor{cyan}{0} \> \textcolor{cyan}{1} \> \textcolor{cyan}{9} \\
$\textcolor{black}{E_{\mathbb{Q}}(116)}$: \> \textcolor{black}{0} \> \textcolor{black}{1} \> \textcolor{black}{4} \> \textcolor{black}{3} \> \textcolor{cyan}{5} \> \textcolor{cyan}{1} \> \textcolor{cyan}{3} \> \textcolor{cyan}{17} \> \textcolor{cyan}{22} \> \textcolor{cyan}{0} \> \textcolor{cyan}{8} \\
$\textcolor{black}{E_{\mathbb{Q}}(117)}$: \> \textcolor{black}{1} \> \textcolor{black}{0} \> \textcolor{black}{3} \> \textcolor{black}{2} \> \textcolor{cyan}{4} \> \textcolor{cyan}{0} \> \textcolor{cyan}{2} \> \textcolor{cyan}{16} \> \textcolor{cyan}{21} \> \textcolor{cyan}{28} \> \textcolor{cyan}{7} \\
$\textcolor{black}{E_{\mathbb{Q}}(118)}$: \> \textcolor{black}{0} \> \textcolor{black}{2} \> \textcolor{black}{2} \> \textcolor{black}{1} \> \textcolor{cyan}{3} \> \textcolor{cyan}{12} \> \textcolor{cyan}{1} \> \textcolor{cyan}{15} \> \textcolor{cyan}{20} \> \textcolor{cyan}{27} \> \textcolor{cyan}{6} \\
$\textcolor{black}{E_{\mathbb{Q}}(119)}$: \> \textcolor{black}{1} \> \textcolor{black}{1} \> \textcolor{black}{1} \> \textcolor{black}{0} \> \textcolor{cyan}{2} \> \textcolor{cyan}{11} \> \textcolor{cyan}{0} \> \textcolor{cyan}{14} \> \textcolor{cyan}{19} \> \textcolor{cyan}{26} \> \textcolor{cyan}{5} \\
$\textcolor{black}{E_{\mathbb{Q}}(120)}$: \> \textcolor{black}{0} \> \textcolor{black}{0} \> \textcolor{black}{0} \> \textcolor{black}{6} \> \textcolor{cyan}{1} \> \textcolor{cyan}{10} \> \textcolor{cyan}{16} \> \textcolor{cyan}{13} \> \textcolor{cyan}{18} \> \textcolor{cyan}{25} \> \textcolor{cyan}{4} \\
$\textcolor{black}{E_{\mathbb{Q}}(121)}$: \> \textcolor{black}{1} \> \textcolor{black}{2} \> \textcolor{black}{4} \> \textcolor{black}{5} \> \textcolor{black}{0} \> \textcolor{cyan}{9} \> \textcolor{cyan}{15} \> \textcolor{cyan}{12} \> \textcolor{cyan}{17} \> \textcolor{cyan}{24} \> \textcolor{cyan}{3} \\
\ldots \\
$\textcolor{black}{E_{\mathbb{Q}}(n)}$: \> $\textcolor{black}{q_{_{1}}}$ \> $\textcolor{black}{q_{_{2}}}$ \> $\textcolor{black}{q_{_{3}}}$ \> $\textcolor{black}{q_{_{4}}}$ \> $\textcolor{black}{q_{_{5}}}$ \> $\textcolor{cyan}{q_{_{6}}}$ \> $\textcolor{cyan}{q_{_{7}}}$ \> $\textcolor{cyan}{q_{_{8}}}$ \> $\textcolor{cyan}{q_{_{9}}}$ \> $\textcolor{cyan}{q_{_{10}}}$ \> $\textcolor{cyan}{q_{_{11}}}$ \> \textcolor{cyan}{\ldots} \\
\ldots \\
\rule{77mm}{.5mm} \\
Divisors: \= 001 \= 002 \= 003 \= 004 \= 005 \= 006 \= 007 \= 008 \= 009 \= 0010 \= 0011 \= \ldots \= 0n-9 \ldots \kill
Prime: \> $\textcolor{black}{p_{_{1}}}$ \> $\textcolor{black}{p_{_{2}}}$ \> $\textcolor{black}{p_{_{3}}}$ \> $\textcolor{black}{p_{_{4}}}$ \> $\textcolor{black}{p_{_{5}}}$ \> $\textcolor{cyan}{p_{_{6}}}$ \> $\textcolor{cyan}{p_{_{7}}}$ \> $\textcolor{cyan}{p_{_{8}}}$ \> $\textcolor{cyan}{p_{_{9}}}$ \> $\textcolor{cyan}{p_{_{10}}}$ \> $\textcolor{cyan}{p_{_{11}}}$ \> \textcolor{cyan}{\ldots} \> $\textcolor{cyan}{p_{_{n}}}$ \textcolor{cyan}{\ldots} \\
Divisor: \> \textcolor{black}{2} \> \textcolor{black}{3} \> \textcolor{black}{5} \> \textcolor{black}{7} \> \textcolor{black}{11} \> \textcolor{cyan}{13} \> \textcolor{cyan}{17} \> \textcolor{cyan}{19} \> \textcolor{cyan}{23} \> \textcolor{cyan}{29} \> \textcolor{cyan}{31} \> \textcolor{cyan}{\ldots} \> $\textcolor{cyan}{p_{n}}$ \textcolor{cyan}{\ldots}
\end{tabbing}

\noindent {\tiny{Fig.11: The output of the prime-number based algorithms $\textcolor{cyan}{E_{\mathbb{P}}}$ and $\textcolor{black}{E_{\mathbb{Q}}}$}}

\section{Appendix III: Comparative values of actual and non-heuristically estimated number of primes $\leq 1500$}
\label{appendix.prim.dist.1500}

\vspace{+2ex}
\noindent \line(1,0){480} \\ \\
\noindent \textbf{\footnotesize{Fig.12: The following table gives comparative values for $\pi(n)$ as approximated non-heuristically by $\pi_{_{L}}(n) = \sum_{j=1}^{n}\prod_{i=1}^{\pi(\sqrt j)}(1-1/p_{_{j}})$, the actual values $\pi(n)$ of the primes less than or equal to $n$, and the values for $\pi(n)$ as estimated non-heuristically by $\pi_{_{H}}(n) = n.\prod_{i=1}^{\pi(\sqrt j)}(1-1/p_{_{j}})$ of $\pi(n)$, for $4 \leq n \leq 1500$.}}\footnote{The downloadable .xlxs source file is accessible \href{http://alixcomsi.com/Primes_less_than_n_1500.xlsx}{here}.} \\
\noindent \line(1,0){480}

\vspace{+5ex}
\tiny{
$
\begin{array}{||r|r||r|r|r|r||r|r|l||}
n & p_{_{n}} & \pi(\sqrt n) & \prod_{j=1}^{n}(1-1/p_{_{j}}) & [\sqrt {n}] & \prod_{j=1}^{\pi(\sqrt n)}(1-1/p_{_{j}}) & \sum_{j=1}^{n}\prod_{i=1}^{\pi(\sqrt j)}(1-1/p_{_{j}}) & \pi(n) &  n.\prod_{i=1}^{\pi(\sqrt j)}(1-1/p_{_{j}}) \\	
 &   &   &   &   &   &   &   & \\		
1	&	2	&	0	&	   1/2  	&	1	&		&	0.0000	&	0	&	0.000	\\
2	&	3	&	0	&	   1/3  	&	1	&	   1/2  	&	0.5000	&	1	&	1.0000	\\
3	&	5	&	0	&	   4/15 	&	1	&	   1/2  	&	1.0000	&	2	&	1.5000	\\
 &   &   &   &   &   &   &   & \\		
4	&	7	&	1	&	   8/35 	&	2	&	   1/3  	&	1.3333	&	2	&	1.3333	\\
5	&	11	&	1	&	  16/77 	&	2	&	   1/3  	&	1.6667	&	3	&	1.6667	\\
\end{array}$
}

\vspace{+2ex}
\tiny{
$\begin{array}{||r|r||r|r|r|r||r|r|l||}
n & p_{_{n}} & \pi(\sqrt n) & \prod_{j=1}^{n}(1-1/p_{_{j}}) & [\sqrt {n}] & \prod_{j=1}^{\pi(\sqrt n)}(1-1/p_{_{j}}) & \sum_{j=1}^{n}\prod_{i=1}^{\pi(\sqrt j)}(1-1/p_{_{j}}) & \pi(n) &  n.\prod_{i=1}^{\pi(\sqrt j)}(1-1/p_{_{j}}) \\		
 &   &   &   &   &   &   &   & \\	
71	&	353	&	4	&	  11/116	&	8	&	  16/77 	&	17.1983	&	20	&	14.7532	\\
72	&	359	&	4	&	  33/349	&	8	&	  16/77 	&	17.4061	&	20	&	14.9610	\\
73	&	367	&	4	&	  43/456	&	8	&	  16/77 	&	17.6139	&	21	&	15.1688	\\
74	&	373	&	4	&	  79/840	&	8	&	  16/77 	&	17.8216	&	21	&	15.3766	\\
75	&	379	&	4	&	  59/629	&	8	&	  16/77 	&	18.0294	&	21	&	15.5844	\\
76	&	383	&	4	&	  45/481	&	8	&	  16/77 	&	18.2372	&	21	&	15.7922	\\
\end{array}$
}

\vspace{+2ex}
\tiny{
$\begin{array}{||r|r||r|r|r|r||r|r|l||}
n & p_{_{n}} & \pi(\sqrt n) & \prod_{j=1}^{n}(1-1/p_{_{j}}) & [\sqrt {n}] & \prod_{j=1}^{\pi(\sqrt n)}(1-1/p_{_{j}}) & \sum_{j=1}^{n}\prod_{i=1}^{\pi(\sqrt j)}(1-1/p_{_{j}}) & \pi(n) &  n.\prod_{i=1}^{\pi(\sqrt j)}(1-1/p_{_{j}}) \\		
 &   &   &   &   &   &   &   & \\		
175	&	1039	&	6	&	  39/485	&	13	&	 165/914	&	37.8506	&	40	&	31.5919	\\
176	&	1049	&	6	&	  67/834	&	13	&	 165/914	&	38.0311	&	40	&	31.7725	\\
177	&	1051	&	6	&	  37/461	&	13	&	 165/914	&	38.2116	&	40	&	31.9530	\\
178	&	1061	&	6	&	  17/212	&	13	&	 165/914	&	38.3921	&	40	&	32.1335	\\
179	&	1063	&	6	&	  29/362	&	13	&	 165/914	&	38.5727	&	41	&	32.3140	\\
180	&	1069	&	6	&	   2/25 	&	13	&	 165/914	&	38.7532	&	41	&	32.4946	\\
\end{array}$
}

\vspace{+2ex}
\tiny{
$\begin{array}{||r|r||r|r|r|r||r|r|l||}
n & p_{_{n}} & \pi(\sqrt n) & \prod_{j=1}^{n}(1-1/p_{_{j}}) & [\sqrt {n}] & \prod_{j=1}^{\pi(\sqrt n)}(1-1/p_{_{j}}) & \sum_{j=1}^{n}\prod_{i=1}^{\pi(\sqrt j)}(1-1/p_{_{j}}) & \pi(n) &  n.\prod_{i=1}^{\pi(\sqrt j)}(1-1/p_{_{j}}) \\		
 &   &   &   &   &   &   &   & \\	
279	&	1801	&	6	&	  47/630	&	16	&	 165/914	&	56.6252	&	59	&	50.3666	\\
280	&	1811	&	6	&	  17/228	&	16	&	 165/914	&	56.8057	&	59	&	50.5471	\\
281	&	1823	&	6	&	  74/993	&	16	&	 165/914	&	56.9862	&	60	&	50.7276	\\
282	&	1831	&	6	&	  61/819	&	16	&	 165/914	&	57.1668	&	60	&	50.9082	\\
283	&	1847	&	6	&	  30/403	&	16	&	 165/914	&	57.3473	&	61	&	51.0887	\\
284	&	1861	&	6	&	  59/793	&	16	&	 165/914	&	57.5278	&	61	&	51.2692	\\
\end{array}$
}

\vspace{+2ex}
\tiny{
$\begin{array}{||r|r||r|r|r|r||r|r|l||}
n & p_{_{n}} & \pi(\sqrt n) & \prod_{j=1}^{n}(1-1/p_{_{j}}) & [\sqrt {n}] & \prod_{j=1}^{\pi(\sqrt n)}(1-1/p_{_{j}}) & \sum_{j=1}^{n}\prod_{i=1}^{\pi(\sqrt j)}(1-1/p_{_{j}}) & \pi(n) &  n.\prod_{i=1}^{\pi(\sqrt j)}(1-1/p_{_{j}}) \\	
 &   &   &   &   &   &   &   & \\	
383	&	2647	&	8	&	  18/253	&	19	&	  62/379	&	74.3262	&	76	&	62.6543	\\
384	&	2657	&	8	&	  50/703	&	19	&	  62/379	&	74.4898	&	76	&	62.8179	\\
385	&	2659	&	8	&	  61/858	&	19	&	  62/379	&	74.6534	&	76	&	62.9815	\\
386	&	2663	&	8	&	  14/197	&	19	&	  62/379	&	74.8169	&	76	&	63.1450	\\
387	&	2671	&	8	&	  13/183	&	19	&	  62/379	&	74.9805	&	76	&	63.3086	\\
388	&	2677	&	8	&	  37/521	&	19	&	  62/379	&	75.1441	&	76	&	63.4722	\\
\end{array}$
}

\vspace{+2ex}
\tiny{
$\begin{array}{||r|r||r|r|r|r||r|r|l||}
n & p_{_{n}} & \pi(\sqrt n) & \prod_{j=1}^{n}(1-1/p_{_{j}}) & [\sqrt {n}] & \prod_{j=1}^{\pi(\sqrt n)}(1-1/p_{_{j}}) & \sum_{j=1}^{n}\prod_{i=1}^{\pi(\sqrt j)}(1-1/p_{_{j}}) & \pi(n) &  n.\prod_{i=1}^{\pi(\sqrt j)}(1-1/p_{_{j}}) \\		
 &   &   &   &   &   &   &   & \\	
487	&	3469	&	8	&	  11/160	&	22	&	  62/379	&	91.3394	&	93	&	79.6675	\\
488	&	3491	&	8	&	  20/291	&	22	&	  62/379	&	91.5029	&	93	&	79.8310	\\
489	&	3499	&	8	&	   9/131	&	22	&	  62/379	&	91.6665	&	93	&	79.9946	\\
490	&	3511	&	8	&	  34/495	&	22	&	  62/379	&	91.8301	&	93	&	80.1582	\\
491	&	3517	&	8	&	  16/233	&	22	&	  62/379	&	91.9937	&	94	&	80.3218	\\
492	&	3527	&	8	&	  67/976	&	22	&	  62/379	&	92.1573	&	94	&	80.4854	\\
\end{array}$
}

\vspace{+2ex}
\tiny{
$\begin{array}{||r|r||r|r|r|r||r|r|l||}
n & p_{_{n}} & \pi(\sqrt n) & \prod_{j=1}^{n}(1-1/p_{_{j}}) & [\sqrt {n}] & \prod_{j=1}^{\pi(\sqrt n)}(1-1/p_{_{j}}) & \sum_{j=1}^{n}\prod_{i=1}^{\pi(\sqrt j)}(1-1/p_{_{j}}) & \pi(n) &  n.\prod_{i=1}^{\pi(\sqrt j)}(1-1/p_{_{j}}) \\		
 &   &   &   &   &   &   &   & \\	
590	&	4297	&	9	&	  65/971	&	24	&	 157/994	&	107.8392	&	107	&	93.1889	\\
591	&	4327	&	9	&	  17/254	&	24	&	 157/994	&	107.9971	&	107	&	93.3468	\\
592	&	4337	&	9	&	  55/822	&	24	&	 157/994	&	108.1551	&	107	&	93.5048	\\
593	&	4339	&	9	&	  39/583	&	24	&	 157/994	&	108.3130	&	108	&	93.6627	\\
594	&	4349	&	9	&	  21/314	&	24	&	 157/994	&	108.4710	&	108	&	93.8207	\\
595	&	4357	&	9	&	  45/673	&	24	&	 157/994	&	108.6289	&	108	&	93.9786	\\
\end{array}$
}

\vspace{+2ex}
\tiny{
$\begin{array}{||r|r||r|r|r|r||r|r|l||}
n & p_{_{n}} & \pi(\sqrt n) & \prod_{j=1}^{n}(1-1/p_{_{j}}) & [\sqrt {n}] & \prod_{j=1}^{\pi(\sqrt n)}(1-1/p_{_{j}}) & \sum_{j=1}^{n}\prod_{i=1}^{\pi(\sqrt j)}(1-1/p_{_{j}}) & \pi(n) &  n.\prod_{i=1}^{\pi(\sqrt j)}(1-1/p_{_{j}}) \\		
 &   &   &   &   &   &   &   & \\
694	&	5227	&	9	&	  41/626	&	26	&	 157/994	&	124.2657	&	125	&	109.6154	\\
695	&	5231	&	9	&	  11/168	&	26	&	 157/994	&	124.4237	&	125	&	109.7733	\\
696	&	5233	&	9	&	  51/779	&	26	&	 157/994	&	124.5816	&	125	&	109.9313	\\
697	&	5237	&	9	&	  18/275	&	26	&	 157/994	&	124.7395	&	125	&	110.0892	\\
698	&	5261	&	9	&	  25/382	&	26	&	 157/994	&	124.8975	&	125	&	110.2472	\\
699	&	5273	&	9	&	  53/810	&	26	&	 157/994	&	125.0554	&	125	&	110.4051	\\
\end{array}$
}

\vspace{+2ex}
\tiny{
$\begin{array}{||r|r||r|r|r|r||r|r|l||}
n & p_{_{n}} & \pi(\sqrt n) & \prod_{j=1}^{n}(1-1/p_{_{j}}) & [\sqrt {n}] & \prod_{j=1}^{\pi(\sqrt n)}(1-1/p_{_{j}}) & \sum_{j=1}^{n}\prod_{i=1}^{\pi(\sqrt j)}(1-1/p_{_{j}}) & \pi(n) &  n.\prod_{i=1}^{\pi(\sqrt j)}(1-1/p_{_{j}}) \\	
 &   &   &   &   &   &   &   & \\
798	&	6121	&	9	&	  29/451	&	28	&	 157/994	&	140.6922	&	139	&	126.0419	\\
799	&	6131	&	9	&	   9/140	&	28	&	 157/994	&	140.8502	&	139	&	126.1998	\\
800	&	6133	&	9	&	   9/140	&	28	&	 157/994	&	141.0081	&	139	&	126.3578	\\
801	&	6143	&	9	&	  59/918	&	28	&	 157/994	&	141.1661	&	139	&	126.5157	\\
802	&	6151	&	9	&	  16/249	&	28	&	 157/994	&	141.3240	&	139	&	126.6737	\\
803	&	6163	&	9	&	  39/607	&	28	&	 157/994	&	141.4820	&	139	&	126.8316	\\
\end{array}$
}

\vspace{+2ex}
\tiny{
$\begin{array}{||r|r||r|r|r|r||r|r|l||}
n & p_{_{n}} & \pi(\sqrt n) & \prod_{j=1}^{n}(1-1/p_{_{j}}) & [\sqrt {n}] & \prod_{j=1}^{\pi(\sqrt n)}(1-1/p_{_{j}}) & \sum_{j=1}^{n}\prod_{i=1}^{\pi(\sqrt j)}(1-1/p_{_{j}}) & \pi(n) &  n.\prod_{i=1}^{\pi(\sqrt j)}(1-1/p_{_{j}}) \\		
 &   &   &   &   &   &   &   & \\
902	&	7013	&	10	&	   5/79 	&	30	&	 142/929	&	156.8028	&	154	&	137.8726	\\
903	&	7019	&	10	&	   5/79 	&	30	&	 142/929	&	156.9557	&	154	&	138.0255	\\
904	&	7027	&	10	&	  46/727	&	30	&	 142/929	&	157.1085	&	154	&	138.1783	\\
905	&	7039	&	10	&	  31/490	&	30	&	 142/929	&	157.2614	&	154	&	138.3312	\\
906	&	7043	&	10	&	  21/332	&	30	&	 142/929	&	157.4142	&	154	&	138.4840	\\
907	&	7057	&	10	&	  53/838	&	30	&	 142/929	&	157.5671	&	155	&	138.6369	\\
\end{array}$
}

\vspace{+2ex}
\tiny{
$\begin{array}{||r|r||r|r|r|r||r|r|l||}
n & p_{_{n}} & \pi(\sqrt n) & \prod_{j=1}^{n}(1-1/p_{_{j}}) & [\sqrt {n}] & \prod_{j=1}^{\pi(\sqrt n)}(1-1/p_{_{j}}) & \sum_{j=1}^{n}\prod_{i=1}^{\pi(\sqrt j)}(1-1/p_{_{j}}) & \pi(n) &  n.\prod_{i=1}^{\pi(\sqrt j)}(1-1/p_{_{j}}) \\	
 &   &   &   &   &   &   &   & \\
1006	&	7963	&	11	&	  48/769	&	31	&	  29/195	&	172.5094	&	169	&	149.6133	\\
1007	&	7993	&	11	&	  44/705	&	31	&	  29/195	&	172.6581	&	169	&	149.7621	\\
1008	&	8009	&	11	&	  41/657	&	31	&	  29/195	&	172.8069	&	169	&	149.9108	\\
1009	&	8011	&	11	&	  38/609	&	31	&	  29/195	&	172.9556	&	169	&	150.0595	\\
1010	&	8017	&	11	&	  35/561	&	31	&	  29/195	&	173.1043	&	169	&	150.2082	\\
1011	&	8039	&	11	&	  33/529	&	31	&	  29/195	&	173.2530	&	169	&	150.3569	\\
\end{array}$
}

\vspace{+2ex}
\tiny{
$\begin{array}{||r|r||r|r|r|r||r|r|l||}
n & p_{_{n}} & \pi(\sqrt n) & \prod_{j=1}^{n}(1-1/p_{_{j}}) & [\sqrt {n}] & \prod_{j=1}^{\pi(\sqrt n)}(1-1/p_{_{j}}) & \sum_{j=1}^{n}\prod_{i=1}^{\pi(\sqrt j)}(1-1/p_{_{j}}) & \pi(n) &  n.\prod_{i=1}^{\pi(\sqrt j)}(1-1/p_{_{j}}) \\	
 &   &   &   &   &   &   &   & \\
1110	&	8929	&	11	&	  32/519	&	33	&	  29/195	&	187.9764	&	186	&	1650803	\\
1111	&	8933	&	11	&	   9/146	&	33	&	  29/195	&	188.1251	&	186	&	165.2290	\\
1112	&	8941	&	11	&	   9/146	&	33	&	  29/195	&	188.2739	&	186	&	165.3778	\\
1113	&	8951	&	11	&	  49/795	&	33	&	  29/195	&	188.4226	&	186	&	165.5265	\\
1114	&	8963	&	11	&	  31/503	&	33	&	  29/195	&	188.5713	&	186	&	165.6752	\\
1115	&	8969	&	11	&	  57/925	&	33	&	  29/195	&	188.7200	&	186	&	165.8239	\\
\end{array}$
}

\vspace{+2ex}
\tiny{
$\begin{array}{||r|r||r|r|r|r||r|r|l||}
n & p_{_{n}} & \pi(\sqrt n) & \prod_{j=1}^{n}(1-1/p_{_{j}}) & [\sqrt {n}] & \prod_{j=1}^{\pi(\sqrt n)}(1-1/p_{_{j}}) & \sum_{j=1}^{n}\prod_{i=1}^{\pi(\sqrt j)}(1-1/p_{_{j}}) & \pi(n) &  n.\prod_{i=1}^{\pi(\sqrt j)}(1-1/p_{_{j}}) \\	
 &   &   &   &   &   &   &   & \\
1214	&	9839	&	11	&	   5/82 	&	34	&	  29/195	&	203.4434	&	198	&	180.5473	\\
1215	&	9851	&	11	&	   5/82 	&	34	&	  29/195	&	203.5921	&	198	&	180.6960	\\
1216	&	9857	&	11	&	   5/82 	&	34	&	  29/195	&	203.7408	&	198	&	180.8448	\\
1217	&	9859	&	11	&	  42/689	&	34	&	  29/195	&	203.8896	&	199	&	180.9935	\\
1218	&	9871	&	11	&	  32/525	&	34	&	  29/195	&	204.0383	&	199	&	181.1422	\\
1219	&	9883	&	11	&	  49/804	&	34	&	  29/195	&	204.1870	&	199	&	181.2909	\\
\end{array}$
}

\vspace{+2ex}
\tiny{
$\begin{array}{||r|r||r|r|r|r||r|r|l||}
n & p_{_{n}} & \pi(\sqrt n) & \prod_{j=1}^{n}(1-1/p_{_{j}}) & [\sqrt {n}] & \prod_{j=1}^{\pi(\sqrt n)}(1-1/p_{_{j}}) & \sum_{j=1}^{n}\prod_{i=1}^{\pi(\sqrt j)}(1-1/p_{_{j}}) & \pi(n) &  n.\prod_{i=1}^{\pi(\sqrt j)}(1-1/p_{_{j}}) \\	
 &   &   &   &   &   &   &   & \\
1318	&	10847	&	11	&	   2/33 	&	36	&	  29/195	&	218.9104	&	214	&	196.0143	\\
1319	&	10853	&	11	&	   2/33 	&	36	&	  29/195	&	219.0591	&	215	&	196.1630	\\
1320	&	10859	&	11	&	   2/33 	&	36	&	  29/195	&	219.2078	&	215	&	196.3117	\\
1321	&	10861	&	11	&	   2/33 	&	36	&	  29/195	&	219.3565	&	216	&	196.4605	\\
1322	&	10867	&	11	&	   2/33 	&	36	&	  29/195	&	219.5053	&	216	&	196.6092	\\
1323	&	10883	&	11	&	   2/33 	&	36	&	  29/195	&	219.6540	&	216	&	196.7579	\\
\end{array}$
}

\vspace{+2ex}
\tiny{
$\begin{array}{||r|r||r|r|r|r||r|r|l||}
n & p_{_{n}} & \pi(\sqrt n) & \prod_{j=1}^{n}(1-1/p_{_{j}}) & [\sqrt {n}] & \prod_{j=1}^{\pi(\sqrt n)}(1-1/p_{_{j}}) & \sum_{j=1}^{n}\prod_{i=1}^{\pi(\sqrt j)}(1-1/p_{_{j}}) & \pi(n) &  n.\prod_{i=1}^{\pi(\sqrt j)}(1-1/p_{_{j}}) \\	
 &   &   &   &   &   &   &   & \\
1422	&	11863	&	12	&	   7/116	&	37	&	 139/958	&	234.1815	&	223	&	206.3232	\\
1423	&	11867	&	12	&	   7/116	&	37	&	 139/958	&	234.3266	&	224	&	206.4683	\\
1424	&	11887	&	12	&	   7/116	&	37	&	 139/958	&	234.4717	&	224	&	206.6134	\\
1425	&	11897	&	12	&	   7/116	&	37	&	 139/958	&	234.6168	&	224	&	206.7585	\\
1426	&	11903	&	12	&	   7/116	&	37	&	 139/958	&	234.7619	&	224	&	206.9036	\\
1427	&	11909	&	12	&	   7/116	&	37	&	 139/958	&	234.9070	&	225	&	207.0487	\\
\end{array}$
}

\vspace{+5ex}
\subsection{Observations and analysis of the error between the actual ($Act\ p$) number $\pi(p_{_{n+1}}^{2}) - \pi(p_{_{n}}^{2})$ and non-heuristically expected ($Exp\ p$) number $\pi_{_{L}}(p_{_{n+1}}^{2}) - \pi_{_{L}}(p_{_{n}}^{2})$ of primes in the interval ($Int$) $(p_{_{n}}^{2}, p_{_{n+1}}^{2})$ for $1 \leq n \leq 11$}
\label{hyp.obs}

\vspace{+2ex}
\footnotesize
Observation and analysis of the error between the actual ($Act\ p$) number, $\pi(p_{_{n+1}}^{2}) - \pi(p_{_{n}}^{2})$, of primes in the interval ($Int$) $(p_{_{n}}^{2}, p_{_{n+1}}^{2})$, and the non-heuristically expected ($Exp\ p$) number, $\pi_{_{L}}(p_{_{n+1}}^{2}) - \pi_{_{L}}(p_{_{n}}^{2})$, of primes in the interval ($Int$) $(p_{_{n}}^{2}, p_{_{n+1}}^{2})$ for $1 \leq n \leq 11$ raises the query:

\vspace{+1ex}
\begin{quote}
Does the ratio:

\begin{quote}
$\mathbb{R} = \frac{CSD}{CExp\ p} = \frac{Cumulative\ standard\ deviation\ of\ the\ cumulative\ sum\ of\ expected\ primes\ in\ the\ interval\ (4,\ p_{_{n+1}}^{2})}{Cumulative\ sum\ \pi_{_{L}}(p_{_{n+1}}^{2})\ =\ \sum_{j=1}^{p_{_{n+1}}^{2}}\prod_{i=1}^{\pi(\sqrt j)}(1-1/p_{_{j}})\ of\ expected\ primes\ in\ the\ interval\ (4,\ p_{_{n+1}}^{2})}$
\end{quote}

tend to a limit?
\end{quote}

\vspace{+1ex}
\noindent \line(1,0){472} \\ \\
\noindent \textbf{\footnotesize{Fig.13: Ratio $CSD/C Exp\ p$ = $\sum_{i=1}^{n}$ Standard Deviation in $(p_{_{i}}^{2}, p_{_{i+1}}^{2})$/$\sum_{i=1}^{n}$ Expected Primes in $(p_{_{i}}^{2}, p_{_{i+1}}^{2})$}} \\
\noindent \line(1,0){472}

\vspace{+5ex}
\tiny{
$\begin{array}{||r|r|r||r|r|r||r|r|r|r||r|r|r||c||}
n & Interval & Int & Int & Int & Int & Cum & Cum & Cum & \% & Int\ density & Int & Cum & Ratio \\	
 & p_{_{n}}^{2} - p_{_{n+1}}^{2}  & Size  & Act\ p  & Exp\ p  & Error  & Act\ p  & Exp\ p  & Error & Error & \prod_{1}^{n}(1 - \frac{1}{p_{_{i}}}) & SD & SD & \mathbb{R} \\
 &   &   &   &   &   &   &   &   &   &   &   &  & \\
1	&	4 - 9	&	5	&	   2 	&	1.6	&	 0.4000	&	2	&	1.6	& 0.4000 & 20.00 & 0.3333 & 1.0541 & 1.0541 & 0.6588 \\
2	&	9 - 25	&	16	&	  5	&	4.2	&	 0.7714	&	7	&	5.8	& 1.1714 & 16.73 & 0.2667 & 1.7689 & 2.8230 & 0.4843 \\
3	&	25 - 49	&	24	&	  6	&	5.5	&	 0.5351	&	13	&	11.3	 & 1.7065 & 13.13 & 0.2286 & 2.0571 & 4.8801 & 0.4321 \\
4	&	49 - 121	&	72	&	  15	&	14.9	&	 0.0549	&	28	&	26.2	 & 1.7614 & 6.29 & 0.2078 & 3.4427 & 8.3228 & 0.3172 \\
5	&	121 - 169	&	48	&	  9	&	9.2	&	 -0.1955	&	37	&	35.4 & 1.5659 & 4.23 & 0.1918 & 2.7278 & 11.0506 & 0.3119 \\
6	&	169 - 289	&	120	&	  22	&	21.7	&	 0.3465	&	59	&	57.1	 & 1.9124 & 3.24 & 0.1805 & 4.2133 & 15.2640 & 0.2674 \\
7	&	289 - 361	&	72	&	  11	&	12.3	&	 -1.3063	&	70	&	69.4	 & 0.6061 & 0.87 & 0.1710 & 3.1950 & 18.4589 & 0.2660 \\
8	&	361 - 529	&	168	&	  29	&	27.5	&	 1.5228	&	99	&	96.9	  & 2.1289 & 2.15 & 0.1636 & 4.7945 & 23.2534 & 0.2400 \\
9	&	529 - 841	&	312	&	  47	&	49.3	&	 -2.2745	&	146	&	146.1 & -0.1456 & -0.10 & 0.1579 & 6.4417 & 29.6951 & 0.2032 \\
10	&	841 - 961	&	120	&	  16	&	18.4	&	 -2.3381	&	162	&	164.5 & -2.4837 & -1.53 & 0.1529 & 3.9419 & 33.6371 & 0.2045 \\
11	&	961 - 1369	&	408	&	  57	&	60.7	&	 -3.6745	&	219	&	225.2 & -6.1582 & -2.81 & 0.1487 & 7.1870 & 40.8241 & 0.1813 \\
12	&	1369 - 1500	&	131	&	  20	&	19.0	&	 0.9927	&	239	&	244.2 & -5.1655 & -2.16 & 0.1451 & 4.0311 & 44.8551 & \textcolor{red}{0.1837} \\
\end{array}$
}

\vspace{+5ex}


\noindent \tiny\textbf{Acknowledgements}: I am indebted to my erstwhile classmate, Professor Chaitanya Kumar Harilan Mehta, for his unqualified encouragement and support for my scholarly pursuits over the years, without which this extension of a 1964 investigation into the nature of divisibility and the structure of the primes---begun whilst yet classmates---would have vanished into some black hole of the informal universe of seemingly self-evident truths. I am also grateful for the extra-ordinary patience and indulgent persistence of Professor William Timothy Gowers in impressing upon me that my original argument conflated the concept of the probability of a number being a prime with that of the density of primes.


\bibliographystyle{amsalpha}

\authoraddresses{
Bhupinder Singh Anand\\
\#1003 B Wing, Lady Ratan Tower\\ Dainik Shivner Marg\\ Gandhinagar, Worli\\ Mumbai - 400 018\\ Maharashtra, India. \\
\email bhup.anand@gmail.com
}

\end{document}